\documentclass{agtart_a}
\pdfoutput=1

\usepackage{pinlabel,bbm}

\title{Surgery untying of coloured knots}

\author{Daniel Moskovich}
\givenname{Daniel}
\surname{Moskovich}
\address{Research Institute for Mathematical Sciences\\
Kyoto University\\\newline
Kyoto 606-8502\\
Japan}
\email{dmoskovich@gmail.com}
\urladdr{http://www.sumamathematica.com/}

\volumenumber{6}
\issuenumber{}
\publicationyear{2006}
\papernumber{26}
\startpage{673}
\endpage{697}

\doi{}
\MR{}
\Zbl{}


\received{25 June 2005}
\revised{}
\accepted{}
\published{24 May 2006}
\publishedonline{24 May 2006}
\proposed{}
\seconded{}
\corresponding{}
\editor{}
\version{}

\arxivreference{math.GT/0506541}

\let\xysavmatrix\xymatrix
\def\xymatrix{\disablesubscriptcorrection\xysavmatrix}
\AtBeginDocument{\let\bar\wbar\let\tilde\wtilde\let\hat\what}

\let\hash\#
\let\mathds\mathbb   

\makeatletter
\def\cnewtheorem#1[#2]#3{\newtheorem{#1}{#3}
\expandafter\let\csname c@#1\endcsname\c@thm}

\newtheorem{thm}{Theorem}
\cnewtheorem{cor}[thm]{Corollary}
\cnewtheorem{lem}[thm]{Lemma}
\newtheorem{conj}{Conjecture}
\newtheorem*{remi}{Reminder}
\cnewtheorem{prop}[thm]{Proposition}
\theoremstyle{definition}
\cnewtheorem{defn}[thm]{Definition}
\theoremstyle{remark}
\newtheorem*{notation}{Notation}

\makeatother  

\newcommand{\abs}[1]{\bigl\vert#1\bigl\vert}
\newcommand{\set}[1]{\left\{#1\right\}}
\newcommand{\To}{\longrightarrow}
\newcommand{\hC}{\hat{C}_{2}}
\newcommand{\seq}{\sim}
\def\ass{\mathrel{\mathop:}=}
\allowdisplaybreaks

\begin{document}

\begin{asciiabstract}
For p=3 and for p=5 we prove that there are exactly p
equivalence classes of p-coloured knots modulo (+/-1)--framed
surgeries along unknots in the kernel of a p-colouring. These
equivalence classes are represented by connect-sums of n left-hand
(p,2)-torus knots with a given colouring when n=1,2,...,p.
This gives a 3-colour and a 5-colour analogue of the
surgery presentation of a knot.
\end{asciiabstract}

\begin{htmlabstract}
For p=3 and for p=5 we prove that there are exactly p equivalence classes
of p&ndash;coloured knots modulo &plusmn;1&ndash;framed surgeries along
unknots in the kernel of a p&ndash;colouring. These equivalence classes
are represented by connect-sums of n left-hand (p,2)&ndash;torus knots
with a given colouring when n=1,2,&hellip;,p.  This gives a 3&ndash;colour
and a 5&ndash;colour analogue of the surgery presentation of a knot.
\end{htmlabstract}

\begin{abstract}
For $p=3$ and for $p=5$ we prove that there are exactly $p$
equivalence classes of $p$--coloured knots modulo $\pm1$--framed
surgeries along unknots in the kernel of a $p$--colouring. These
equivalence classes are represented by connect-sums of $n$ left-hand
$(p,2)$--torus knots with a given colouring when $n=1,2,\ldots,p$.
This gives a $3$--colour and a $5$--colour analogue of the
surgery presentation of a knot.
\end{abstract}

\maketitle

\section{Introduction}

A \emph{$p$--colouring} of a knot $K$ is a surjective homomorphism
$\rho$ from its knot group $\mathfrak{G}:=\pi_{1}(S^{3}-K)$ to
$D_{2p}:=\{t,s|t^{2}=s^{p}=1, tst=s^{p-1}\}$ the dihedral group of
order $2p$, when $p$ is any odd integer. The pair $(K,\rho)$ is
called a \textit{$p$--coloured knot}. It is a well-known fact that
we can encode $\rho$ as a colouring of arcs of a knot diagram by
elements of $\mathds{Z}_{p}$ (the cyclic group of order $p$),
subject to the `colouring rule' that at least two colours are used,
and that at each crossing half the sum of the labels of the
under-crossing arcs equals the label of the over-crossing arc modulo
$p$. By labeling an arc by an element $n\in\mathds{Z}_{p}$, we are
indicating that $\rho$ maps its meridian to the element $ts^{n}\in
D_{2p}$. If a knot $K$ admits such a colouring the knot is said to
be \textit{$p$--colourable}, and $\rho$ is said to be a
\textit{$p$--colouring} of $K$.  These definitions may be extended
to links and to tangles. A necessary and sufficient condition for a
knot to be $p$--colourable is that its determinant be divisible by
$p$. Whether or not a knot is $p$--colourable is the simplest
non-trivial invariant which detects non-commutativity of the knot
group. For more about $p$--colourability we refer the reader to
Fox's original paper \cite{Fox62}.

In this paper we investigate the concept of untying a $p$--coloured
knot by $\pm1$--framed surgery. Any knot may be untied by surgery
along $\pm 1$--framed null-homotopic unknots. But such surgeries may
not preserve $p$--colourability. The remedy to this is only to allow
surgeries that preserve a $p$--colouring $\rho$ of $K$--- surgeries
by $\pm 1$--framed loops in the kernel of $\rho$. We call such
operations \emph{surgery in $\ker \rho$}.

It is natural to ask what the equivalence classes are of
$p$--coloured knots modulo such surgeries. Our main theorem is:

\begin{thm}\label{T:3case}
There are exactly $p$ equivalence classes of $p$--coloured knots
modulo surgery in $\ker\rho$ for $p=3,5$. These equivalence classes
are represented by connect-sums of $n$ left-hand $(p,2)$--torus
knots with a given colouring, for $n=1,2,\ldots p$.
\end{thm}

The above theorem tells us in particular that any $p$--coloured knot
$(K,\rho)$ is equivalent modulo $\ker\rho$ to a connect-sum of $n$
left-hand $(p,2)$--torus knots with a given colouring, with some
$n=1,2,\ldots p$.

That there are at most $p$ equivalence classes of $p$--coloured
knots modulo surgery in $\ker\rho$ is proved by taking a disc-band
presentation of the knot and first unlinking the bands (\fullref{S:Smn}) to get a connect-sum of genus $1$ knots. The number of
twists in each of the bands may then be $p$--reduced, and the
remaining knots may be reduced to $p$--coloured $(p,2)$--torus knots
explicitly. The final stage is to reduce the number of
connect-summands. These steps are carried out in \fullref{S:3Col} for $p=3$ and in \fullref{S:5Col} for $p=5$. That
there are at least $p$ such equivalence classes is proved in
\fullref{S:noteq} by finding a non-trivial invariant of $p$--coloured
knots which is invariant under surgery in $\ker\rho$. Connect-sums
of $n$ copies of a $(p,2)$--torus knot with a given colouring are
separated by this invariant, for $n=1,2,\ldots p$.

\begin{conj} \label{C:main}
The invariant of \fullref{S:noteq} is an complete invariant on
the set of equivalence classes modulo $\ker\rho$ for any odd prime
$p$. In other words \fullref{T:3case} holds any odd prime $p$.
\end{conj}

It is interesting to consider what the corresponding conjecture
should be when $p$ is not prime.

Note that the converse of this conjecture is easy--- if $K$ is a
$p$--coloured knot with $p$--colouring $\rho$ and $K'$ is a knot
obtained from $K$ by surgery by an element $C$ in $\ker \rho$, then
$K'$ is $p$--colourable with $p$--colouring $\rho'$ induced by
$\rho$. A presentation for the fundamental group of $S^{3}-K'$ is
given by the Wirtinger presentation of $K\cup C$ (see Kawauchi \cite{Kaw96}),
where the relations added to the set of generators of
$\pi_{1}(S^{3}-K)$ when we push it forward by the surgery map into
$\pi_{1}(S^{3}-K')$ are in the kernel of $\rho$. Thus $K'$ is
$p$--colourable, with the representation $\rho'$ induced by the
image of $\rho$ on the push-forward of the generators of
$\mathfrak{G}$.

Two applications of \fullref{T:3case} are given in \fullref{S:appl}.

Firstly, as pointed out to the author by Andrew Kricker, a surgery
presentation for $3$-- and $5$--coloured knots gives a surgery
presentation for irregular branched dihedral coverings of knots associated
to the dihedral groups $D_{6}$ and $D_{10}$. In this way they form the
basis for research into invariants of $p$--coloured knots in the spirit
of Garoufalidis and Kricker's work \cite{GK03,GK04}.
In \cite{Mos07a}, the non-commutative surgery presentation of coloured knots
given in this paper is used to define a non-commutative analogue of
the Kontsevich invariant for $3$-- and $5$--coloured knots.
Examining the loop expansion of such a Kontsevich invariant is expected to
provide dihedral analogues of such classical invariants as the Alexander
polynomial.

A second application of \fullref{T:3case}, which was pointed out to the
author by Makoto Sakuma, is that for a certain class of 3--manifolds they
enable us to expand a result of Przytycki and Sokolov \cite{PrzS01} (see
also Sakuma \cite{Sak01}) concerning surgery presentations of periodic
$3$--manifolds to the case of dihedral periods.  It is an interesting
problem to expand these results further to a wider class of manifolds
and periods.

\subsection*{Acknowledgements}
This paper owes a huge amount to a number of people. The author
would like to sincerely thank Tomotada Ohtsuki for his support and
encouragement, for suggesting substantial simplifications to some of
the proofs and reorganization of the material, for pointing out
various minor gaps, and for his topological interpretation of the
invariant in \fullref{S:noteq}.  The author would also like to
thank Andrew Kricker for providing motivation for the research that
lead to this paper and for his valuable comments and advice
throughout. Special thanks also to Makoto Sakuma for most useful
discussions including pointing out a gap in an earlier version of
the proof of \fullref{P:Bunlink}, to Kazuo Habiro for
important comments, encouragement and advice, especially for
pointing out an error in \fullref{S:5Col} in a previous version,
and to Steven Wallace for his careful reading and useful comments.
The author would finally like to thank the anonymous referee for
numerous comments and suggestions for improvements.

\section{Basic moves}

We begin by identifying moves which relate $p$--coloured tangles
$(T,\rho)$ through surgery in $\ker{\rho}$.

\begin{notation}
If $(T,\rho)$ and $(T',\rho')$ are related by a sequence of
$\pm1$--framed surgeries by components in $\ker{\rho}$, we write
$(T,\rho)\seq(T',\rho')$.
\end{notation}

\begin{notation}
It is convenient to denote half twists by numbers in boxes. By
convention, we set:
$$
\labellist\tiny
\pinlabel {$\tfrac12$} at 215 245
\endlabellist
\includegraphics[width=42pt]{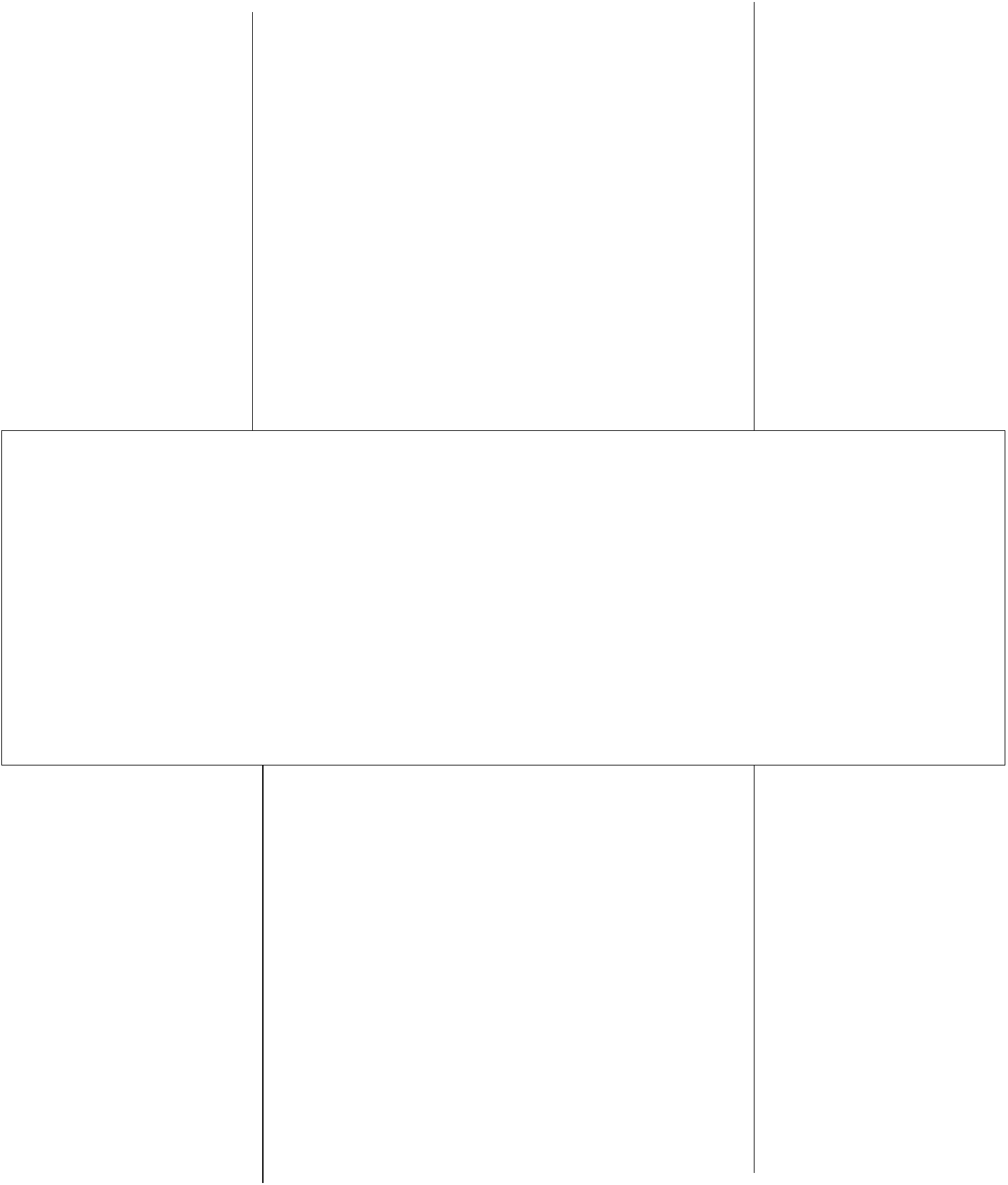}
\begin{picture}(20,42) \put(10,21){$=$} \end{picture}
\begin{picture}(42,42)
  \put(42,0){\line(-1,1){42}}
  \qbezier(0,0)(0,0)(14,14)
  \qbezier(28,28)(42,42)(42,42)
  \put(42,42){\line(1,1){0}}
\end{picture}
\qquad\qquad
\labellist\tiny
\pinlabel {$-\tfrac12$} at 215 245
\endlabellist
\includegraphics[width=42pt]{\figdir/twist}
\begin{picture}(20,42) \put(10,21){$=$} \end{picture}
\begin{picture}(42,42)
  \put(0,0){\line(1,1){42}}
  \qbezier(42,0)(42,0)(28,14)
  \qbezier(14,28)(0,42)(0,42)
  \put(0,42){\line(-1,1){0}}
\end{picture}$$
\end{notation}

The next proposition lists some basic moves which we would like to
use in our proofs.

\begin{prop}
The $RR$ move:
$$\begin{picture}(60,60)
    \put(3,60){$a$}
    \put(50,60){$a$}
    \put(0,0){\line(1,1){60}}
    \qbezier(60,0)(60,0)(40,20)
    \qbezier(20,40)(0,60)(0,60)
    \put(0,60){\line(-1,1){0}}
\end{picture}
\begin{picture}(20,60) \put(10,30){$\seq$} \end{picture}
\begin{picture}(60,60)
\put(3,60){$a$}
\put(50,60){$a$}
    \put(60,0){\line(-1,1){60}}
    \qbezier(0,0)(0,0)(20,20)
    \qbezier(40,40)(60,60)(60,60)
    \put(60,60){\line(1,1){0}}
\end{picture}$$
The $R2G$ move:
$$\labellist\small
\pinlabel {$a$} [b] at 11 437
\pinlabel {$b$} [b] at 103 437
\pinlabel {$b$} [b] at 193 437
\pinlabel {$\seq$} at 344 218
\endlabellist
\includegraphics[width=80pt]{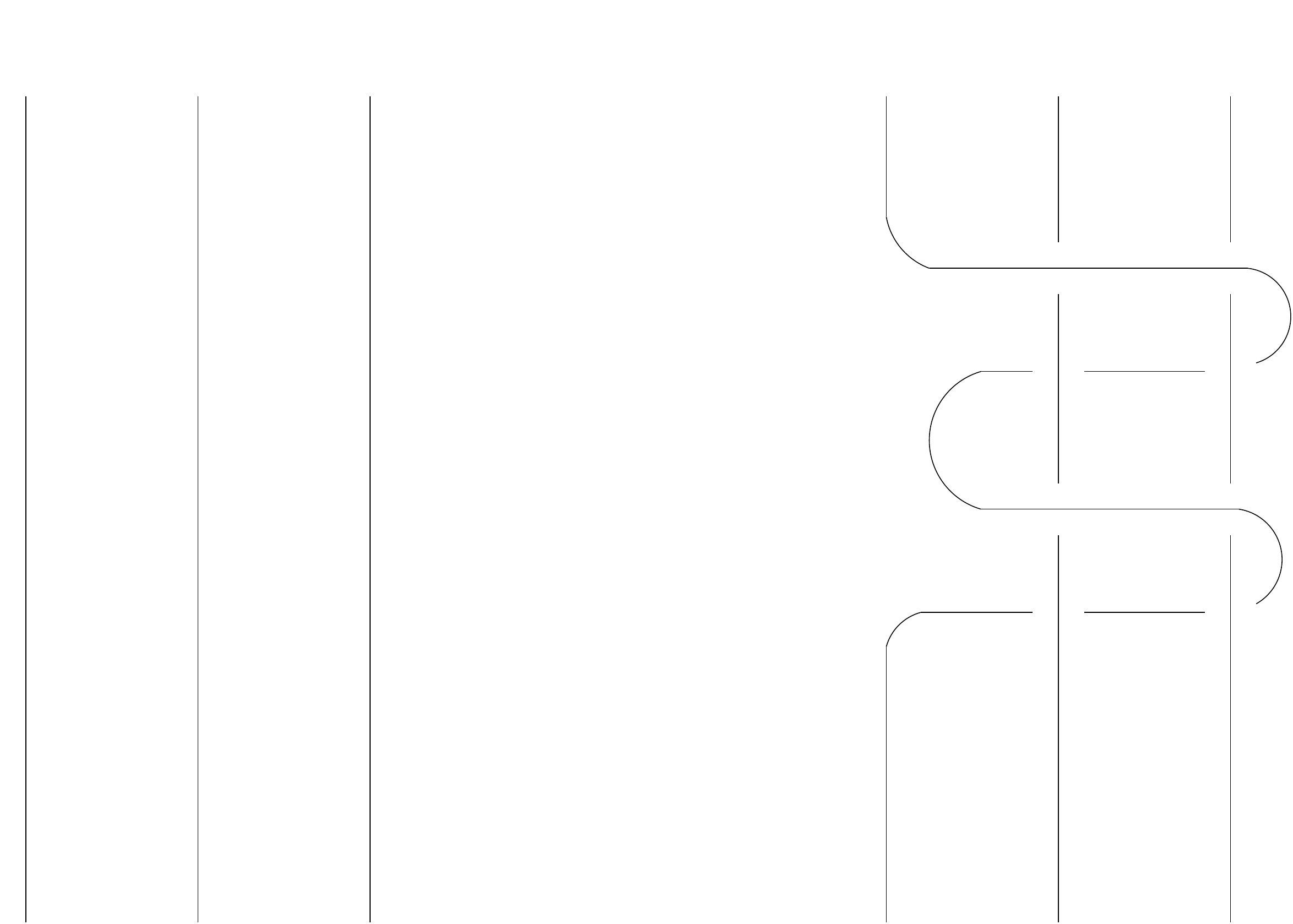}$$
\end{prop}

\begin{proof}
Note that since an element $ts^{i}$ is its own inverse, we
may ignore the orientation of the strands. Changing the orientation
of the surgery component reverses linkage.
The $RR$ move is the result of $1$--framed surgery on a component
which loops around both strands.

The $R2G$ move is the result of $1$--framed surgery on a component
which loops twice around the strand labeled $a$ and once around the
two stands labeled $b$. Various $RR$'s are necessary to clean up.
Graphically:
$$
\labellist\tiny
\pinlabel {$a$} [b] at 40 544
\pinlabel {$b$} [b] at 328 544
\pinlabel {$b$} [b] at 370 544
\endlabellist
\includegraphics[width=60pt]{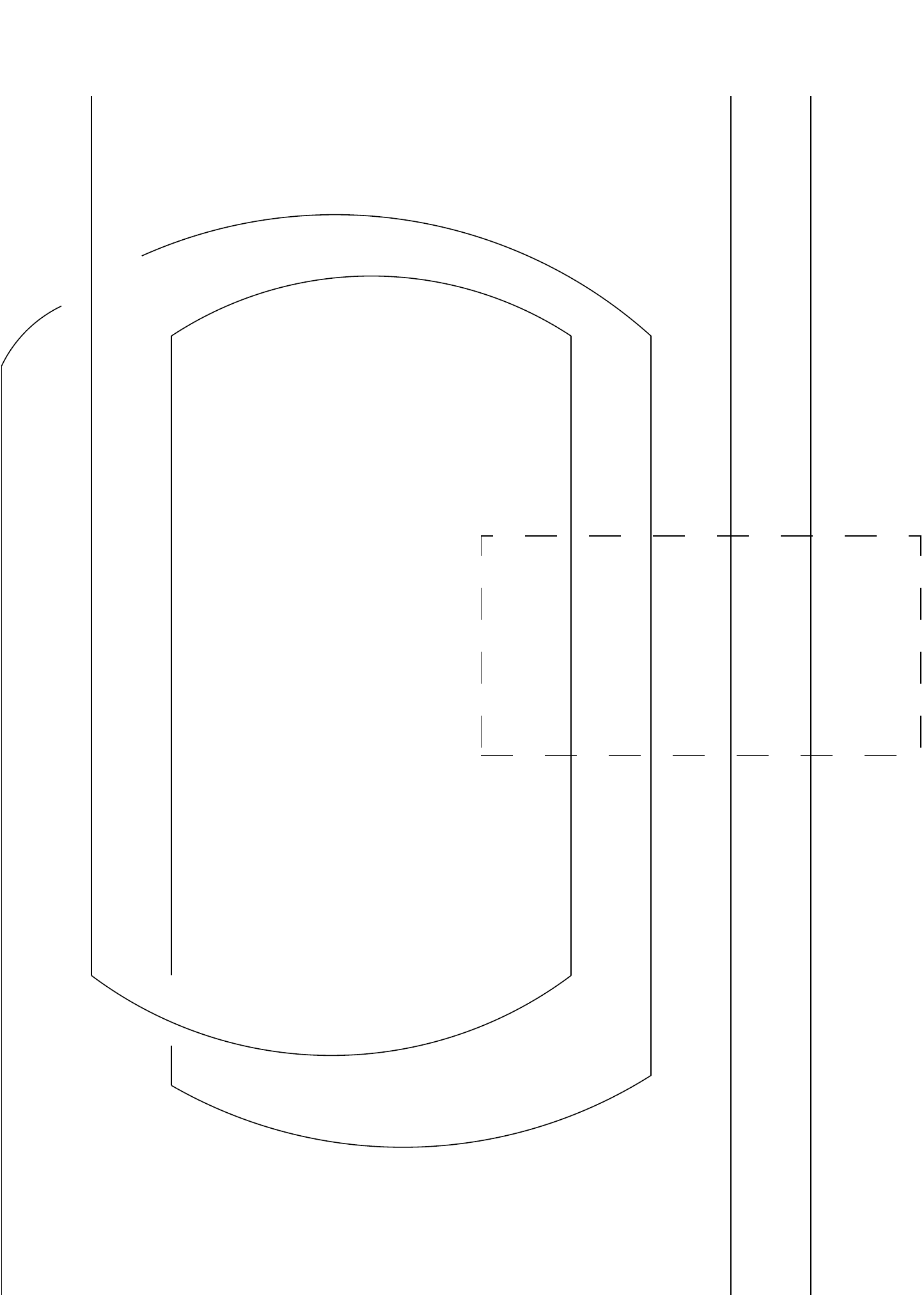}
\begin{picture}(40,70)
  \put(8,35){$\overset{\text{surgery}}{\sim}$}
\end{picture}
\labellist\tiny
\pinlabel {$a$} [b] at 40 544
\pinlabel {$b$} [b] at 328 544
\pinlabel {$b$} [b] at 370 544
\endlabellist
\includegraphics[width=60pt]{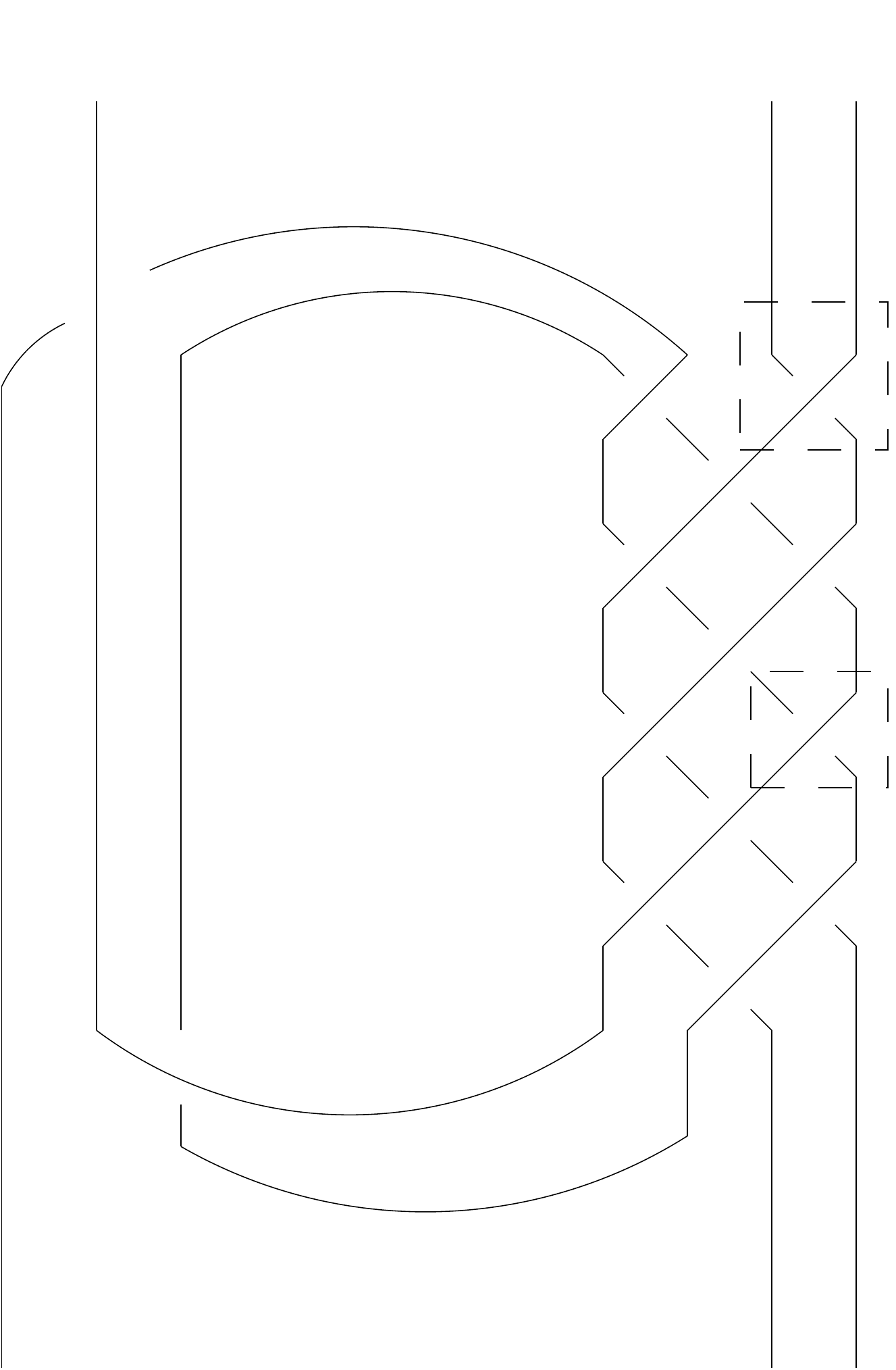}
\begin{picture}(40,70)
  \put(8,35){$\overset{2 \times RR}{\sim}$}
\end{picture}
\labellist\tiny
\pinlabel {$a$} [b] at 40 544
\pinlabel {$b$} [b] at 328 544
\pinlabel {$b$} [b] at 370 544
\endlabellist
\includegraphics[width=60pt]{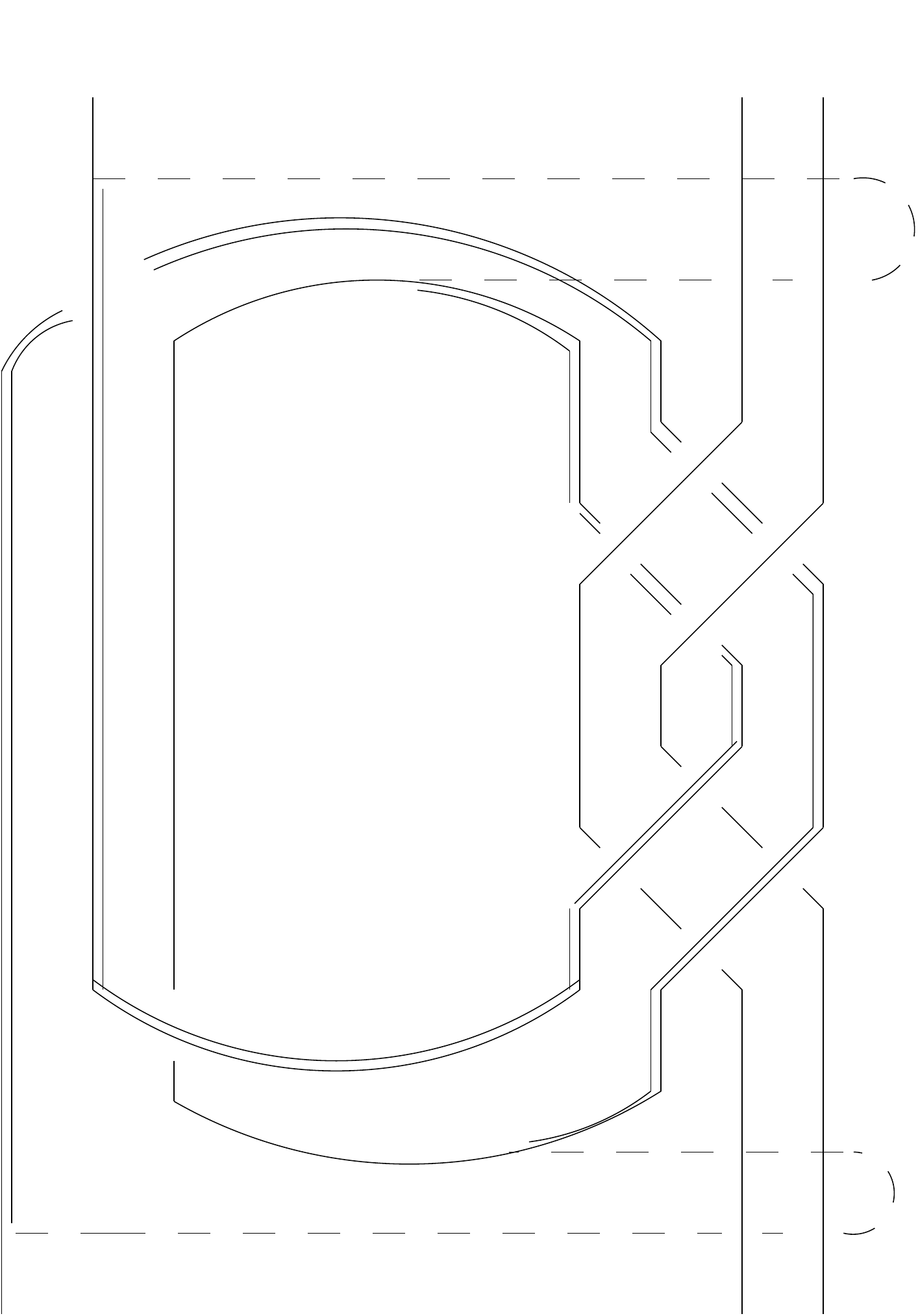}
\begin{picture}(30,70)
  \put(8,35){$=$}
\end{picture}
\labellist\tiny
\pinlabel {$a$} [b] at 2 435
\pinlabel {$b$} [b] at 92 435
\pinlabel {$b$} [b] at 182 435
\endlabellist
\includegraphics[width=40pt]{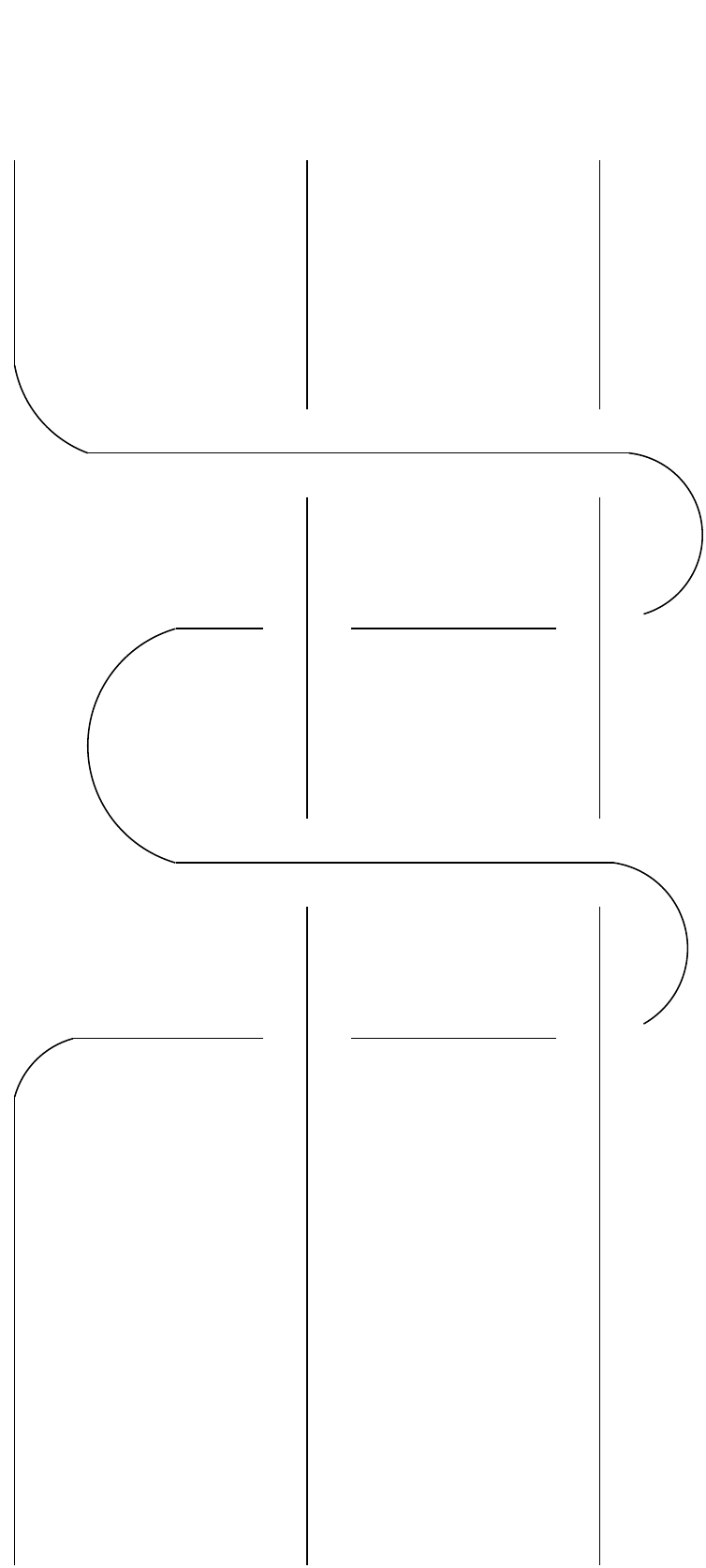}
$$
This completes the proof.
\end{proof}

The following move, which is a combination of $RR$ and of $R2B$, is
also useful.

\begin{cor}[$LT$ move]\label{C:LT}
$$\labellist\small
\pinlabel {$\seq$} at 260 230
\pinlabel {$3$} at 440 245
\endlabellist
\includegraphics[width=80pt]{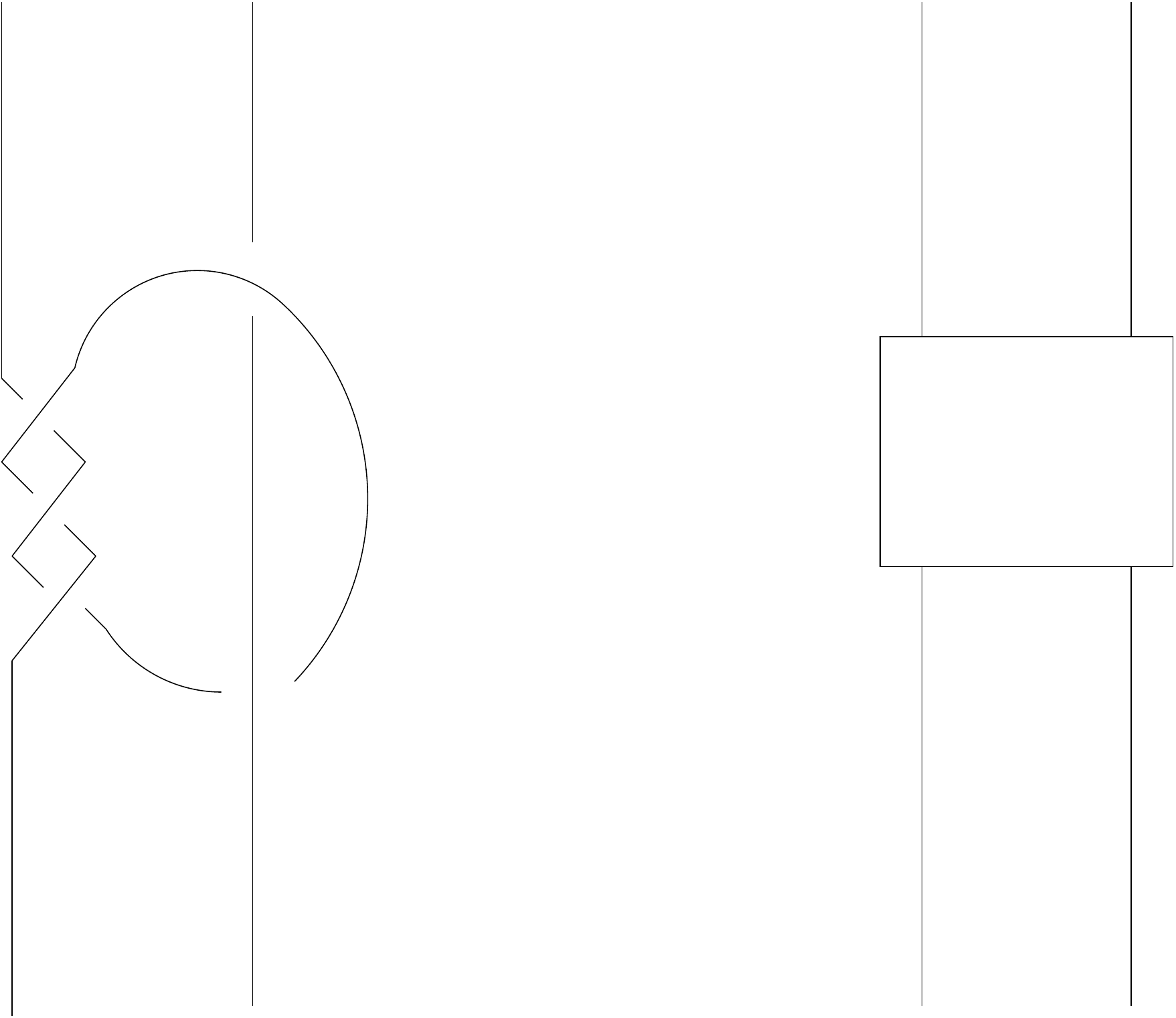}$$
\end{cor}

\begin{proof}
The colours below are to be read modulo $p$.
\begin{multline*}
\labellist\small
\pinlabel {$1$} [b] at 20 505
\pinlabel {$0$} [b] at 235 505
\endlabellist
\includegraphics[width=42pt]{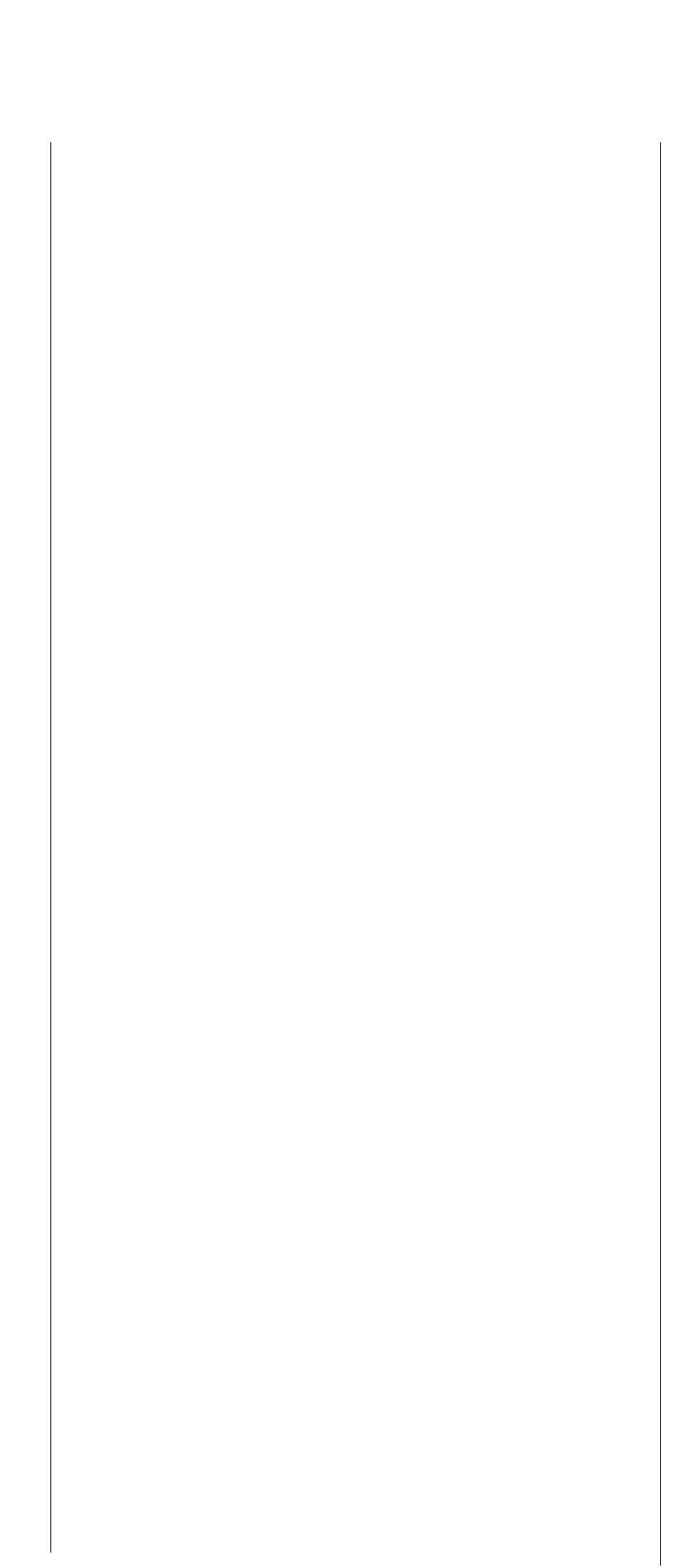}
\begin{picture}(40,70) \put(10,40){$\overset{R2G}{\sim}$} \end{picture}
\labellist\tiny
\hair=2pt
\pinlabel {$1$} [l] at 0 486
\pinlabel {$0$} [l] at 468 486
\pinlabel {$0$} [t] at 435 410
\pinlabel {$0$} [b] at 220 446
\pinlabel {$2$} [b] at 282 315
\pinlabel {$2$} [r] at 136 324
\pinlabel {$3$} [b] at 185 285
\pinlabel {$1$} [b] at 220 222
\pinlabel {$-1$} [t] at 176 144
\pinlabel {$0$} [t] at 306 38
\pinlabel {$1$} [l] at 0 20
\pinlabel {$0$} [l] at 570 20
\endlabellist
\includegraphics[width=120pt]{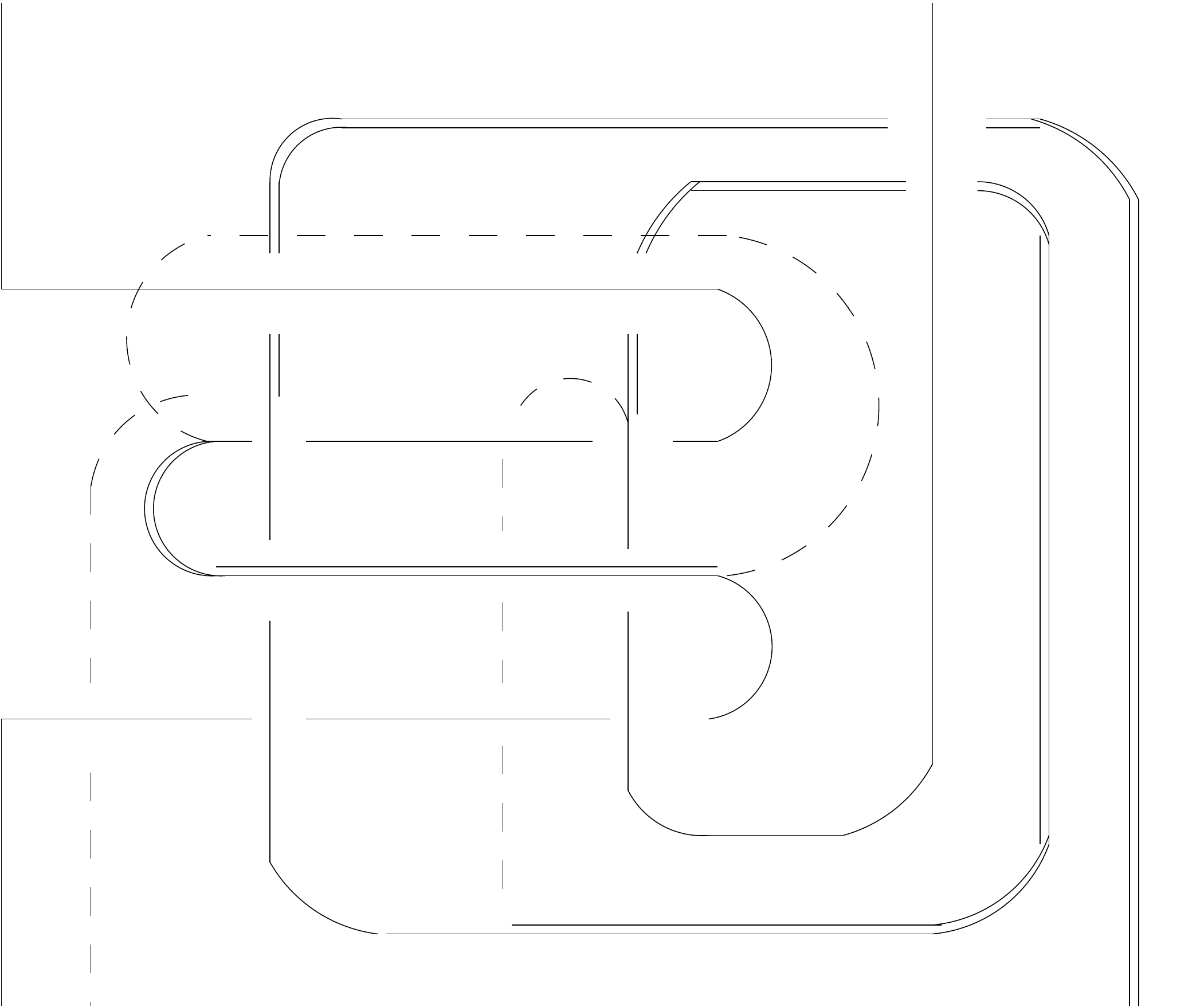}
\begin{picture}(40,70) \put(15,40){$=$} \end{picture}
\labellist\tiny
\hair=2pt
\pinlabel {$1$} [l] at 5 573
\pinlabel {$0$} [l] at 333 573
\pinlabel {$2$} [r] at 112 134
\pinlabel {$1$} [l] at 5 32
\pinlabel {$0$} [l] at 112 36
\pinlabel {$1$} [r] at 40 300
\pinlabel {$-1$} [l] at 480 260
\pinlabel {$1$} [bl] at 164 344
\pinlabel {$-1$} [b] at 395 263
\pinlabel {$-2$} [tl] at 312 96
\pinlabel {$0$} [b] at 136 224
\endlabellist
\includegraphics[width=90pt]{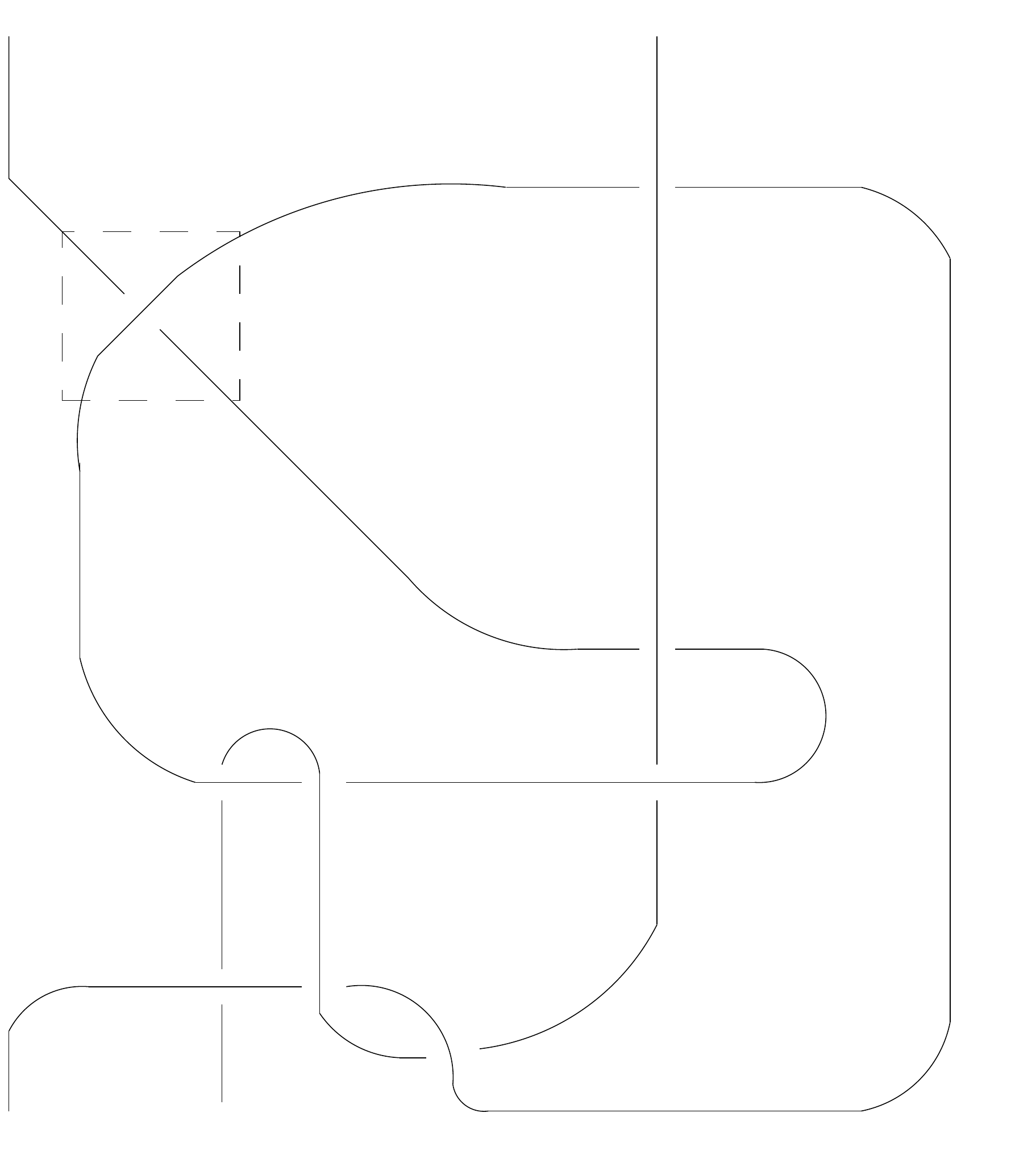} \\
\begin{picture}(40,70) \put(10,50){$\overset{RR}{\sim}$} \end{picture}
\labellist\small
\pinlabel {$0$} [l] at 248 18
\pinlabel {$1$} [r] at 70 18
\pinlabel {$3$} [l] at 247 180
\pinlabel {$4$} [r] at 72 180
\pinlabel {$0$} [l] at 252 396
\pinlabel {$1$} [r] at 40 396
\pinlabel {$2$} [t] at 157 255
\pinlabel {$3$} [r] at 40 277
\endlabellist
\includegraphics[width=80pt]{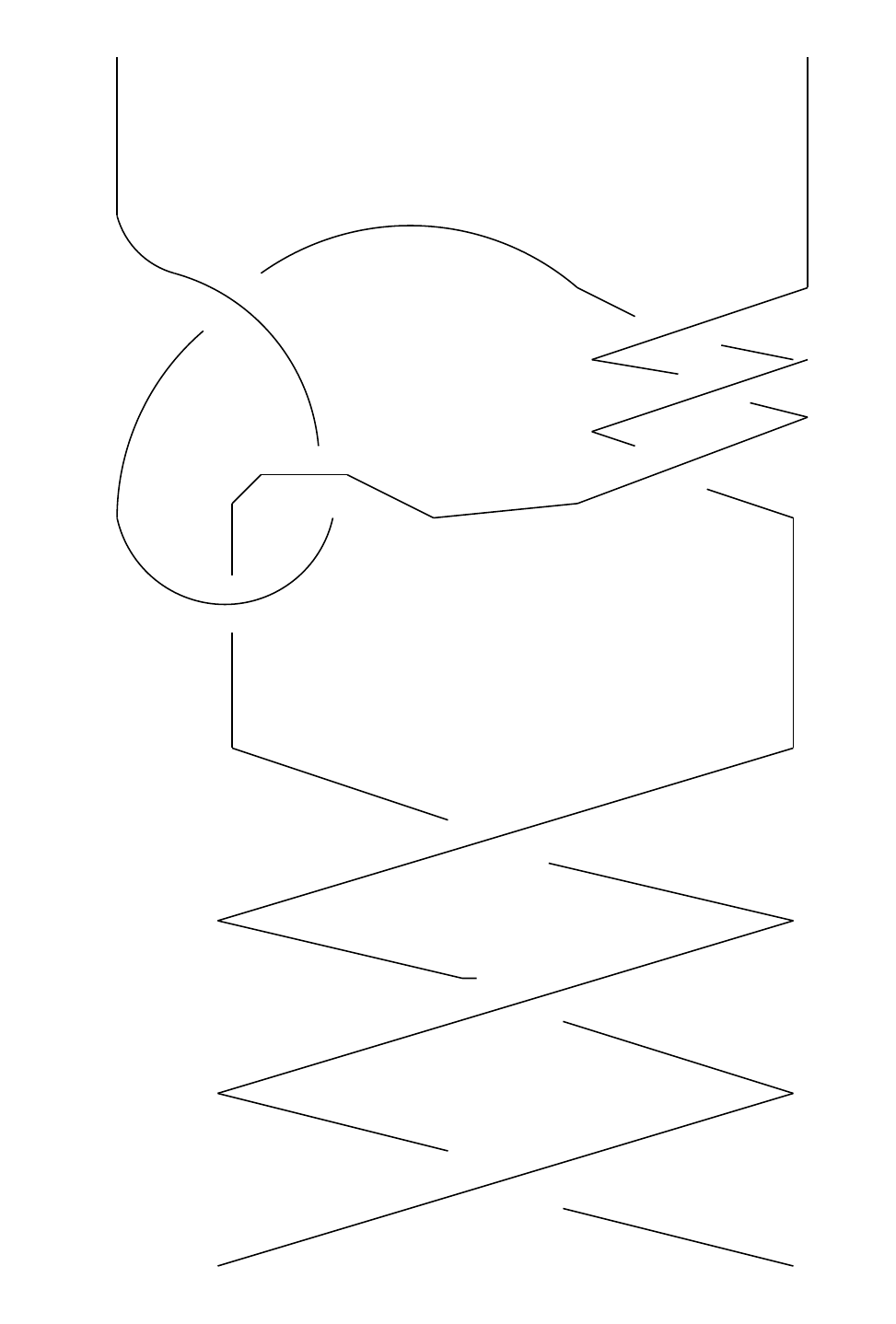}
\begin{picture}(40,70) \put(15,50){$\overset{RR}{\sim}$} \end{picture}
\labellist\small
\pinlabel {$1$} [l] at 31 15
\pinlabel {$0$} [r] at 175 15
\pinlabel {$-3$} at 105 105
\pinlabel {$1$} [r] at 31 405
\pinlabel {$0$} [l] at 175 413
\pinlabel {$3$} [b] at 110 345
\pinlabel {$7$} [r] at 36 200
\pinlabel {$6$} [l] at 176 178
\endlabellist
\includegraphics[width=60pt]{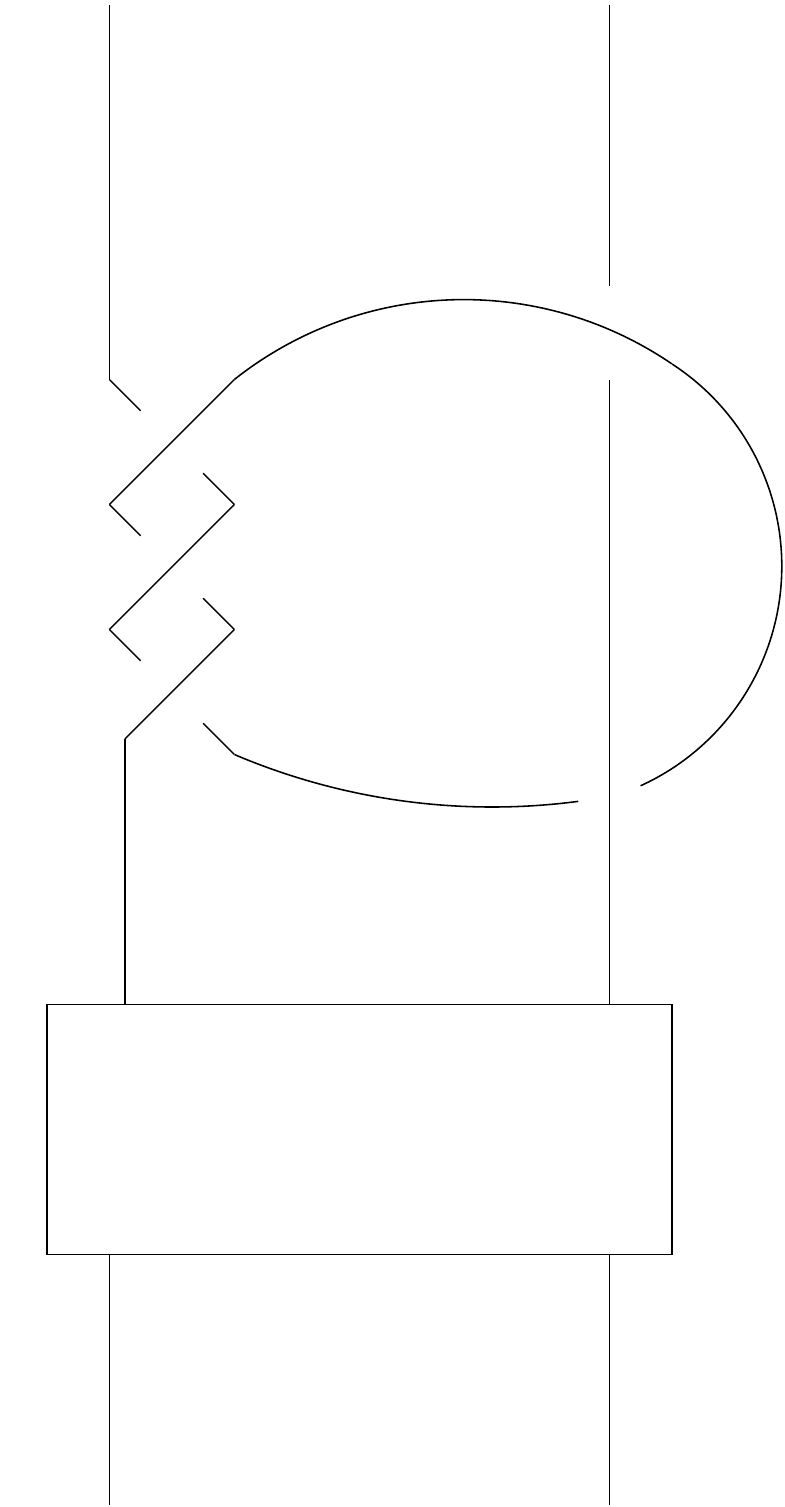}
\end{multline*}
This proves the corollary.
\end{proof}

Before we state the next corollary, we define a connect-sum for
$p$--coloured knots. Let $(K_{1},\rho_{1})$ and $(K_{2},\rho_{2})$
be two $p$--coloured knots. We equip both with orientations and
basepoints denoted $\star_{1}$ and $\star_{2}$ correspondingly. With
respect to this information the connect-sum
$(K_{1},\rho_{1})\hash_{(\star_{1},\star_{2})} (K_{2},\rho_{2})$ is
defined by cutting out neighbourhoods of $\star_{1}$ in $K_{1}$ and
of $\star_{2}$ in $K_{2}$ and gluing the ends by two parallel
overpaths. If the $p$--colouring of the arc containing $\star_{1}$
is different from the $p$--colouring of the arc containing
$\star_{2}$, we change the $p$--colouring of $K_{2}$ by a rotation
of the set of colours $\mathds{Z}_{p}$ before gluing. The colouring
of the connect-sum is then uniquely induced by the colouring of the
connect-summands.

\begin{lem}\label{L:consum}
A connect-sum of $p$--coloured knots is independent of the choices
of basepoints and for connect-sums where one of the connect-summands
is a $(p,2)$--torus knot is independent of the orientation of that
connect-summand.
\end{lem}

\begin{proof}
Note that an equivalent definition of connect-summing $p$--coloured
knots would be to cut out a neighbourhood of $\star_{2}$ to obtain a
$p$--coloured tangle $T$, to make $T$ `very small', and then to glue
$T$ into $K_{1}$ in place of a neighbourhood of $\star_{1}$
according to the orientations of $K_{1}$ and $T$, with
$p$--colourings as before.

Let $\star_{1}$ and $\star'_{1}$ be basepoints of
$(K_{1},\rho_{1})$, and $\star_{2}$ be a basepoint of
$(K_{2},\rho_{2})$. Begin with
$(K_{1},\rho_{1})\hash_{(\star_{1},\star_{2})} (K_{2},\rho_{2})$ and
making $T$ `very small', slide it along $(K_{1},\rho_{1})$ to a
neighbourhood of $\star'_{1}$. As $T$ moves under an over-crossing
arc of $K_{1}$ coloured $c\in\mathds{Z}_{p}$, its colouring changes
by a reflection by $c$ of $\mathds{Z}_{p}$. We may reverse this
reflection (so that we remain with only a rotation on the set of
colours) by the following sequence in which a loop is added by a
Reidemeister \textrm{I} move, $T$ is pushed through it by ambient
isotopy, and the loop is deleted by another Reidemeister \textrm{I}
move:
\begin{equation}\label{E:Loopit}
\begin{minipage}{60pt}
\labellist\tiny
\pinlabel {$T$} at 155 25
\endlabellist
\includegraphics[width=60pt]{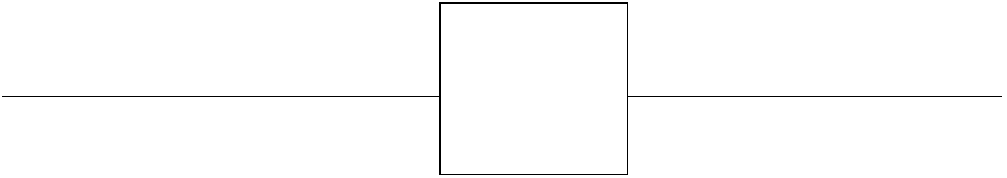}
\end{minipage}
\qua = \qua
\begin{minipage}{70pt}
\labellist\tiny
\pinlabel {$T$} at 85 55
\endlabellist
\includegraphics[width=70pt]{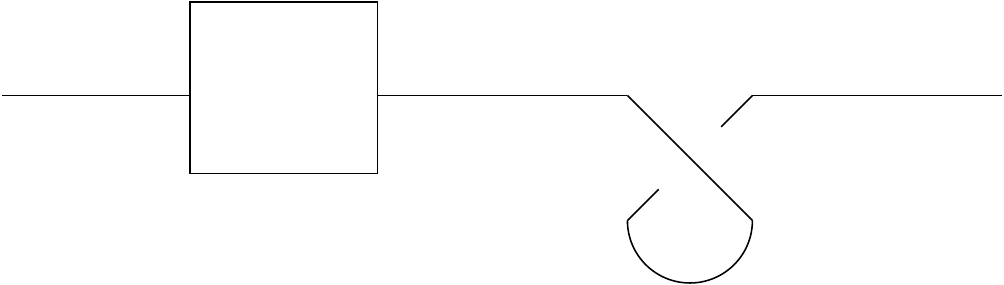}
\end{minipage}
\qua = \qua
\begin{minipage}{70pt}
\labellist\tiny
\pinlabel {$T$} at 210 55
\endlabellist
\includegraphics[width=70pt]{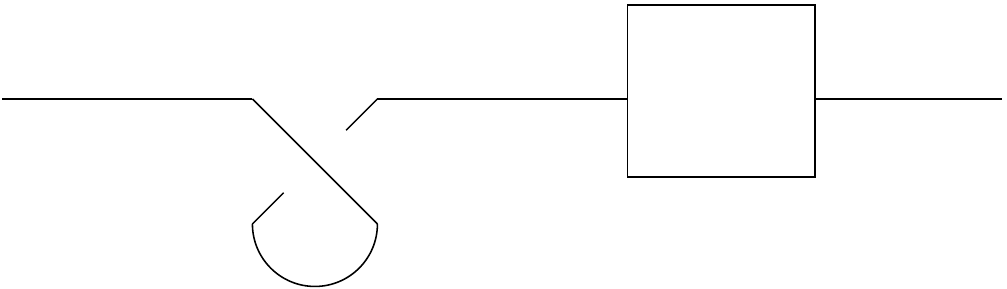}
\end{minipage}
\qua = \qua
\begin{minipage}{60pt}
\labellist\tiny
\pinlabel {$T$} at 155 25
\endlabellist
\includegraphics[width=60pt]{\figdir/reflpcol3}
\end{minipage}
\end{equation}
Thus $(K_{1},\rho_{1})\hash_{(\star_{1},\star_{2})}
(K_{2},\rho_{2})=(K_{1},\rho_{1})\hash_{(\star'_{1},\star_{2})}(K_{2},\rho_{2})$
and the connect-sum of $p$--coloured knots does not depend on the
choice of basepoint.

If $(K_{2},\rho_{2})$ is a $(p,2)$--torus knot, then $K_{2}$ is
invertible. The ambient isotopy from $K_{2}$ to $K_{2}$ with the
opposite orientation is given by taking a diagram of $K_{2}$ as
below and rotating the plane with respect to which the diagram is
taken by $\pi$ (this is the order $2$ symmetry of the $(p,2)$--torus
knot--- see Kawauchi \cite{Kaw96}).
$$\begin{minipage}{80pt}
\centering
\labellist\small
\pinlabel {$\tfrac{p}{2}$} at 120 130
\endlabellist
\includegraphics[width=80pt]{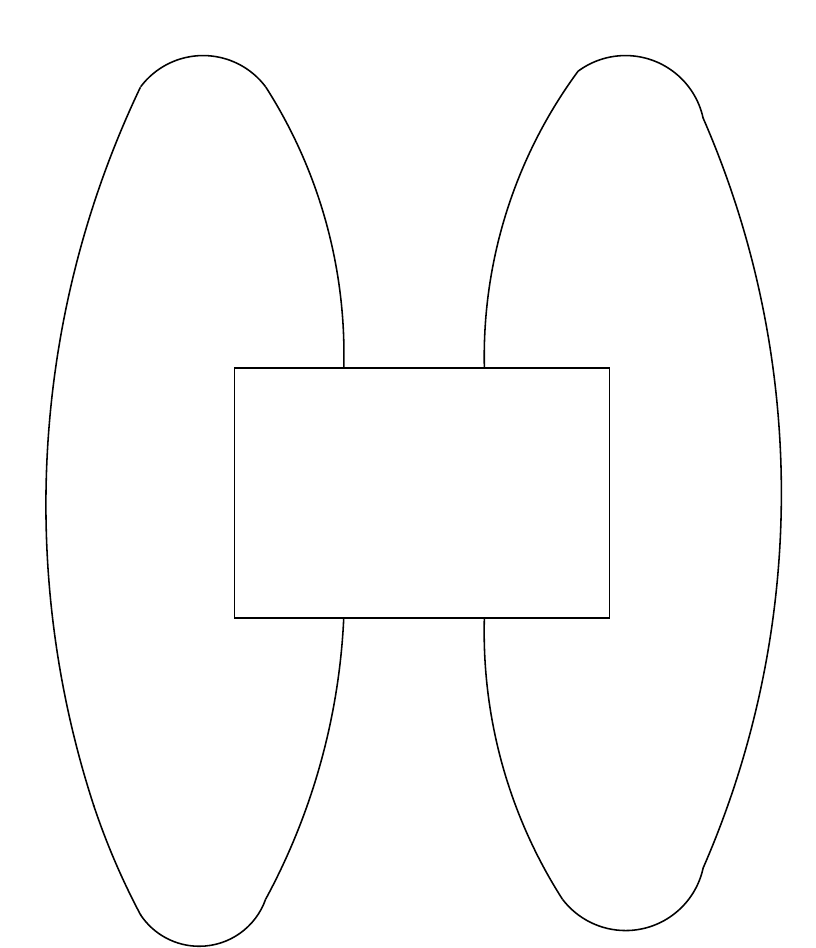}
\end{minipage}$$
Such an ambient isotopy reflects the set of
colours $\mathds{Z}_{p}$. As before, we may reverse this reflection
by sequence \ref{E:Loopit}, giving us the connect-sum of
$(K_{2},\rho_{2})$ with the opposite orientation with
$(K_{1},\rho_{1})$.
\end{proof}

From now on let $(K_{1},\rho_{1})\hash (K_{2},\rho_{2})$ denote the
connect-sum of $(K_{1},\rho_{1})$ and $(K_{2},\rho_{2})$. When one
of the connect summands is a $(p,2)$--torus knot this is independent
of orientation, and in other cases orientations of connect-summands
are mentioned explicitly in the text.

If $p=3$ the statement of \fullref{C:LT} simplifies even
further, giving us what may be thought of as the key lemma in the
proof for the three coloured case.

\begin{cor}\label{N:pTW}
For $p=3$, adding $3$ twists in a band is equivalent modulo surgery
in $\ker{\rho}$ to taking a connect sum with a left-hand trefoil.
\end{cor}

\begin{proof}
\begin{gather*}
\begin{minipage}{60pt}
\labellist\small
\pinlabel {$1$} [l] at 31 15
\pinlabel {$0$} [r] at 175 15
\pinlabel {$-3$} at 105 105
\pinlabel {$1$} [r] at 31 405
\pinlabel {$0$} [l] at 175 413
\pinlabel {$0$} [b] at 110 345
\pinlabel {$1$} [r] at 36 200
\pinlabel {$0$} [l] at 176 178
\endlabellist
\includegraphics[width=60pt]{\figdir/linkedtref}
\end{minipage}
\quad\overset{RR}{\sim}\quad
\begin{minipage}{60pt}
\labellist\small
\pinlabel {$1$} [l] at 30 15
\pinlabel {$0$} [r] at 175 15
\pinlabel {$-3$} at 105 107
\pinlabel {$1$} [l] at 32 195
\pinlabel {$1$} [r] at 32 370
\pinlabel {$0$} [r] at 175 365
\pinlabel {$0$} [bl] at 100 350
\endlabellist
\includegraphics[width=60pt]{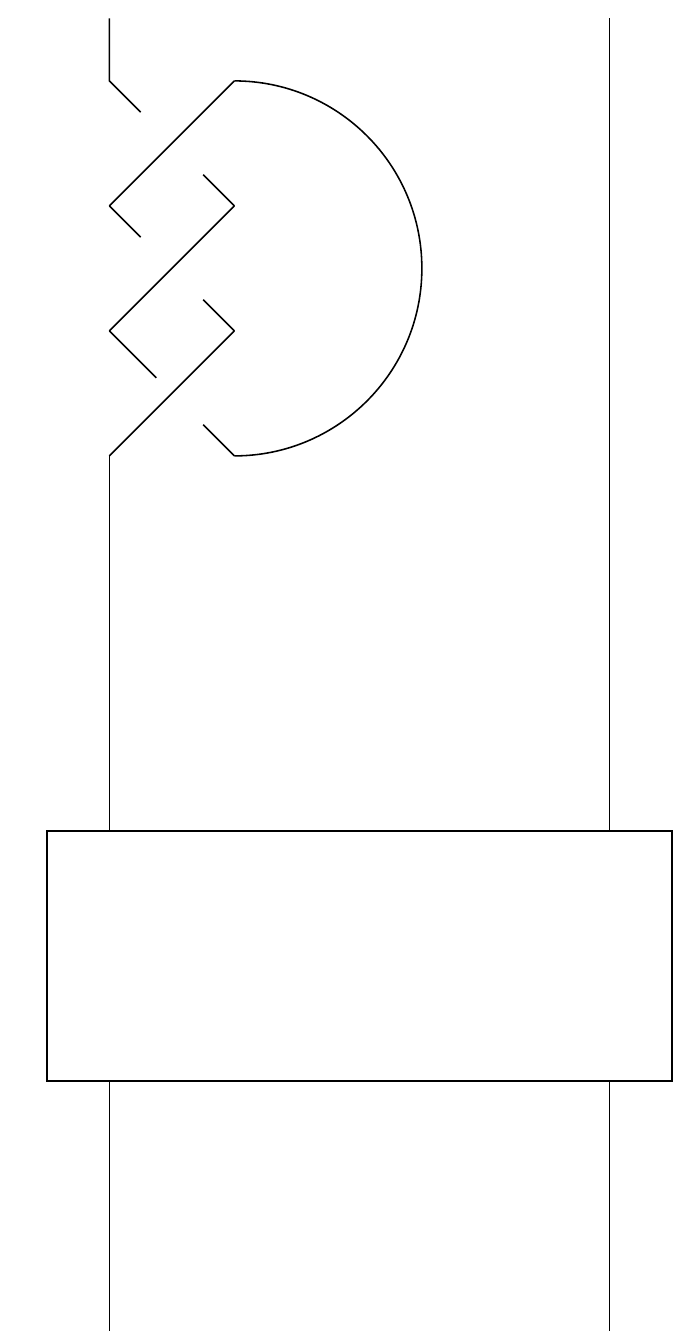}
\end{minipage}
\end{gather*}
This proves the corollary.
\end{proof}

\section{Presenting a $p$--coloured knot as a connect-sum of
$p$--col\-our\-ed knots of the form $S(m,n)$}\label{S:Smn}

Let $S(m,n)$ (S for Seifert) denote a genus $1$ knot with a single
negative half-twist between its bands, $m$ full twists in its left
band, and $n$ in its right band. Note that by sliding the leftmost
band around the disk and reflecting we see that the mirror image of
$S(m,n)$ is $S(-m,-n)$. In this section we show that any
$p$--coloured knot may be reduced to a connect-sum of $p$--coloured
knots of the form $S(m,n)$.

We now state the condition for an $S(m,n)$ knot to be
$p$--colourable.

\begin{prop}\label{P:prop9}
$S(m,n)$ is $p$--colourable if and only if $4mn\equiv 1\mod p$, and
if $4mn\not\equiv 1\mod p^{2}$ then this $p$--colouring is unique (up
to automorphisms of $D_{2p}$).
\end{prop}

\begin{proof}
The Seifert matrix of $S(m,n)$ is $M=\smash{\bigl(\begin{smallmatrix}m &
1\\ 0 & n
\end{smallmatrix}\bigr)}$, and thus its determinant is
$\det(M{+}M^{T})=\abs{H_{1}(\hC(S(m,n)))}=4mn{-}1$, where
$\abs{H_{1}(\hC(S(m,n)))}$ denotes the order of the first homology
group of the $2$--fold covering of $S^{3}$ branched over $S(m,n)$.
If $p$ divides this number but $p^{2}$ does not, then
$H_{1}(\hC(S(m,n)))$ contains a single copy of $\mathds{Z}/p\mathds{Z}$, and
therefore $S(m,n)$ has a $p$--colouring if and only if $4mn{-}1\equiv
0\mod p$, and the colouring is unique if $4mn{-}1\not\equiv 0\mod
p^{2}$. See Fox \cite{Fox62} for details.
\end{proof}

\begin{figure}
\label{fig1}
\labellist\small
\pinlabel {$m$} at 33 150
\pinlabel {$n$} at 320 150
\endlabellist
\centering
\includegraphics[width=2in]{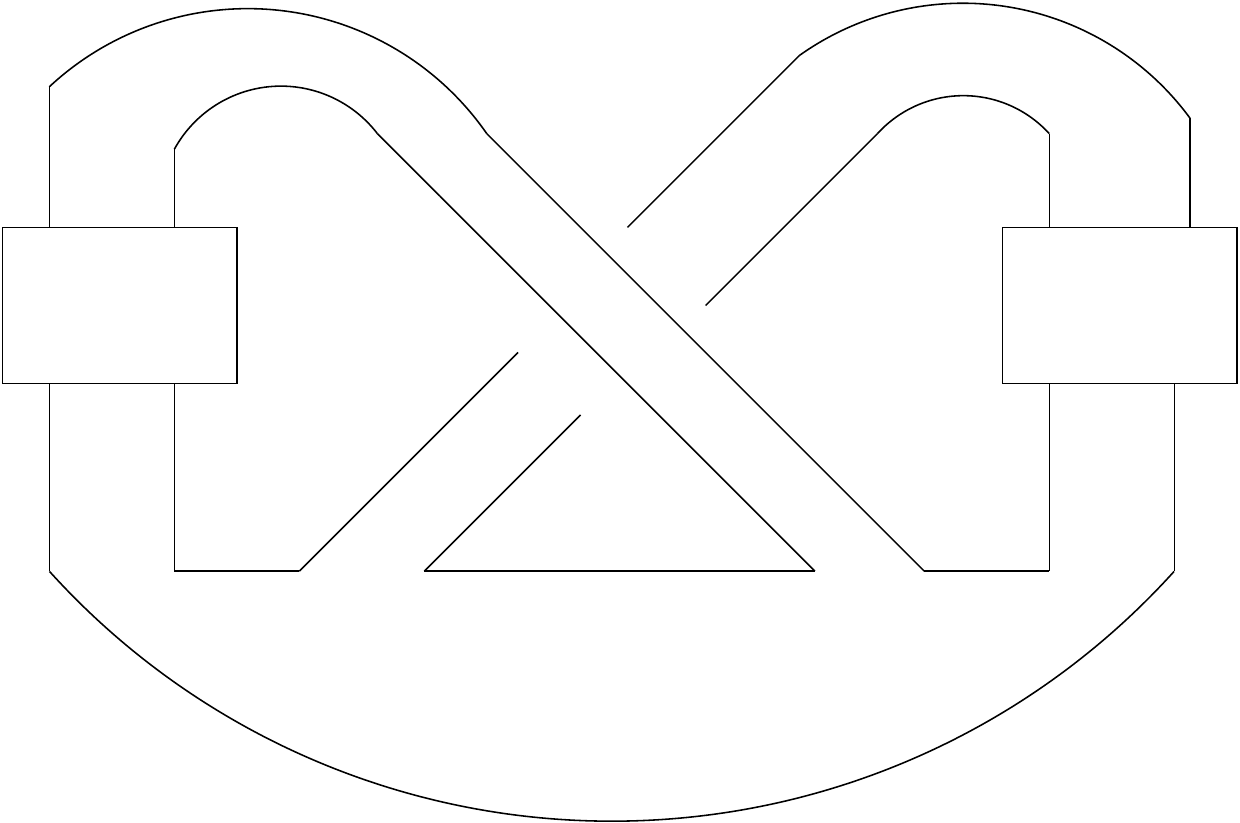}
\caption{The knot $S(m,n)$}
\end{figure}

All knots are assumed to be given in their band projections (see for
instance Burde--Zieschang \cite[Proposition 8.2]{BZ03}). We remind the reader that
this is defined as a projection of a knot $K$ which is represented
in $S^{3}$ as the boundary of an orientable surface $F$ in $S^{3}$
with the following properties:

\begin{enumerate}
\item $S=D^{2}\cup B_{1}\cup\cdots \cup B_{2n}$ where $D^{2}$ and each
$B_{j}$ is a disc.
\item $B_{i}\cap B_{j}=\emptyset$ for $i\neq j$, $\partial
B_{i}=\alpha_{i}\gamma_{i}\beta_{i}(\tilde{\gamma}_{i})^{-1}$,
$D^{2}\cap B_{i}=\alpha_{i}\cup\beta_{i}$, $\partial
D^{2}=\alpha_{1}\delta_{1}\beta_{2}^{-1}\delta_{2}\beta_{1}^{-1}\delta_{3}\alpha_{2}\delta_{4}\ldots
\alpha_{2n-1}\delta_{4n-3}\beta_{2n}^{-1}\delta_{4n-2}\beta_{2n-1}^{-1}\delta_{4n-1}\alpha_{2n}\delta_{4n}$.
\end{enumerate}

In this section, twists in the bands themselves are ignored. We show
that for any $p$, the bands of a knot may be unlinked by
$\pm1$--framed surgery in $\ker{\rho}$.  This implies:

\begin{prop} \label{P:Bunlink}
For any $p$--coloured knot $(K,\rho)$ of genus $n$ there exists a
natural number $g$ and integers $m_{1},m_{2},\ldots,m_{g}$ and
$n_{1},n_{2},\ldots,n_{g}$ such that
$$K\seq S(m_{1},n_{1})\hash S(m_{2},n_{2})\hash\cdots\hash S(m_{g},n_{g})$$
\noindent and each connect-summand has a non-trivial $p$--colouring
induced by $\rho$, and the connect-sum is with respect to some
orientation of the connect-summands.
\end{prop}

The way in which we unlink depends on the colouring of the arcs
which border the bands. We thus define the following `invariant' of
bands.

\begin{defn}\label{D:index}
Let $a,b,c,d\in\mathds{Z}_{p}$ be the colours of the arcs of the
band $A$ where it connects to the disc, read from left to right (in
the notation above, for $A=\gamma_{i}\cup(\tilde{\gamma}_{i})^{-1}$,
$a$ is the colour of the arc which contains the point
$\gamma_{i}\cap\alpha_{i}$, $b$ of the arc which contains the point
$\tilde{\gamma}_{i}\cap\alpha_{i}$, $c$ of
$\gamma_{i}\cap\beta_{i}$, and $d$ of
$\tilde{\gamma}_{i}\cap\beta_{i}$). Let
$\underline{a},\underline{b},\underline{c},\underline{d}\in\set{0,1,2\ldots,p-1}$
represent $a,b,c$ and $d$. We define
$$\abs{A}:=
  \min(\abs{\underline{a}-\underline{b}},p-\abs{\underline{a}-\underline{b}})
=\min(\abs{\underline{d}-\underline{c}},p-\abs{\underline{c}-\underline{d}})$$
\noindent the number $\abs{A}$ is independent of twists in $A$ and
of crossings of $A$ with other bands. It is called the \emph{index}
of $A$. It corresponds to the image under $\rho$ of a path which
loops once around the band.
\end{defn}

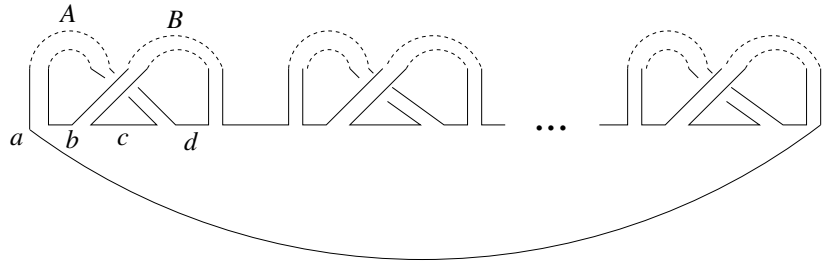
\begin{figure}\label{Fi:genseif}
\centering
\labellist\small
\pinlabel {$A$} [b] at 38 221
\pinlabel {$B$} [b] at 140 215
\pinlabel {$a$} [tr] at 0 127
\pinlabel {$b$} [tl] at 30 130
\pinlabel {$c$} [t] at 90 130
\pinlabel {$d$} [t] at 155 130
\endlabellist
\includegraphics[width=300pt]{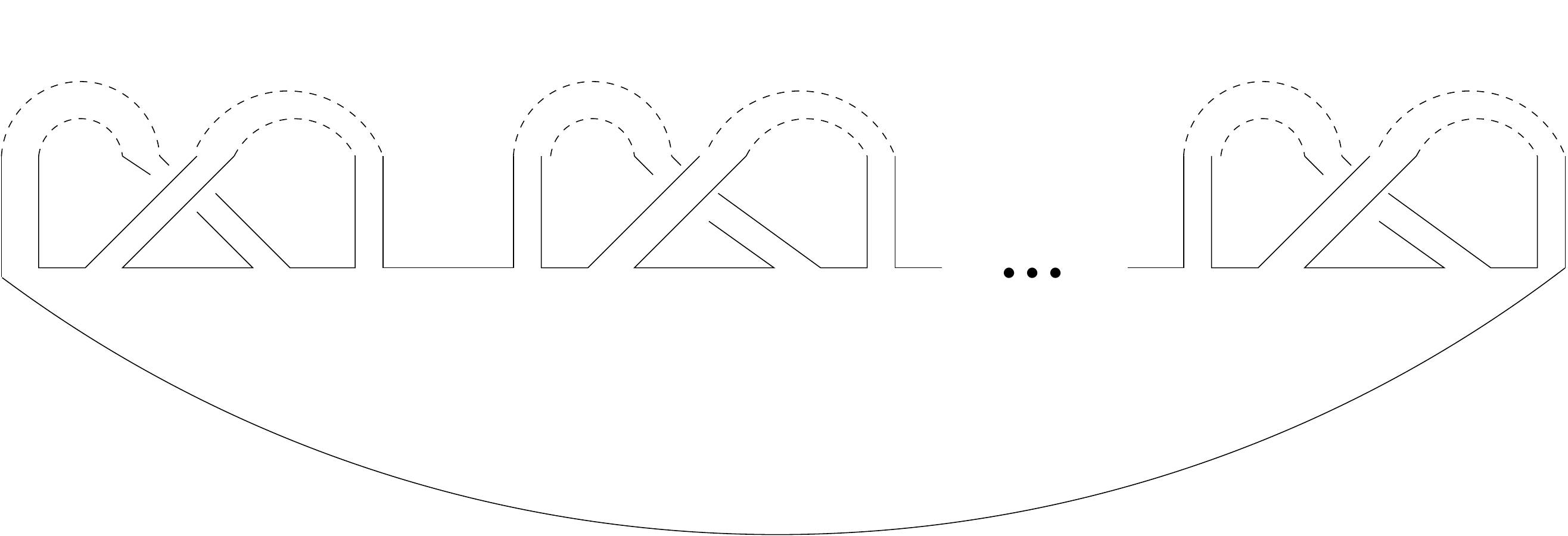}
\caption{A band projection. Each band may twist and link with any
other band.}
\end{figure}

In the simplest case, bands may be unlinked by surgery along a
single component which circles them once.

\begin{lem}\label{L:aeqb}
If $\abs{A}=\abs{B}$ then $A$ and $B$ may be unlinked in
$\ker{\rho}$.
\end{lem}

\begin{proof}
If $\abs{A}=\abs{B}$, let $C$ be a $1$--framed unknot which circles
$A$ and $B$ once. As we may add and subtract a single twist in
either band, we may assume that the strands ringed by $C$ are
coloured $a$, $a+r$, $b$, and $b-r$, and that therefore
$C\in\ker{\rho}$. Surgery along $C$ unlinks the bands.
\end{proof}

In the next simplest case, to unlink the bands we require more
components and some basic number theory.

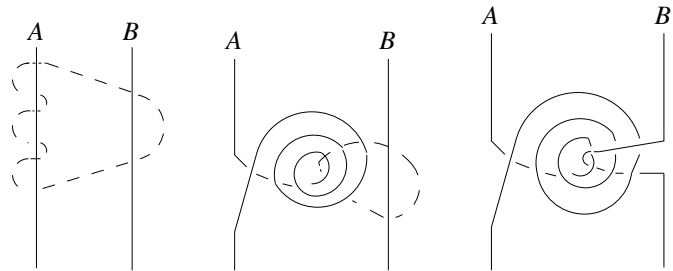
\begin{figure}\label{Fi:swirl} \centering
\labellist\small
\pinlabel {$A$} [b] at 55 505
\pinlabel {$B$} [b] at 270 505
\endlabellist
\includegraphics[width=0.8in]{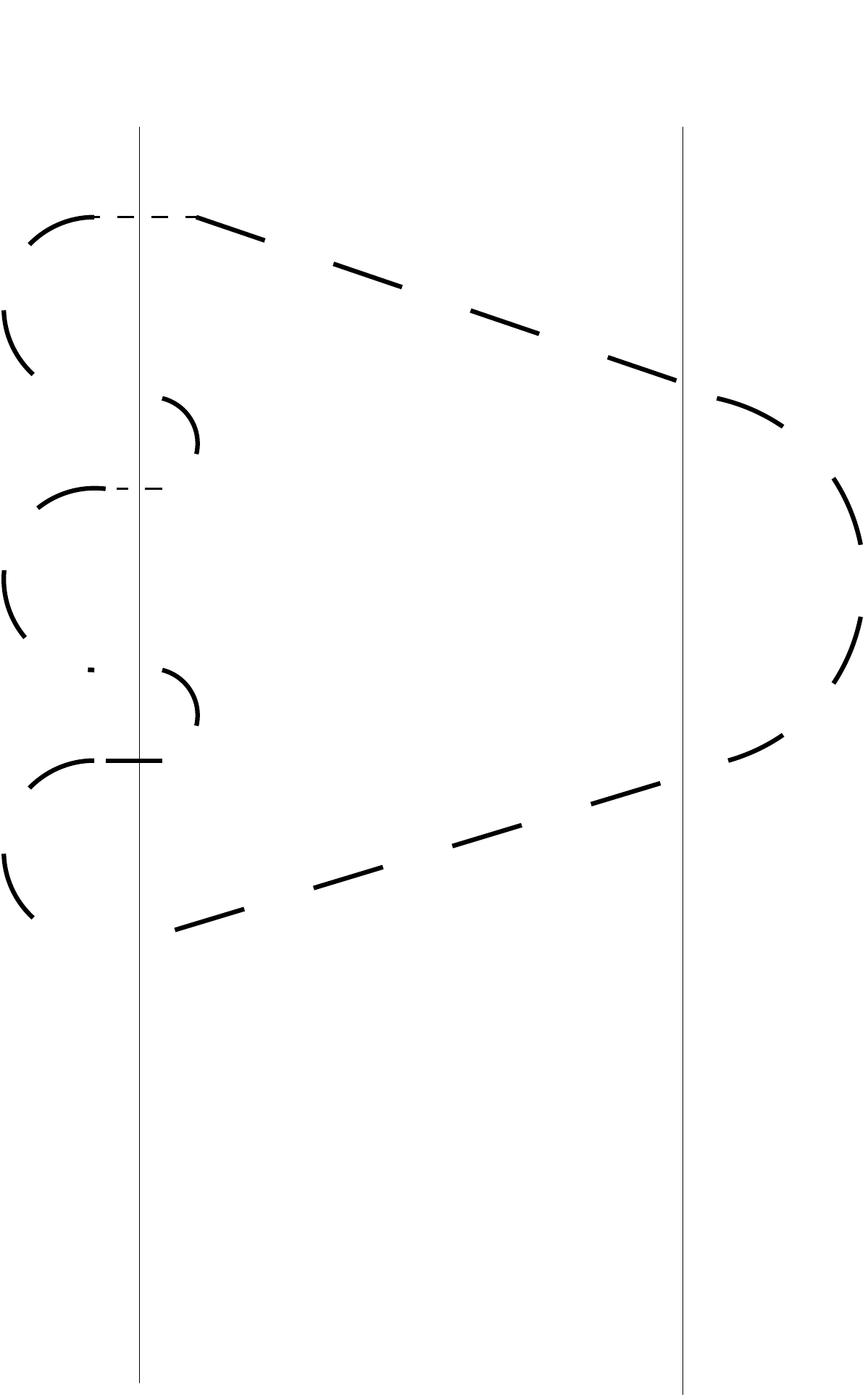}
\qquad
\labellist\small
\pinlabel {$A$} [b] at 8 400
\pinlabel {$B$} [b] at 296 400
\endlabellist
\includegraphics[width=1in]{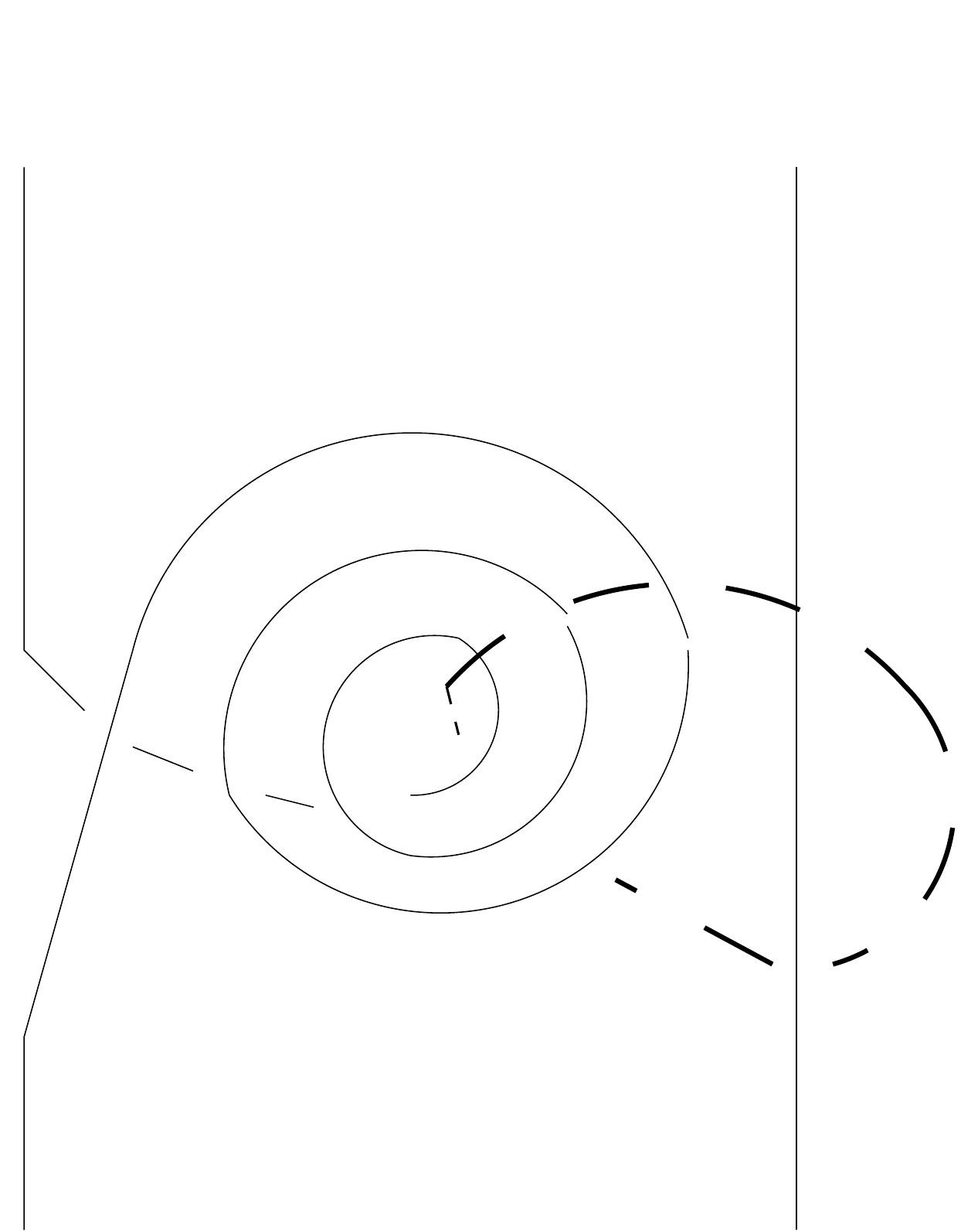}
\qquad
\labellist\small
\pinlabel {$A$} [b] at 8 400
\pinlabel {$B$} [b] at 296 400
\endlabellist
\includegraphics[width=1in]{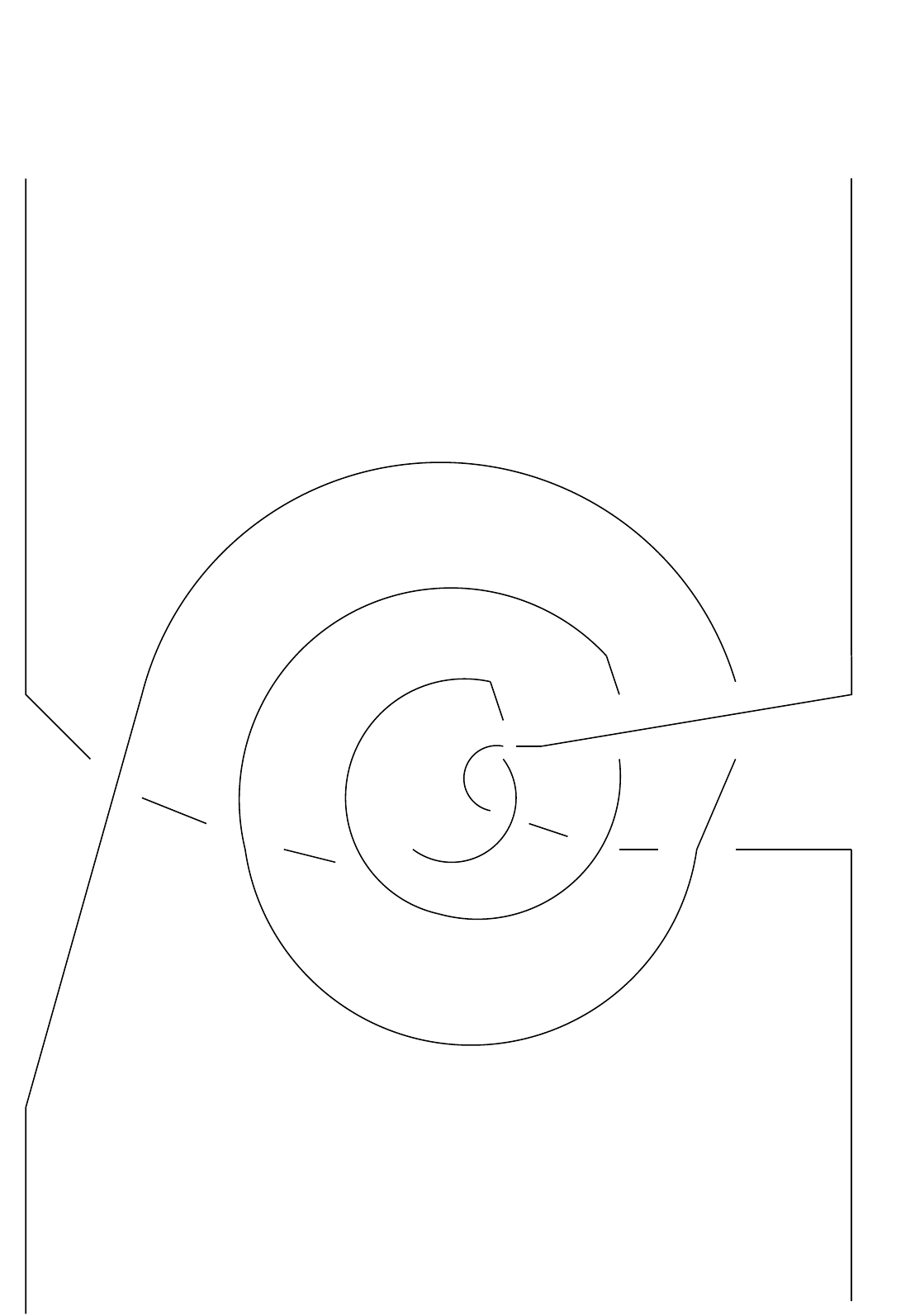}
\caption{Adding $3$ twists between bands $A$ and $B$ for which
$3\abs{A}+\abs{B}\equiv 0\mod{p}$. The solid lines represent bands,
which we allow to pass through themselves freely by \fullref{L:aeqb},
while the dotted line represents the surgery component.}
\end{figure}

\begin{lem}\label{L:aneqb}
If $\abs{A}\neq\abs{B}$ are both nonzero, then $A$ and $B$ may be
unlinked in $\ker{\rho}$.
\end{lem}

\begin{proof}
Pick $a$ and $b$ coprime to $p$ such that $a\abs{A}+b\abs{B}\equiv
0\mod{p}$. Let $C$ be the component in $H_{1}(S^{3}-K)$ which (from some
choice of base-point) circles $a$ times around $A$ and then $b$ times
around $B$. The surgery component $C$ is in $\ker{\rho}$, and the effect
of surgery by $C$, by applying the previous lemma (which allows us to
pass a band through itself) is to add $ab$ twists between the bands $A$
and $B$. The example $a=3$ and $b=1$ is illustrated in \fullref{Fi:swirl}.

Because $p$ is prime, we may assume by an automorphism of $D_{2p}$
that $\abs{A}=a=1$. We have $1+b\abs{B}\equiv 0\mod{p}$ and also
$1+(b-p)\bmod{p}\abs{B}\equiv 0\bmod{p}$, hence we may add or
subtract $b$ twists and $p-b$ twists between $A$ and $B$. Since $p$
is prime, $b$ and $p-b$ are coprime and therefore there exist
integers $m$ and $n$ such that $mb+n(p-b)=1$. So by subtracting $b$
twists between the bands $m$ times and subtracting $p-b$ twists
between the bands $n$ times, we may subtract a single twist between
$A$ and $B$.
\end{proof}

\begin{proof}[Proof of \fullref{P:Bunlink}]
Assume first that $\abs{A}=\abs{B}=0$ and that $A$ and $B$ link with
some other bands of non-zero index. Assume also that $A$ and $B$ are
the leftmost and the second-to-leftmost bands of $K$ (this can be
achieved by rotating the band projection of the knot around the
disc). The following sequence allows us to reduce the genus of $K$
by surgery in $\ker{\rho}$.
$$
\begin{minipage}{100pt}
\labellist\small
\pinlabel {$\wbar{C}$} [br] at 178 226
\pinlabel {$\wbar{D}$} [b] at 315 206
\pinlabel {$a$} [t] at 136 38
\pinlabel {$a$} [b] at 64 86
\pinlabel {$a$} [b] at 231 86
\endlabellist
\includegraphics[width=100pt]{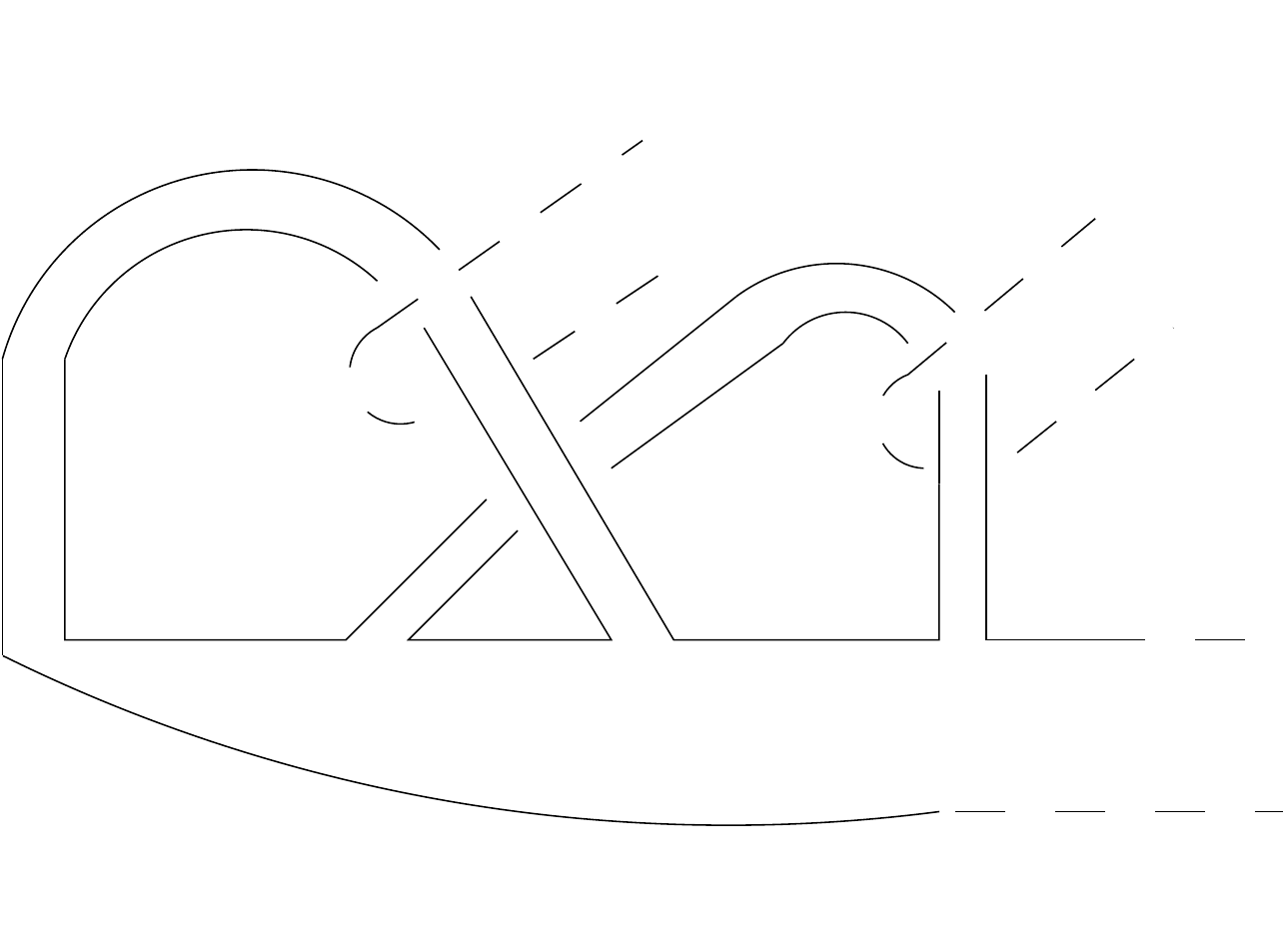}
\end{minipage}
\quad\seq\quad
\begin{minipage}{100pt}
\labellist\small
\pinlabel {$a$} [t] at 136 38
\pinlabel {$a$} [b] at 64 86
\pinlabel {$a$} [b] at 231 86
\pinlabel {$\wbar{C}$} [br] at 181 229
\endlabellist
\includegraphics[width=100pt]{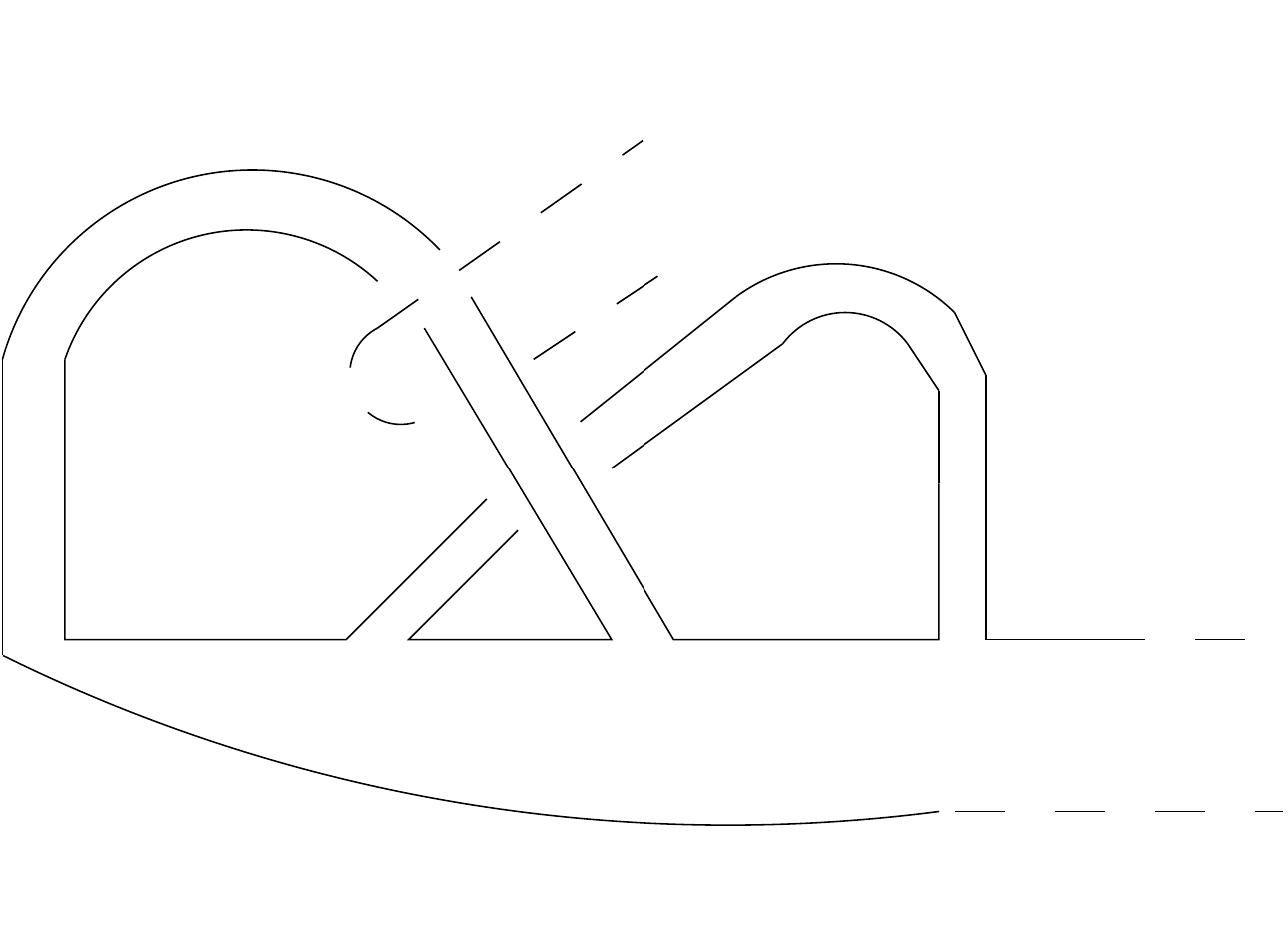}
\end{minipage}
$$
The dotted line denotes any linkage with any other bands,
collectively called $\bar{C}$ and $\bar{D}$. To get from the first
diagram to the second, notice that conjugation of $a$ by $\bar{C}$
and by $\bar{D}$ gives $a$. Thus the component which loops once
around $\bar{D}$ and once around $B$ is in $\ker\rho$, and surgery
on this component followed by $RR$ to undo the twist we get in $B$
allows us to unlink $B$ and $\bar{D}$ and then to isotopy $\bar{D}$
out of the picture. The resulting diagram represents a $p$--coloured
knot of lower genus that the one which we started with, because we
may slide $\bar{C}$ out of the picture by ambient isotopy (or
alternatively repeat the above surgery argument to unlink $\bar{C}$
from $A$) and then eliminate $A$ and $B$ by ambient isotopy.

We are left now only with the case $\abs{A}\neq 0$ but $\abs{B}=0$.
In this case we choose a new disc-band presentation for $K$. The
procedure we carry out below is equivalent to sliding $A$ over $B$.
In the notation at the beginning of the section, we take $A=B_{1}$
and $B=B_{2}$. Choose $D'^{2},B'_{1},\ldots,B'_{2n}$ such that

\begin{enumerate}
\item $S=D'^{2}\cup B'_{1}\cup\cdots \cup B'_{2n}$ where $D'^{2}$ and each
$B'_{j}$ is a disc, and $B'_{i}=B_{i}$ for $i=1$ and for $3\leq
i\leq 2n$;
\item $\partial
B'^{2}=\alpha'_{2}\gamma'_{2}\beta'_{2}(\tilde{\gamma}'_{2})^{-1}$
where $\alpha'_{2}\subset \delta_{2}$, $\beta'_{2}\subset
\delta_{2n}$, and $\gamma'_{2}$ and $\tilde{\gamma}'_{2}$ are
non-intersecting simple curves in $D^{2}$ connecting the leftmost
endpoint of $\alpha'_{2}$ with the leftmost endpoint of $\beta'_{2}$
and the rightmost endpoint of $\alpha'_{2}$ with the rightmost
endpoint of $\beta'_{2}$ correspondingly.
\end{enumerate}

Note that the above data also uniquely characterizes $D'^{2}$. Let
$B'_{2}$ be renamed $C$. We get a new disc-band presentation of the
Seifert surface, in which both $\abs{A}$ and $\abs{C}$ are nonzero.
See \fullref{Fi:genseif1} (note that $A$ and $B$ need not be in leftmost
position).

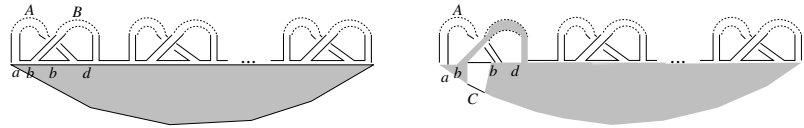
\begin{figure}\label{Fi:genseif1}
\centering
\labellist\tiny
\hair=1pt
\pinlabel {$A$} [b] at 40 226
\pinlabel {$B$} [b] at 140 221
\pinlabel {$a$} [t] at 11 121
\pinlabel {$b$} [t] at 41 127
\pinlabel {$b$} [t] at 89 127
\pinlabel {$d$} [t] at 160 127
\endlabellist
\includegraphics[width=1.9in]{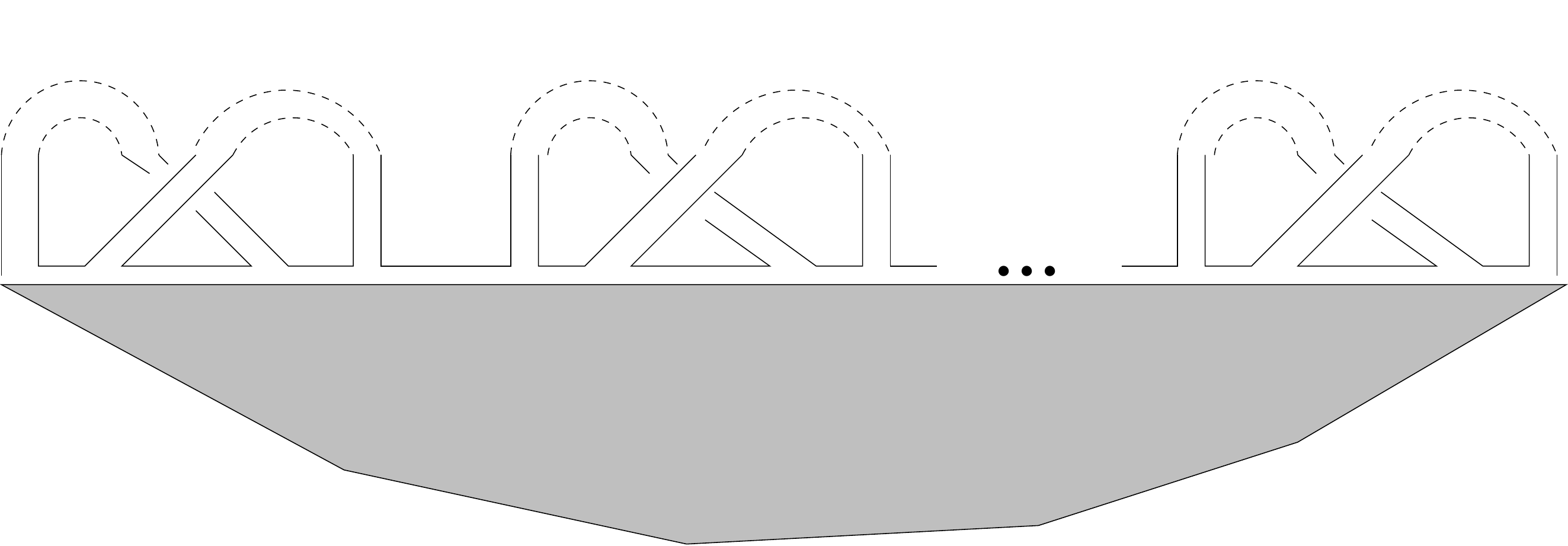}
\qquad
\labellist\tiny
\hair=1pt
\pinlabel {$A$} [b] at 40 226
\pinlabel {$C$} [t] at 70 72
\pinlabel {$a$} [tr] at 25 110
\pinlabel {$b$} [t] at 39 126
\pinlabel {$b$} [t] at 114 131
\pinlabel {$d$} [t] at 160 131
\endlabellist
\includegraphics[width=1.9in]{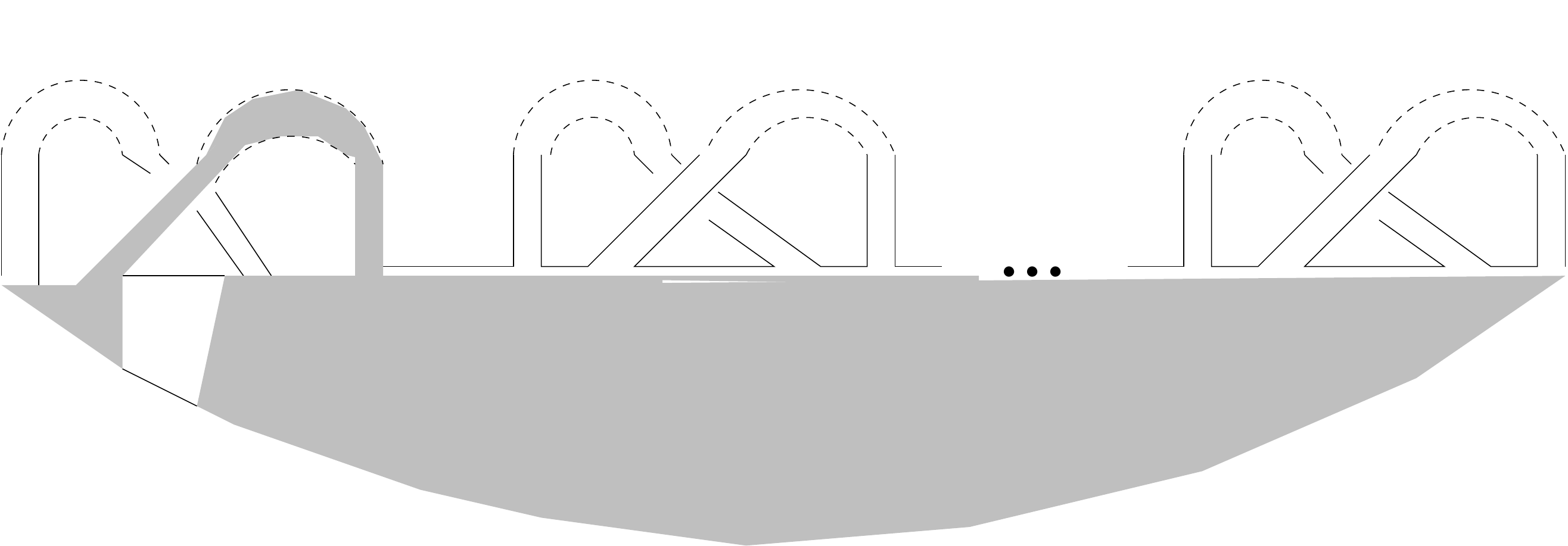}
\caption{Two choices of cut systems for a Seifert surface when
$\abs{B}=0$ but $\abs{A}\neq0$}
\end{figure}

Repeating for all bands whose index is zero, we obtain a knot whose
bands all have non-zero index, reducing the problem to a previously
considered case.

That the $K_{i}$ are $p$--colourable follows from examining the
$p$--colouring of $K_{i}$ as a connect summand of $K$. Assume that
$K_{i}$ has two bands $A$ and $B$ with $A$ to the left of $B$.
Assume $A$ has $m$ twists and $B$ has $n$ twists. Then if $a$ is $0$
and if $b$ is $1$ (notation as in \fullref{D:index}) then $c$
is $2m$ and $d$ is $2m+1$. Labeling the colours of the corresponding
strands of $B$ left to right by $a', b', c', d'$ we have $a'=2m$,
$b'=1$, and $c'=2m+1$. This forces $d'$ to be $0$, proving that
$K_{i}$ is indeed $p$--colourable with colouring $\rho_{i}$ induced
by $\rho$.

If any of the connect-summands are trivially coloured, they may be
untied by the $RR$--move, and thus we may assume that all
connect-summands come equipped with a non-trivial $p$--colouring.
\end{proof}

\section{Three colours}\label{S:3Col}

In this section we prove the first part of \fullref{T:3case} for
$p=3$. We recall that this states that there are at most $3$
equivalence classes of $3$--coloured knots modulo surgery in
$\ker\rho$, and that these equivalence classes are represented by
connect-sums of $n$ left-hand trefoils for $n=1,2,3$. The remaining
part of the theorem for $p=3$, that there are at least $3$ such
equivalence classes, is proved in \fullref{S:noteq}.

\begin{proof}[Proof of Theorem 1 (first part, $p=3$ case)]
Let $(K,\rho)$ be a $3$--coloured knot. By \fullref{P:Bunlink} we have
$$
(K,\rho)\seq S(m_{1},n_{1})\hash S(m_{2},n_{2})\hash\cdots\hash
S(m_{g},n_{g})
$$
\noindent where the $3$--colourings of the connect-summands are
induced by $\rho$, and the connect-sum is with respect to some
orientation of the connect-summands. By the $LT$ move (\fullref{N:pTW}) we may $3$--reduce the number of twists in the bands of
each connect-summand, allowing us to express $K$ as a connect-sum of
$3$--coloured $S(m,n)$ knots with $-1\leq n,m \leq 1$. By
\fullref{P:prop9}, there are two such $3$--colourable
knots--- the left-hand trefoil $S(1,1)$ which we denote $3_{1}$, and
the right-hand trefoil $S(-1,-1)$ which we denote
$\overline{3_{1}}$--- and thus $(K,\rho)$ is equivalent to a
connect-sum of $3_{1}$ knots and $\overline{3_{1}}$ knots modulo
surgery in $\ker \rho$.

By $LT$ on the three half-twists of the trefoil, we have
$$
3_{1}=
\begin{minipage}{50pt}
\labellist\small
\pinlabel {$0$} [b] at 106 253
\pinlabel {$2$} [b] at 73 163
\pinlabel {$1$} [b] at 145 163
\endlabellist
\includegraphics[width=50pt]{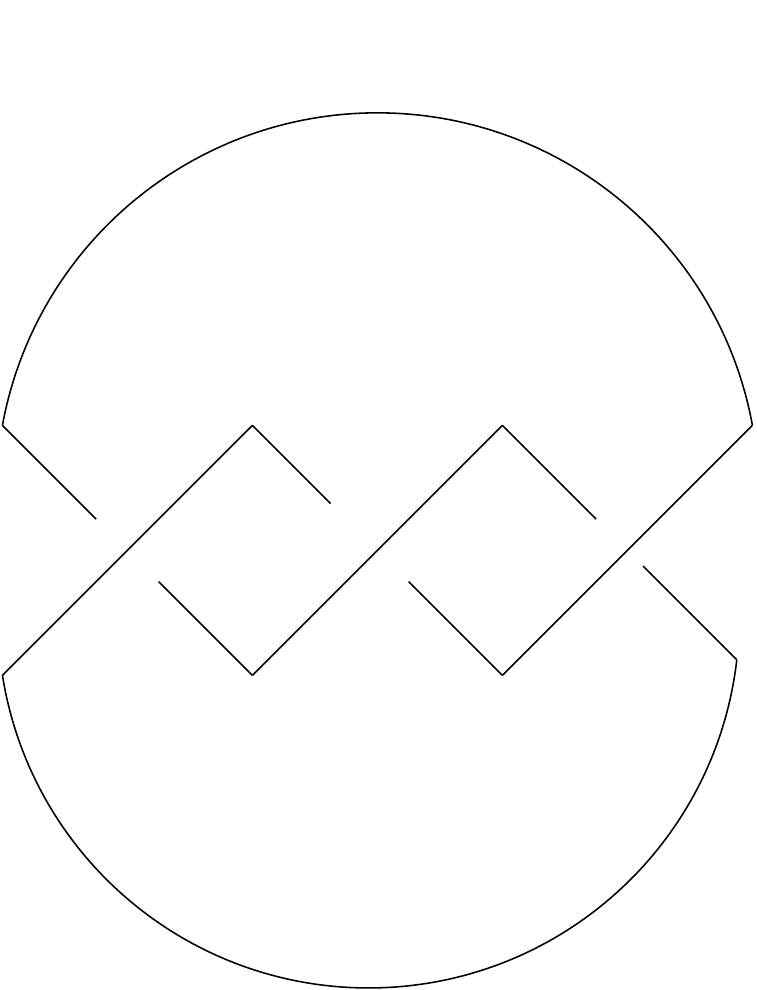}
\end{minipage}
\quad\overset{LT}{\sim}\quad
\begin{minipage}{50pt}
\labellist\small
\pinlabel {$0$} [b] at 106 253
\pinlabel {$1$} [b] at 73 163
\pinlabel {$2$} [b] at 145 163
\endlabellist
\includegraphics[width=50pt]{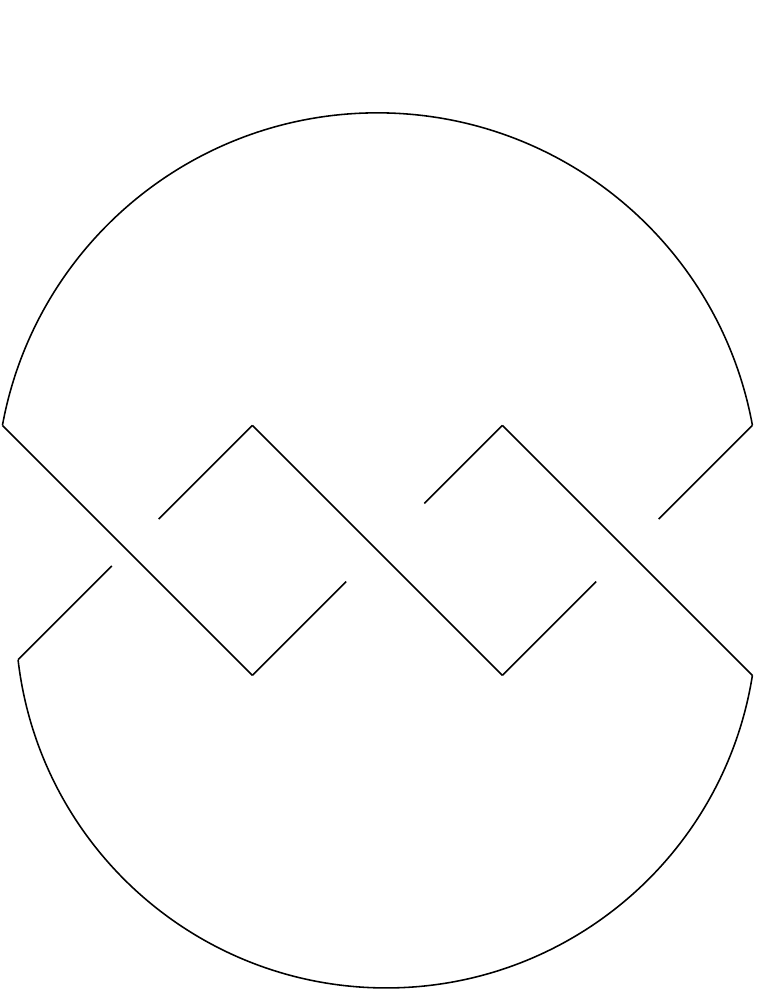}
\end{minipage}
\quad\hash\quad
\begin{minipage}{50pt}
\labellist\small
\pinlabel {$0$} [b] at 106 253
\pinlabel {$1$} [b] at 73 163
\pinlabel {$2$} [b] at 145 163
\endlabellist
\includegraphics[width=50pt]{\figdir/l23torus}
\end{minipage}
\quad=\overline{3_{1}}\hash\overline{3_{1}}
$$
By reflection we also have $\overline{3_{1}}\seq 3_{1}\hash 3_{1}$.
This shows that $(K,\rho)$ is equivalent to a connect-sum of
left-hand trefoils modulo surgery in $\ker \rho$. The number of
connect-summands may be $3$--reduced as follows:
$$
(3_{1})^{\hash 4}\seq(\overline{3_{1}})^{\hash 2}\seq 3_{1}
$$
Thus $(K,\rho)$ is equivalent to $3_{1}$ or $(3_{1})^{\hash 2}$ or
$(3_{1})^{\hash 3}$.
\end{proof}

\section{Five colours}\label{S:5Col}

In this section we prove the first part of \fullref{T:3case} for
$p=5$. This says that there are at most $5$ equivalence classes of
$5$--coloured knots modulo surgery in $\ker\rho$, and that these
equivalence classes are represented by connect-sums of $n$ left-hand
$5_{1}$ knots for $n=1,2,3,3,4,5$. The remaining part of the theorem
for $p=5$, that there are at least $5$ such equivalence classes, is
proved in \fullref{S:noteq}.

We give first the proof, and then fill in the necessary lemmata in
the coming sections.

\begin{proof}[Proof of Theorem 1 (first part, $p=5$ case)]
Let $(K,\rho)$ be a $5$--coloured knot. By \fullref{P:Bunlink} know that
$$
(K,\rho)\seq S(m_{1},n_{1})\hash S(m_{2},n_{2})\hash\cdots\hash
S(m_{g},n_{g})
$$
\noindent where the $5$--colourings of the connect-summands are
induced by $\rho$, and the connect-sum is with respect to some
orientation of the connect-summands.

We may $5$--reduce the number of twists in the bands of each
$S(m,n)$ connect-summand by the $10$--move (\fullref{P:10move}), allowing us to express $K$ as a connect-sum of
$5_{1}$ knots with knots of the form $S(m,n)$ with $-2\leq n,m \leq
2$ whose unique $5$--colouring is induced by $\rho$. A list of such
knots is derived from \fullref{P:prop9}--- these are the
knots $S(2,2)$ and $S(-1,1)$. These knots are each reduced to a
connect-sum of $5_{1}$ knots by \fullref{P:Knotequiv}. We
thus obtain:
$$
(K,\rho)\seq (5_{1},\rho_{1})^{\hash g_{1}}\hash
(5_{1},\rho_{2})^{\hash g_{2}}\hash
(\overline{5_{1}},\overline{\rho_{1}})^{\hash g_{3}}\hash
(\overline{5_{1}},\overline{\rho_{2}})^{\hash g_{4}}
$$
\noindent where $(5_{1},\rho_{1})$, $(5_{1},\rho_{2})$,
$(\overline{5_{1}},\overline{\rho_{1}})$, and
$(\overline{5_{1}},\overline{\rho_{2}})$ denote the $5_{1}$ knots
with different handedness and colourings which are not ambient
isotopic. See \fullref{Fi:solomon}.

\begin{figure}\label{Fi:solomon}
\centering
\labellist\tiny
\pinlabel {$0$} [r] at 35 258
\pinlabel {$1$} [b] at 93 216
\pinlabel {$2$} [b] at 165 216
\pinlabel {$3$} [b] at 237 216
\pinlabel {$4$} [b] at 309 216
\endlabellist
\includegraphics[width=0.8in]{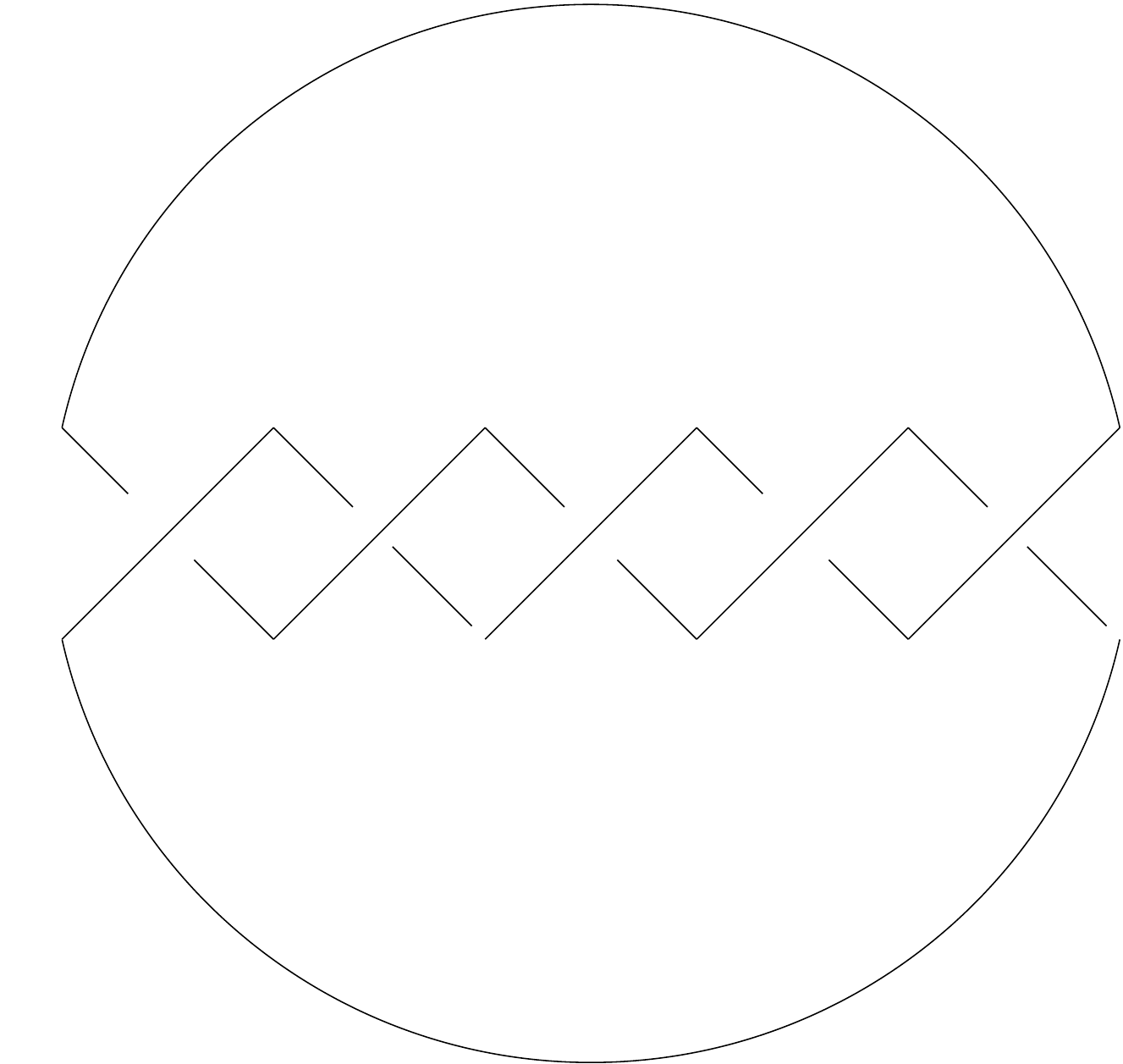}
\qquad
\labellist\tiny
\pinlabel {$0$} [r] at 35 258
\pinlabel {$2$} [b] at 93 216
\pinlabel {$4$} [b] at 165 216
\pinlabel {$1$} [b] at 237 216
\pinlabel {$3$} [b] at 309 216
\endlabellist
\includegraphics[width=0.8in]{\figdir/r25torus}
\qquad
\labellist\tiny
\pinlabel {$0$} [r] at 35 258
\pinlabel {$1$} [b] at 93 216
\pinlabel {$2$} [b] at 165 216
\pinlabel {$3$} [b] at 237 216
\pinlabel {$4$} [b] at 309 216
\endlabellist
\includegraphics[width=0.8in]{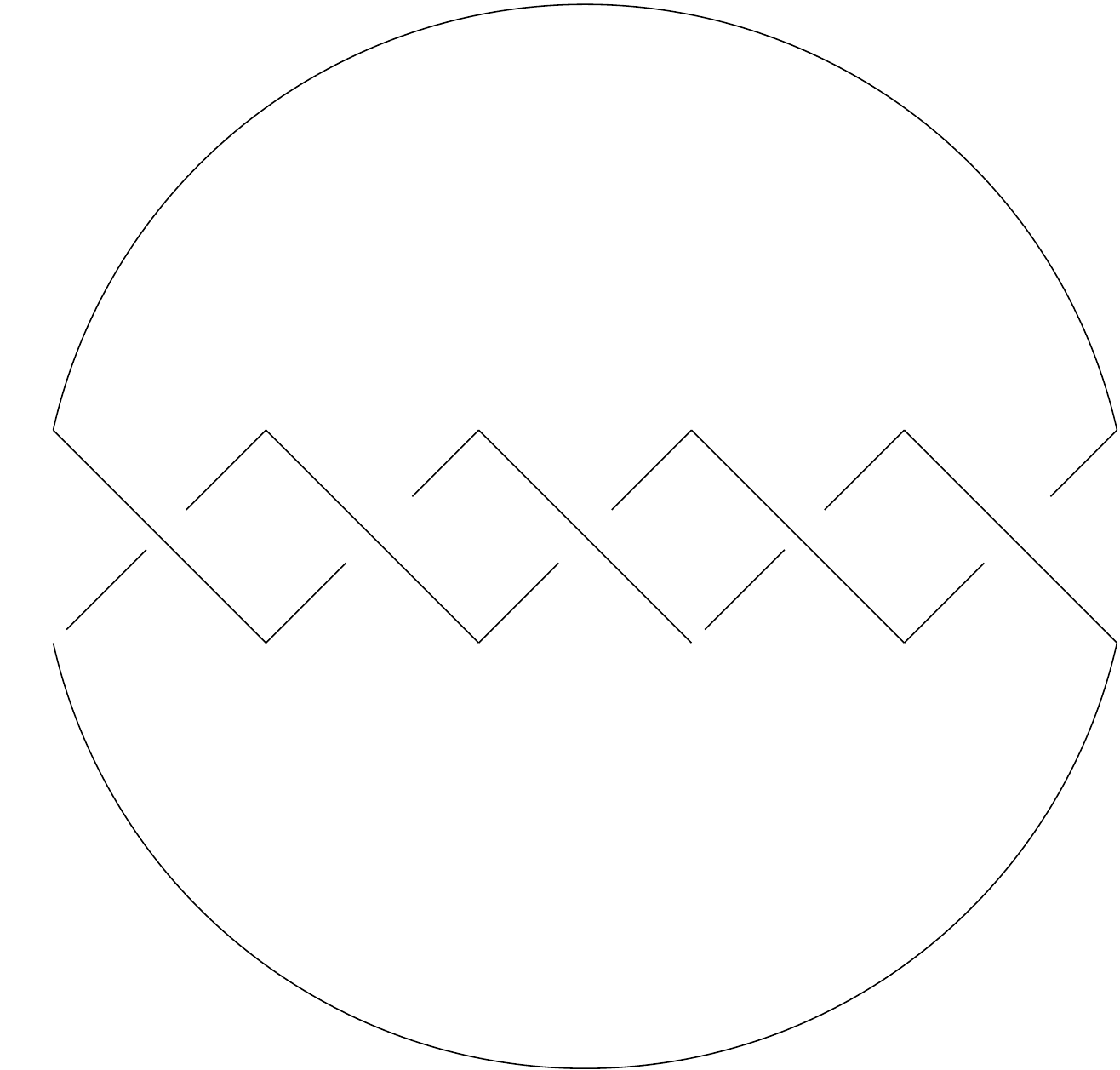}
\qquad
\labellist\tiny
\pinlabel {$0$} [r] at 35 258
\pinlabel {$2$} [b] at 93 216
\pinlabel {$4$} [b] at 165 216
\pinlabel {$1$} [b] at 237 216
\pinlabel {$3$} [b] at 309 216
\endlabellist
\includegraphics[width=0.8in]{\figdir/l25torus}
 \caption{From left to right: $(5_{1},\rho_{1})$, $(5_{1},\rho_{2})$,
 $(\overline{5_{1}},\overline{\rho_{1}})$, and
 $(\overline{5_{1}},\overline{\rho_{2}})$}
\end{figure}

The theorem now follows from \fullref{P:cs}, which tells us
that
$$
(K,\rho)\seq (5_{1},\rho_{1})^{\hash n}\qquad\mbox{for some
$n=1,2,3,4,5$.}
\proved
$$
\end{proof}

The key steps of the proof are thus:

\begin{itemize}
\item To reduce the number of twists in the bands via the $10$--move.
\item To reduce $S(m,n)$ to a connect-sum of $5_{1}$ knots for $-2\leq m,n \leq 2$.
\item To reduce a connect-sum of $5_{1}$ knots to a connect-sum of between one and five $(5_{1},\rho_{1})$ knots.
\end{itemize}

The first stage is carried out in \fullref{SS:10move}, the
second in \fullref{SS:Smnreduce}, and the third in \fullref{SS:51alg}. \fullref{SS:prelim} is concerned with
preliminary lemmata for use in these sections.

\subsection{Preliminary lemmata}\label{SS:prelim}

In the present section we collect together several preliminary
lemmata whose main purpose is to relate $5$--coloured $5_{1}$ knots
to knots and tangles which are technically simpler to compute
with.

We begin with some notation which is used only in this section. We
work in the category of $(2,2)$--tangles rotated by $\frac{\pi}{2}$,
so that the composition of two tangles is given by stacking from
left to right
$$
T_{1}T_{2}\ass\framebox{\makebox[\totalheight]{$T_{1}$}}\framebox{\makebox[\totalheight]{$T_{2}$}}
$$
Let $T$ be a $5$--coloured $(2,2)$--tangle. For $a$ the
$5$--colouring of the left endpoint of the top strand of $T$ and $b$
the $5$--colouring of the left endpoint of the top strand of $T$,
let $\underline{a-b}\in\set{-2,-1,0,1,2}$ represent $a-b$. The
absolute value of the difference $d=\abs{\underline{a-b}}$ is
identical to the absolute value of the difference of the
$5$--colourings of the right endpoints, defined the same way.

We now define several `basic' tangles. Let $\mathbbm{1}$ denote the
trivial $(2,2)$--tangle and
\begin{align*}
l&=
\quad
\begin{minipage}{75pt}
\rotatebox{270}{\includegraphics[width=42pt]{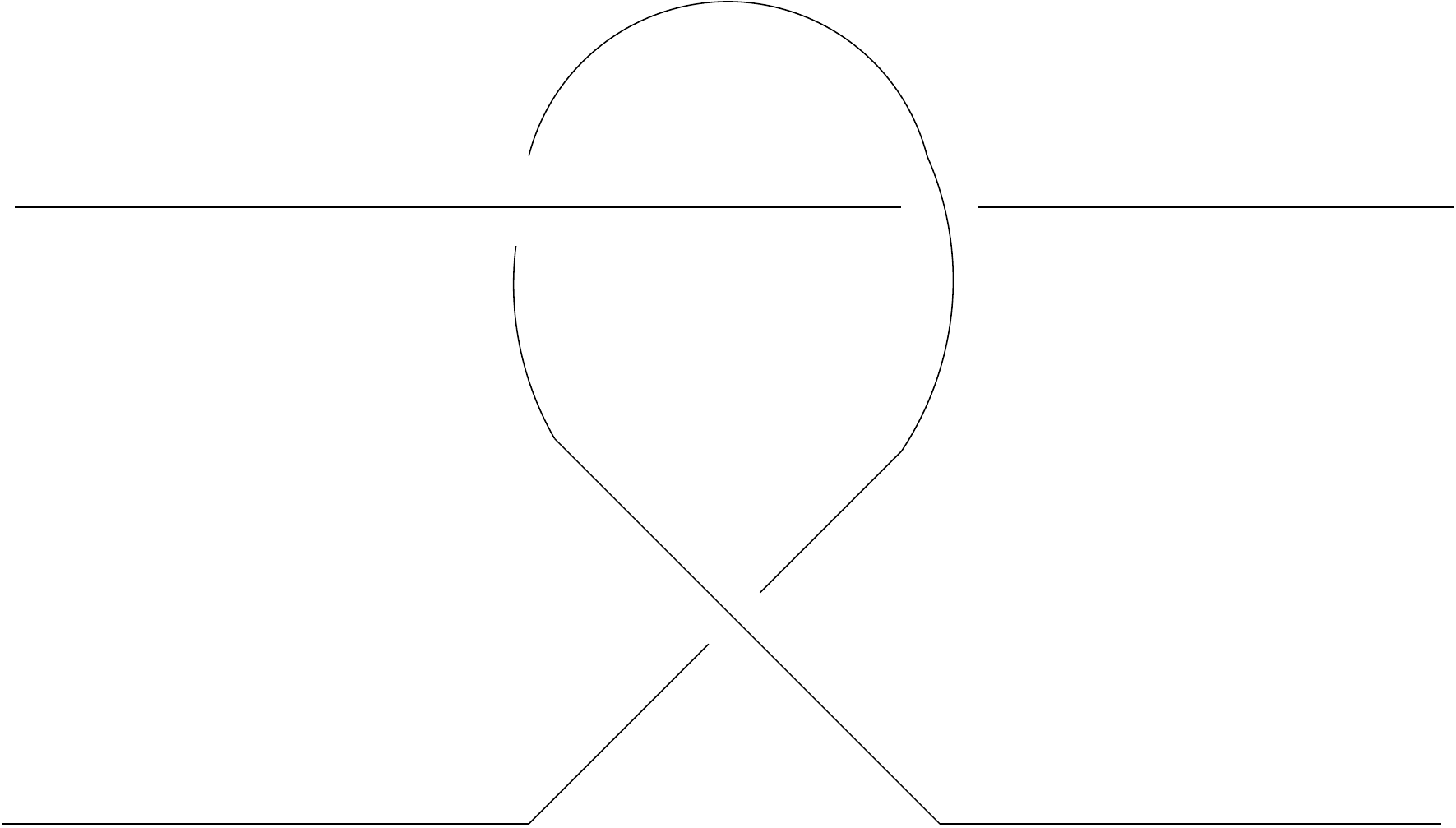}}
\end{minipage}
& r&=
\quad
\begin{minipage}{75pt}
\rotatebox{270}{\includegraphics[width=42pt]{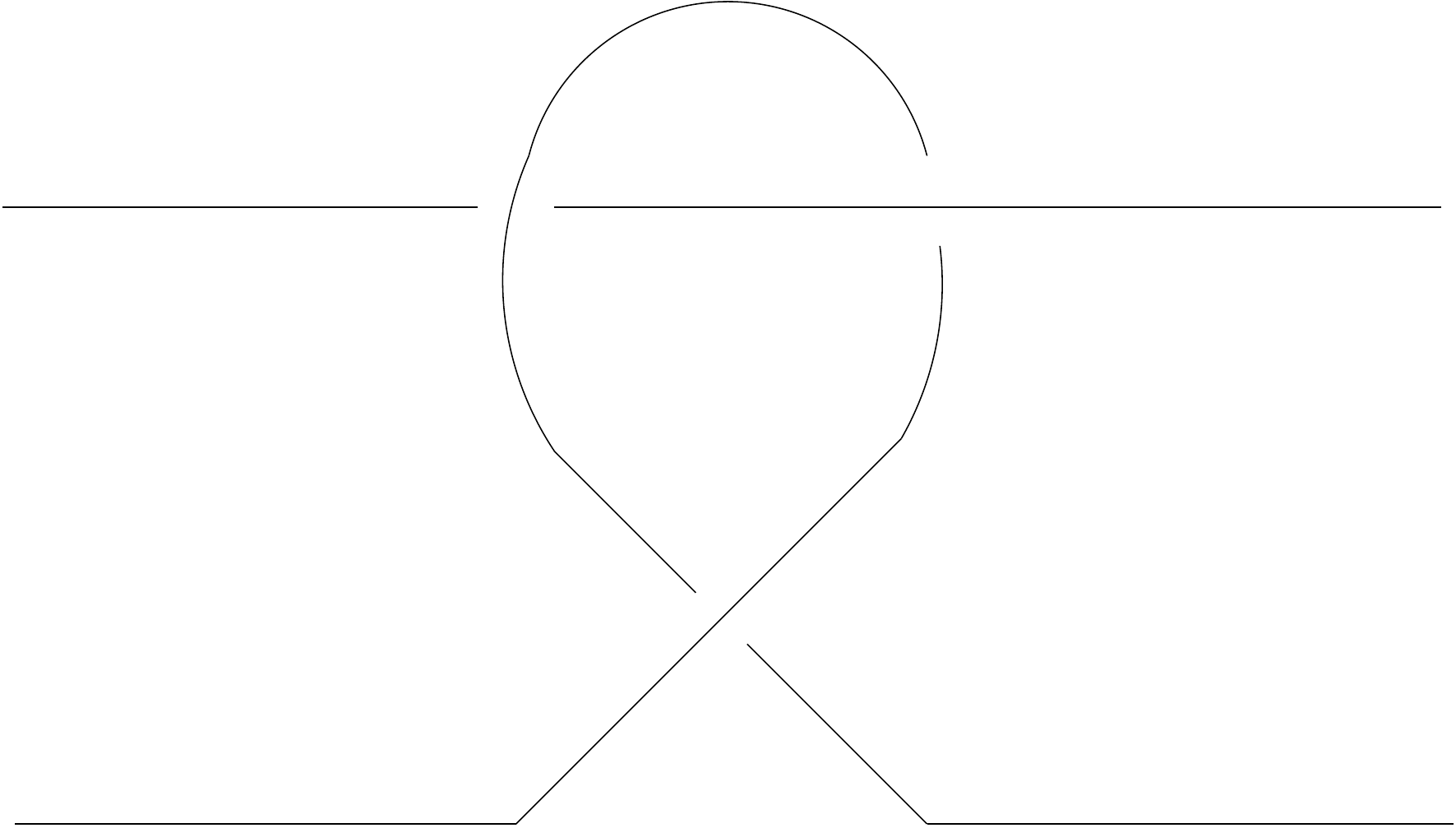}}
\end{minipage}\\
\tau&= \hspace{10pt}\begin{minipage}{84pt}
    \begin{picture}(84,42)
        \put(0,27){}\put(0,8){}
        \put(42,0){\line(-1,1){42}}
        \qbezier(0,0)(0,0)(14,14)
        \qbezier(28,28)(42,42)(42,42)
        \put(42,42){\line(1,1){0}}
        \put(84,0){\line(-1,1){42}}
        \qbezier(42,0)(42,0)(56,14)
        \qbezier(70,28)(84,42)(84,42)
        \put(84,42){\line(1,1){0}}
    \end{picture}
    \end{minipage}
& \tau^{-1}&= \hspace{10pt}\begin{minipage}{84 pt}
    \begin{picture}(84,42)
        \put(0,27){}\put(0,8){}
        \put(0,0){\line(1,1){42}}
        \qbezier(42,0)(42,0)(28,14)
        \qbezier(14,28)(0,42)(0,42)
        \put(0,42){\line(-1,1){0}}
        \put(42,0){\line(1,1){42}}
        \qbezier(84,0)(84,0)(70,14)
        \qbezier(56,28)(42,42)(42,42)
        \put(42,42){\line(-1,1){0}}
    \end{picture}
\end{minipage}
\end{align*}
Let $(5_{1},\rho_{0})$ and $(\overline{5_{1}},\overline{\rho_{0}})$
denote the trivially coloured left-hand and right-hand
$5_{1}$--knots correspondingly.

We prove a number of equivalences between the tangles defined above.
Define a connect-sum of a $p$--coloured knot with a $p$--coloured
$(2,2)$--tangle to be the connect-sum of that knot with the bottom
strand of the tangle (this is well-defined by \fullref{L:consum}).

\begin{lem}\label{L:taus}
Let $d=\abs{\underline{a-b}}$  denote the absolute value of the
difference between the $5$--colouring of the upper strand and the
lower strand as before. Then we have:
\begin{enumerate}
\item $l^{2}\seq r\cdot \tau$ and $r^{2}\seq l\cdot
\tau^{-1}$. \label{P:tau1}
\item $l^{2}\cdot\tau\seq \mathbbm{1}\hash (5_{1},\rho_{d})$ and $r^{2}\cdot
\tau^{-1}\seq \mathbbm{1}\hash
(\overline{5_{1}},\overline{\rho_{d}})$.\label{P:tau2}
\end{enumerate}
\end{lem}

\begin{proof}Part (1)\qua
\begin{multline*}
\quad\begin{minipage}{80pt}
\labellist\small
\pinlabel {$0$} [r] at 23 10
\pinlabel {$1$} [r] at 23 118
\pinlabel {$3$} [b] at 130 155
\pinlabel {$2$} [b] at 235 155
\pinlabel {$4$} [b] at 330 118
\pinlabel {$3$} [b] at 330 10
\pinlabel {$4$} [b] at 170 10
\pinlabel {$0$} [b] at 180 118
\endlabellist
\includegraphics[width=80pt]{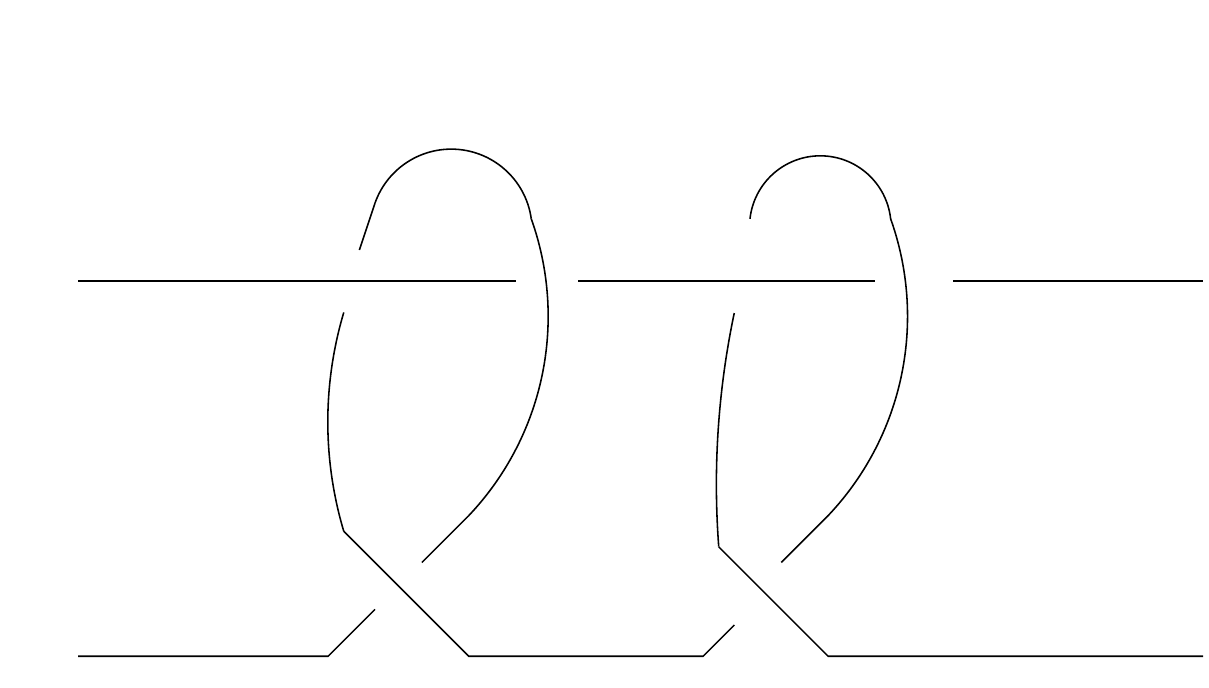}
\end{minipage}
\quad = \quad
\begin{minipage}{80pt}
\labellist\tiny
\hair=1pt
\pinlabel {$1$} [b] at 25 365
\pinlabel {$0$} [t] at 25 55
\pinlabel {$4$} [b] at 710 380
\pinlabel {$3$} [t] at 710 55
\pinlabel {$4$} [t] at 250 45
\pinlabel {$4$} [br] at 245 295
\endlabellist
\includegraphics[width=80pt]{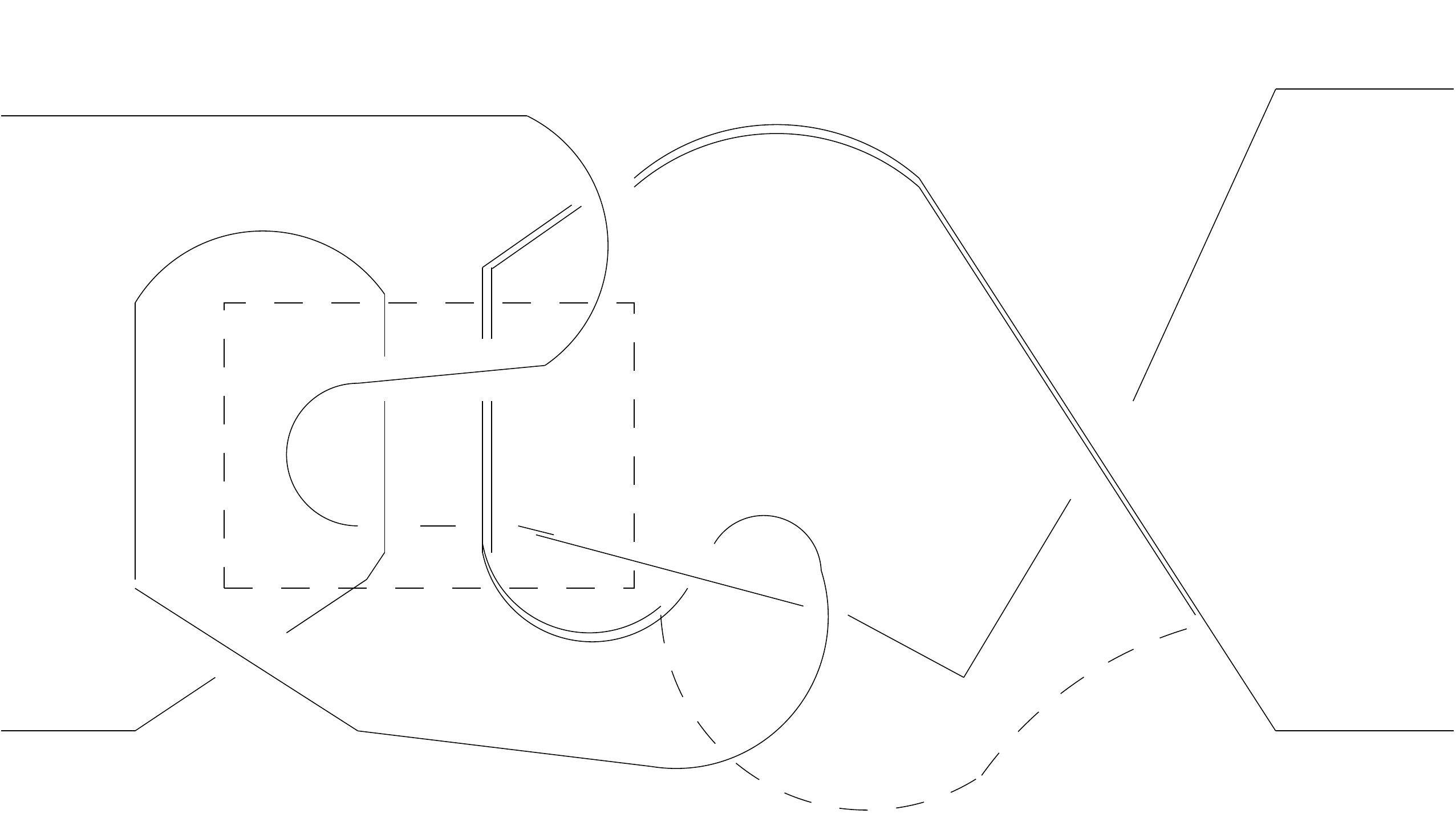}
\end{minipage}
\quad\overset{R2G}{\sim} \\
\begin{minipage}{80pt}
\labellist\tiny
\pinlabel {$1$} [b] at 20 360
\pinlabel {$0$} [b] at 20 50
\pinlabel {$4$} [b] at 725 375
\pinlabel {$3$} [b] at 725 50
\pinlabel {$4$} [tl] at 390 50
\endlabellist
\includegraphics[width=80pt]{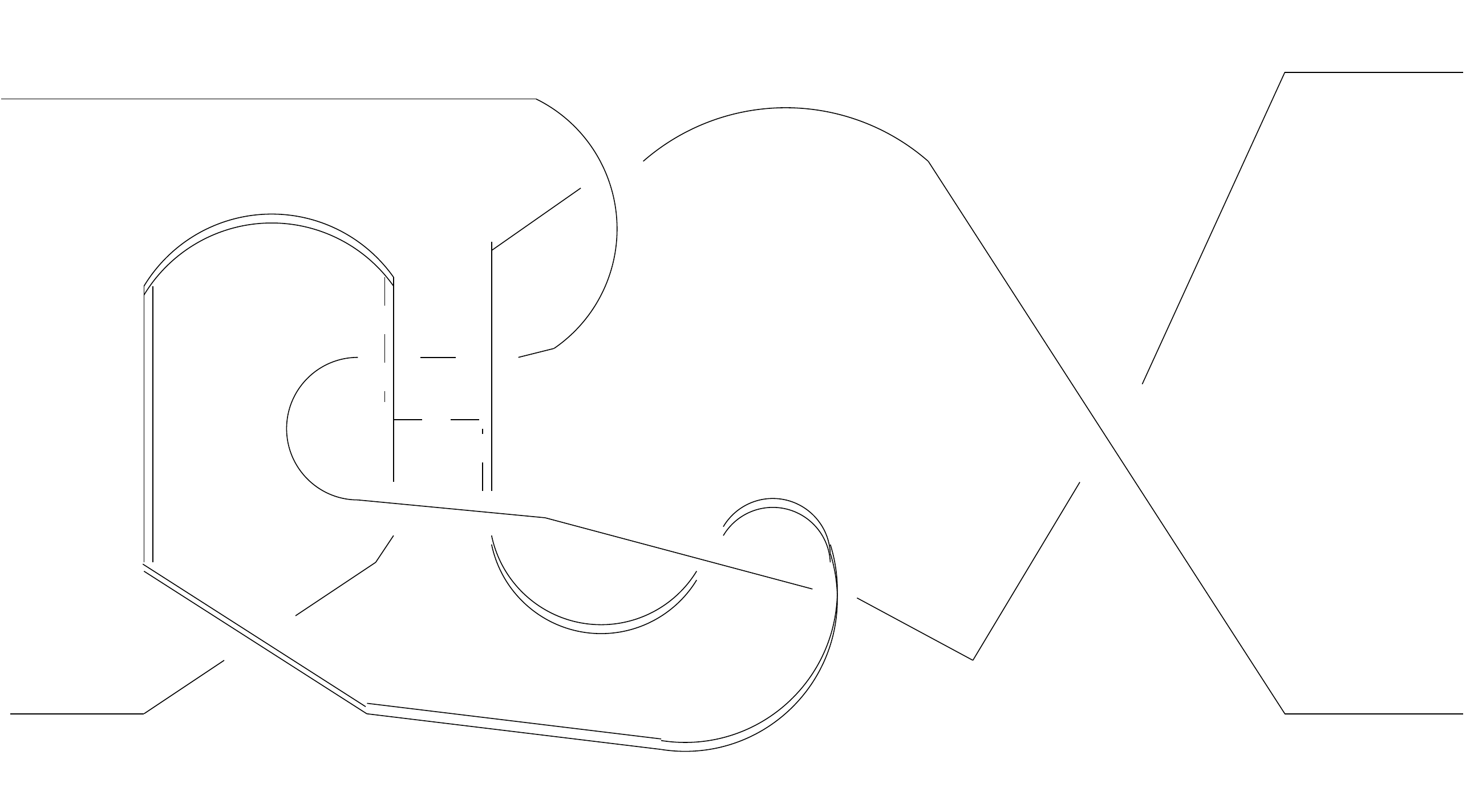}
\end{minipage}
\quad=\quad
\begin{minipage}{80pt}
\labellist\tiny
\pinlabel {$1$} [b] at 15 280
\pinlabel {$0$} [b] at 15 65
\pinlabel {$4$} [b] at 520 280
\pinlabel {$3$} [b] at 520 65
\pinlabel {$2$} [tl] at 260 110
\endlabellist
\includegraphics[width=80pt]{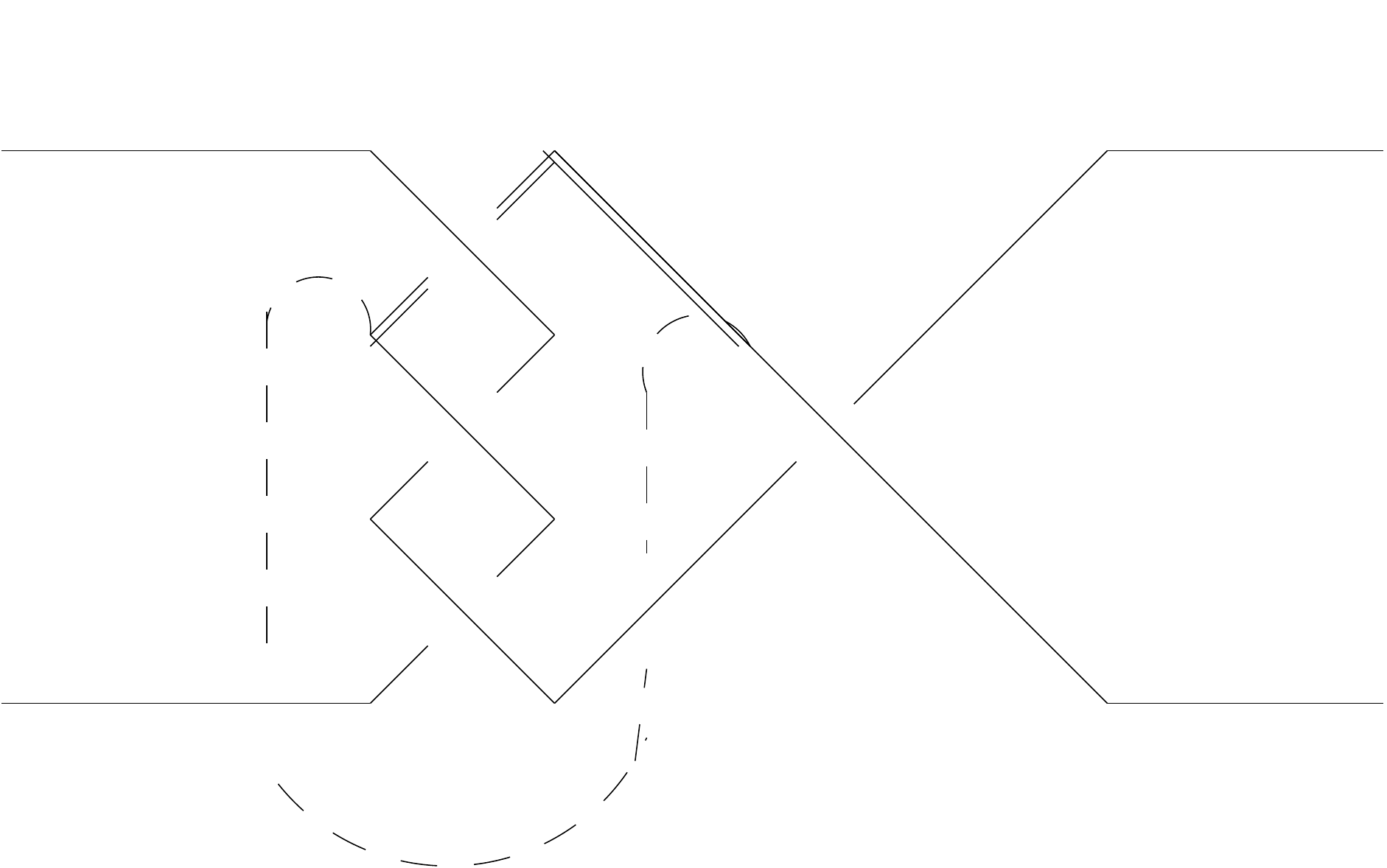}
\end{minipage}
\quad=\quad
\begin{minipage}{80pt}
\labellist\small
\pinlabel {$1$} [t] at 20 140
\pinlabel {$0$} [b] at 20 0
\pinlabel {$4$} [br] at 115 150
\endlabellist
\includegraphics[width=80pt]{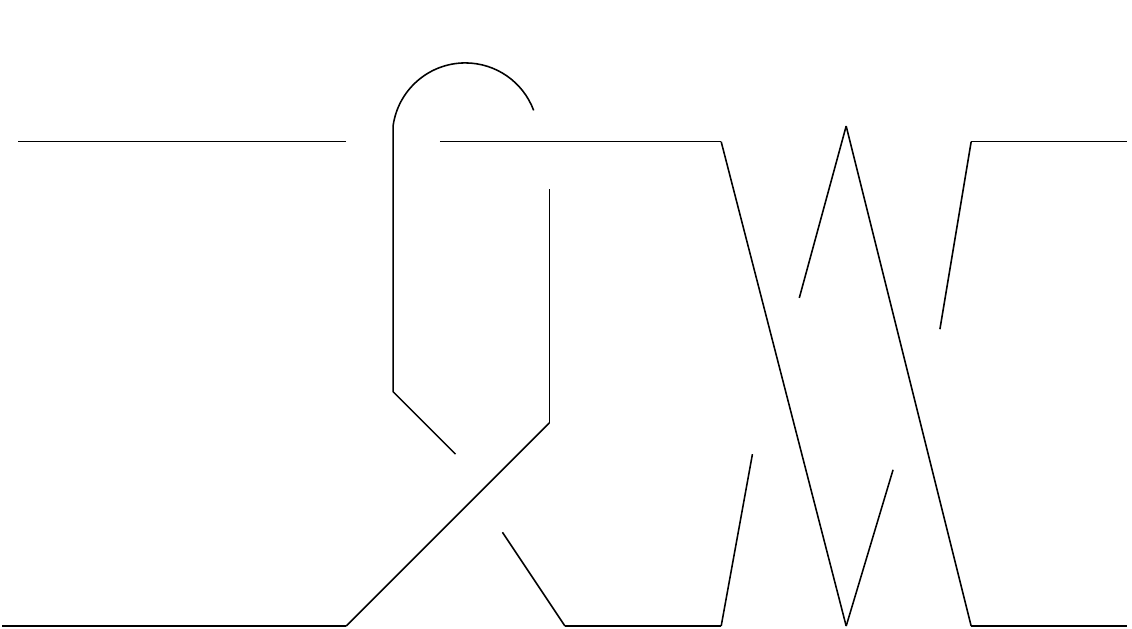}
\end{minipage}
\end{multline*}
To get from the second diagram to the first diagram, from the fourth
diagram to the fifth diagram, and from the fifth diagram to the
final diagram, isotopy the thick line-segment to the position
specified by the dotted line. To get from the second diagram to the
third, use $R2G$ inside the dotted box on the arc coloured $0$
through the two arcs coloured $3$.

The $d=0$ case is trivial.  For $d=2$ instead of $d=1$, in the same way
we obtain
\begin{align*}
\begin{minipage}{80pt}
\labellist\small
\pinlabel {$0$} [r] at 23 10
\pinlabel {$2$} [r] at 23 118
\pinlabel {$1$} [b] at 130 155
\pinlabel {$4$} [b] at 235 155
\pinlabel {$3$} [b] at 330 118
\pinlabel {$1$} [b] at 330 10
\pinlabel {$3$} [b] at 170 10
\pinlabel {$0$} [b] at 180 118
\endlabellist
\includegraphics[width=80pt]{\figdir/pure2loop}
\end{minipage} \quad\seq\quad
\begin{minipage}{80pt}
\labellist\small
\pinlabel {$2$} [t] at 20 140
\pinlabel {$0$} [b] at 20 0
\pinlabel {$3$} [br] at 115 150
\endlabellist
\includegraphics[width=80pt]{\figdir/linkedloop}
\end{minipage}
\end{align*}
The proof of $r^{2}\seq l\cdot \tau^{-1}$ is the same up to
reflection.

Part (2)\qua
\begin{multline*}
\begin{minipage}{80pt}
\labellist\tiny
\hair=1pt
\pinlabel {$0$} [b] at 23 10
\pinlabel {$1$} [b] at 23 118
\pinlabel {$3$} [b] at 130 155
\pinlabel {$2$} [b] at 235 155
\pinlabel {$4$} [b] at 296 118
\pinlabel {$3$} [b] at 272 14
\pinlabel {$4$} [b] at 540 118
\pinlabel {$0$} [b] at 540 10
\pinlabel {$1$} [b] at 170 10
\pinlabel {$0$} [b] at 180 118
\endlabellist
\includegraphics[width=80pt]{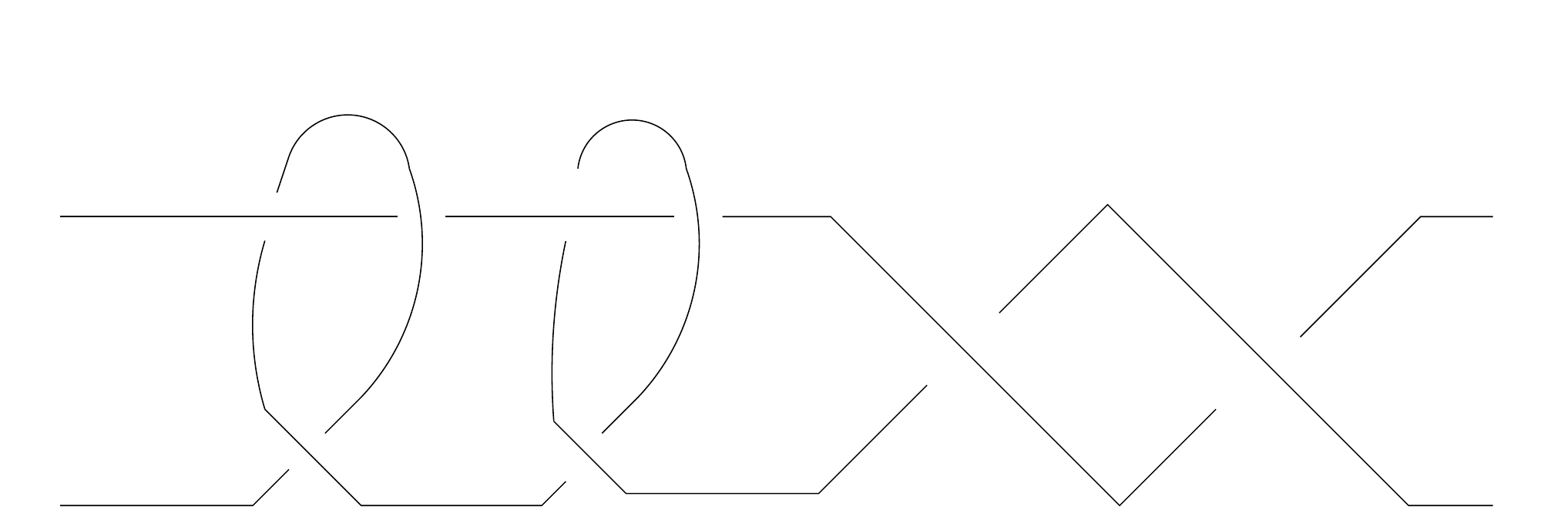}%
\end{minipage}
\quad\overset{R2G}{\sim}\quad
\begin{minipage}{100pt}
\labellist\tiny
\pinlabel {$0$} [r] at 23 10
\pinlabel {$1$} [r] at 23 118
\pinlabel {$3$} [b] at 64 136
\pinlabel {$2$} [l] at 255 140
\pinlabel {$0$} [bl] at 115 118
\pinlabel {$3$} [tl] at 440 30
\pinlabel {$1$} [l] at 557 118
\pinlabel {$0$} [l] at 557 10
\endlabellist
\includegraphics[width=100pt]{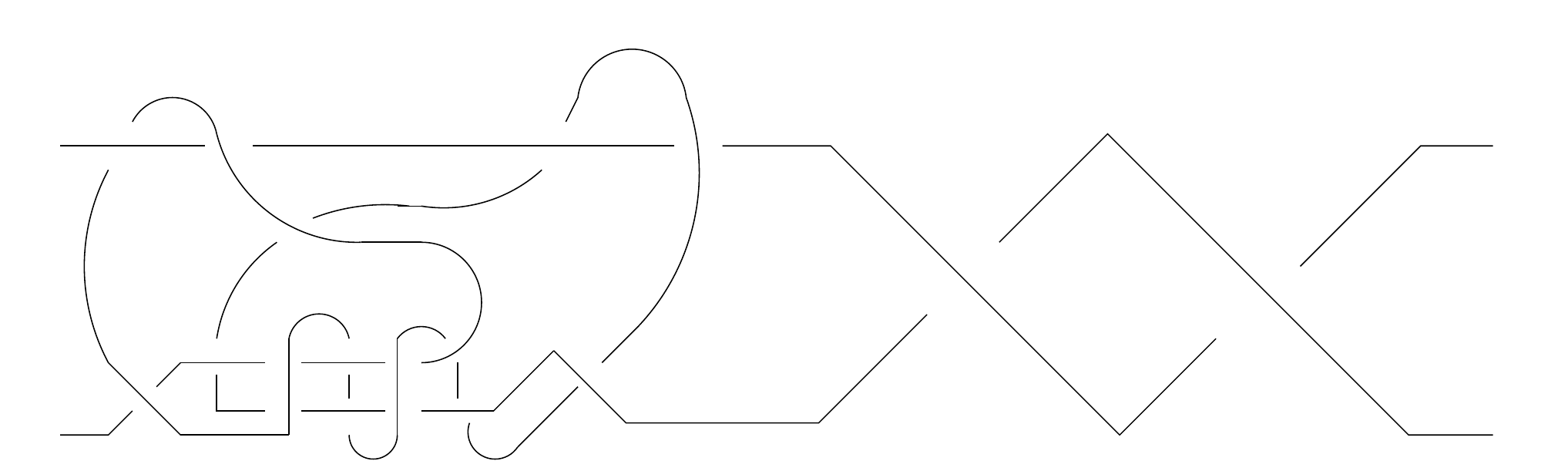}
\end{minipage}
\hspace{5pt}=\hspace{5pt}
\begin{minipage}{80pt}
\labellist\tiny
\hair=1pt
\pinlabel {$1$} [b] at 23 204
\pinlabel {$0$} [b] at 30 10
\pinlabel {$3$} [b] at 122 243
\pinlabel {$2$} [b] at 236 239
\pinlabel {$0$} [b] at 180 204
\pinlabel {$4$} [b] at 306 203
\pinlabel {$1$} [b] at 538 204
\pinlabel {$0$} [b] at 538 96
\pinlabel {$3$} [l] at 310 70
\endlabellist
\includegraphics[width=80pt]{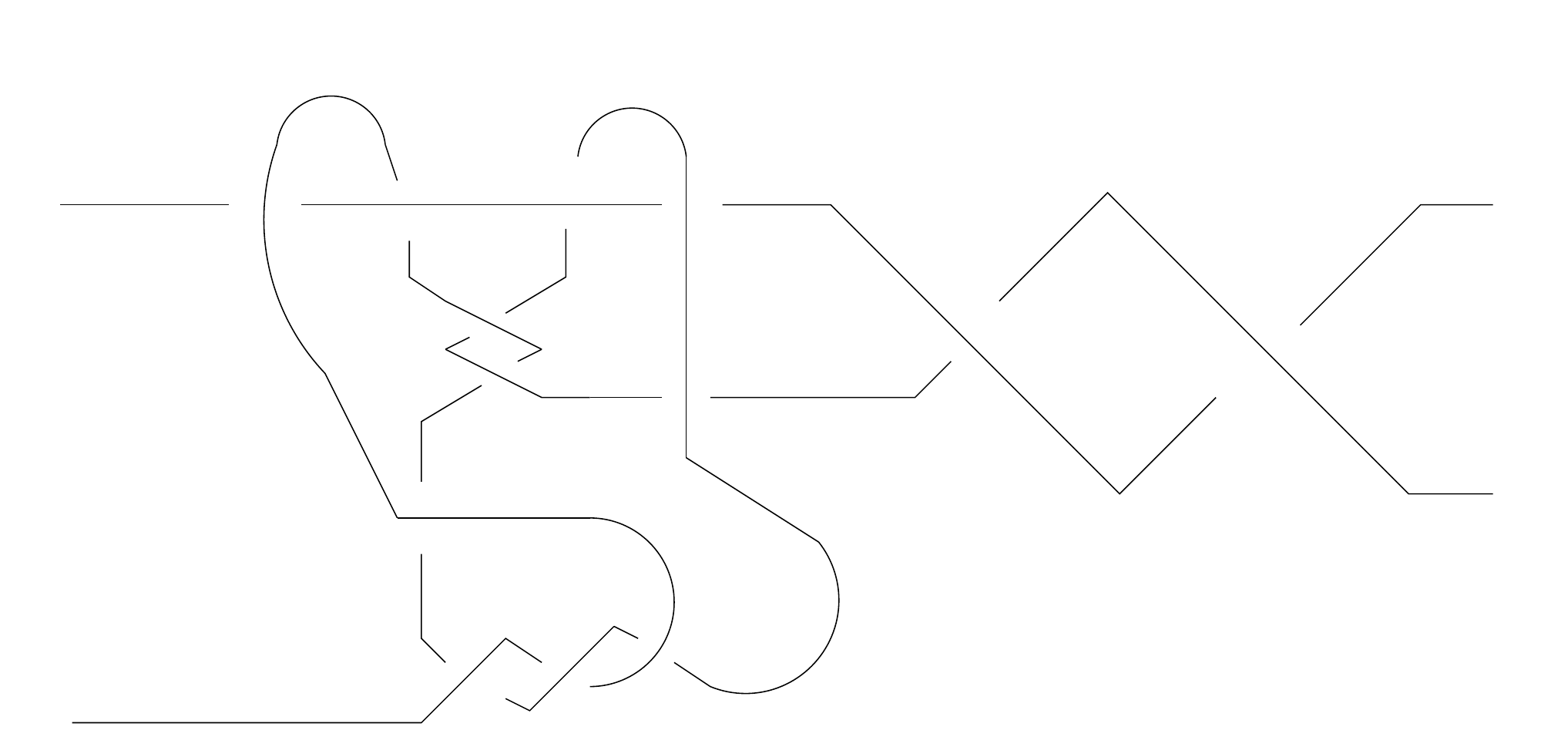}%
\end{minipage} \\
=\quad\begin{minipage}{100pt}
\labellist\tiny
\hair=2pt
\pinlabel {$1$} [r] at 18 74
\pinlabel {$0$} [r] at 24 24
\pinlabel {$3$} [b] at 118 257
\pinlabel {$2$} [b] at 232 252
\pinlabel {$1$} [l] at 553 217
\pinlabel {$0$} [l] at 553 109
\pinlabel {$4$} [tl] at 380 70
\endlabellist
\includegraphics[width=100pt]{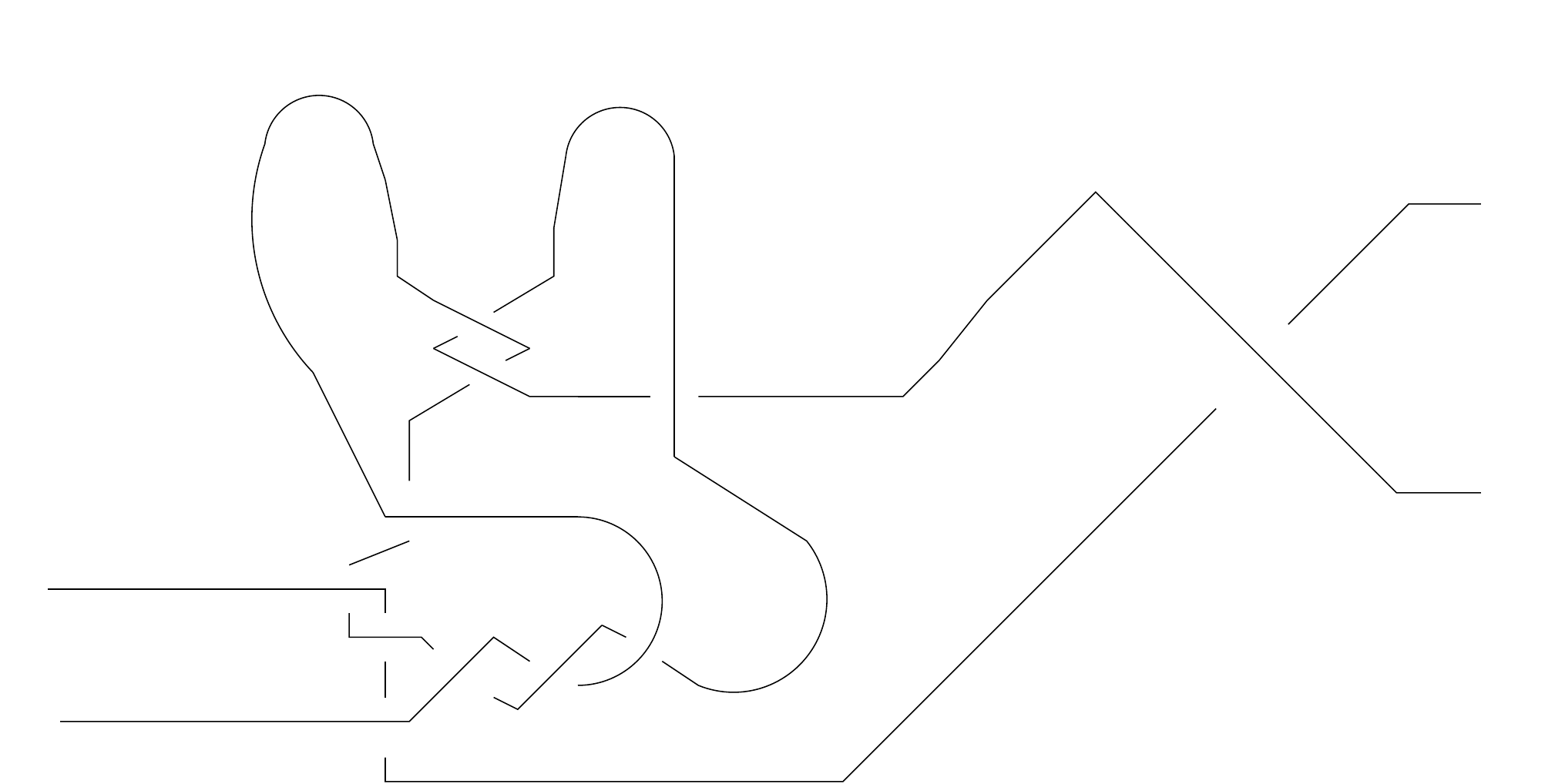}%
\end{minipage}
\quad\overset{RR}{\sim}\quad
\begin{minipage}{100pt}
\labellist\tiny
\pinlabel {$0$} [b] at 20 2
\pinlabel {$1$} [b] at 20 110
\pinlabel {$1$} [l] at 324 37
\endlabellist
\includegraphics[width=100pt]{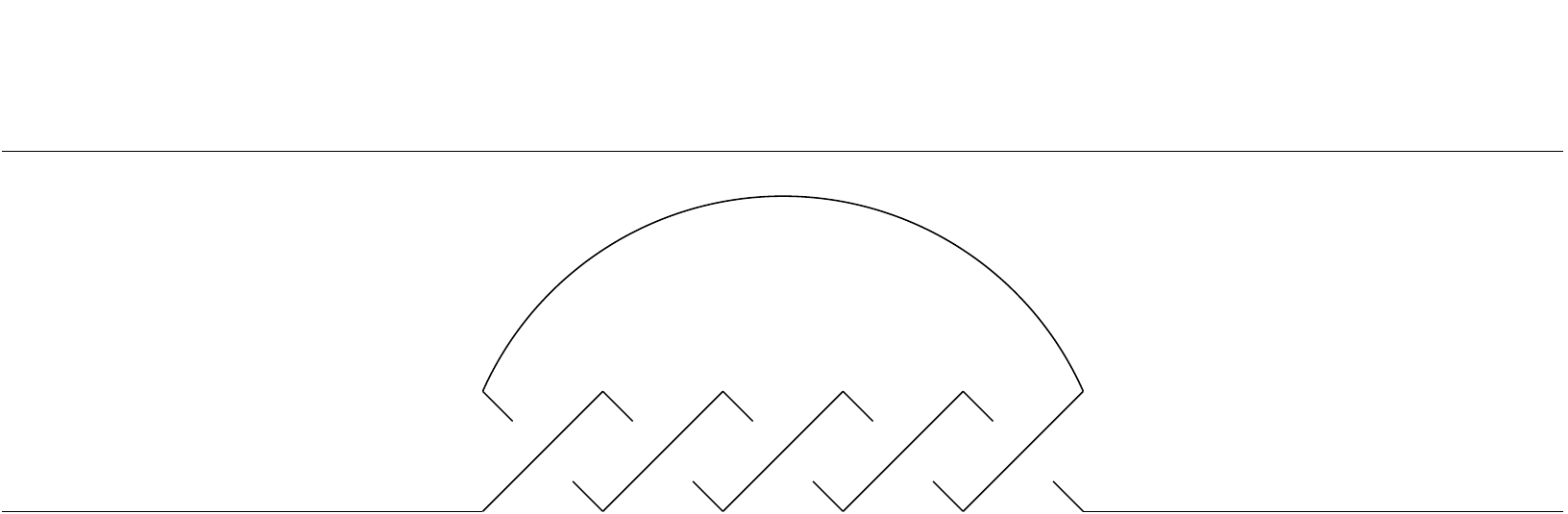}
\end{minipage}
\end{multline*}
To get from the first diagram to the second diagram, use the $R2G$
move on the arc coloured $4$ through the two arcs coloured $3$.

The $d=0$ case is trivial. For $d=2$ instead of $d=1$, in the same
way we obtain
$$\begin{minipage}{80pt}
\labellist\tiny
\hair=1pt
\pinlabel {$0$} [b] at 23 10
\pinlabel {$2$} [b] at 23 118
\pinlabel {$1$} [b] at 130 155
\pinlabel {$4$} [b] at 235 155
\pinlabel {$3$} [b] at 296 118
\pinlabel {$1$} [b] at 272 14
\pinlabel {$3$} [b] at 540 118
\pinlabel {$0$} [b] at 540 10
\pinlabel {$2$} [b] at 170 10
\pinlabel {$0$} [b] at 180 118
\endlabellist
\includegraphics[width=80pt]{\figdir/2loops}%
\end{minipage}\hspace{10pt}\seq\hspace{10pt}
\begin{minipage}{100pt}
\labellist\tiny
\pinlabel {$0$} [b] at 20 2
\pinlabel {$2$} [b] at 20 110
\pinlabel {$2$} [l] at 324 37
\endlabellist
\includegraphics[width=100pt]{\figdir/lktorusfinal1}
\end{minipage}$$
The proof of $r^{2}\cdot \tau^{-1}\seq \mathbbm{1}\hash
(\overline{5_{1}},\overline{\rho_{d}})$ is the same up to
reflection.
\end{proof}

This implies:
\begin{equation}\label{E:10move}
\mathbbm{1}\hash (5_{1},\rho_{d})^{2}\seq l^{4}\cdot \tau^{2}\seq
r^{2}\cdot \tau^{4}= r^{2}\cdot\tau^{-1}\cdot\tau^{5}\seq
\tau^{5}\hash (\overline{5_{1}},\overline{\rho_{d}})
\end{equation}

\subsection{The $10$--move}\label{SS:10move}

By Equation \ref{E:10move} on the five half-twists of the
$(5,2)$--torus knot, we have
$$
\begin{minipage}{50.6pt}
\labellist\tiny
\pinlabel {$0$} [r] at 35 258
\pinlabel {$1$} [b] at 93 216
\pinlabel {$2$} [b] at 165 216
\pinlabel {$3$} [b] at 237 216
\pinlabel {$4$} [b] at 309 216
\endlabellist
\includegraphics[width=0.7in]{\figdir/l25torus}
\end{minipage} \hspace{10pt}\hash\hspace{7pt}
\begin{minipage}{50.6pt}
\labellist\tiny
\pinlabel {$0$} [r] at 35 258
\pinlabel {$1$} [b] at 93 216
\pinlabel {$2$} [b] at 165 216
\pinlabel {$3$} [b] at 237 216
\pinlabel {$4$} [b] at 309 216
\endlabellist
\includegraphics[width=0.7in]{\figdir/l25torus}
\end{minipage} \hspace{10pt}\hash\hspace{7pt}
\begin{minipage}{50.6pt}
\labellist\tiny
\pinlabel {$0$} [r] at 35 258
\pinlabel {$1$} [b] at 93 216
\pinlabel {$2$} [b] at 165 216
\pinlabel {$3$} [b] at 237 216
\pinlabel {$4$} [b] at 309 216
\endlabellist
\includegraphics[width=0.7in]{\figdir/l25torus}
\end{minipage}
\hspace{15pt}\overset{\ref{E:10move}}{\sim}\hspace{10pt}
\begin{minipage}{50.6pt}
\labellist\tiny
\pinlabel {$0$} [r] at 35 258
\pinlabel {$1$} [b] at 93 216
\pinlabel {$2$} [b] at 165 216
\pinlabel {$3$} [b] at 237 216
\pinlabel {$4$} [b] at 309 216
\endlabellist
\includegraphics[width=0.7in]{\figdir/r25torus}
\end{minipage}\hspace{10pt}\hash\hspace{7pt} \begin{minipage}{50.6pt}
\labellist\tiny
\pinlabel {$0$} [r] at 35 258
\pinlabel {$1$} [b] at 93 216
\pinlabel {$2$} [b] at 165 216
\pinlabel {$3$} [b] at 237 216
\pinlabel {$4$} [b] at 309 216
\endlabellist
\includegraphics[width=0.7in]{\figdir/r25torus}
\end{minipage}
$$
\noindent Repeating for all different colourings of the $5_{1}$ knot
we obtain:
\begin{align}
(5_{1},\rho_{1})^{\hash 3}&\seq (\overline{5_{1}},\overline{\rho_{1}})^{\hash
2}\label{E:51-32-1}\\
(5_{1},\rho_{2})^{\hash 3}&\seq (\overline{5_{1}},\overline{\rho_{2}})^{\hash
2}\label{E:51-32-2}\\
(\overline{5_{1}},\overline{\rho_{1}})^{\hash 3}&\seq (5_{1},\rho_{1})^{\hash
2}\\
(\overline{5_{1}},\overline{\rho_{2}})^{\hash 3}&\seq (5_{1},\rho_{2})^{\hash
2}\label{E:51-32-4}
\end{align}
\noindent where the exponent is with respect to the connect-sum
operation.

\begin{lem}\label{L:51generator}
Any $5$--coloured knot $K$ is equivalent modulo surgery in
$\ker\rho$ to itself connect-summed with five $5_{1}$ knots of any
handedness and any $5$--colouring.
\end{lem}

\begin{proof}
Let $(K_{0},\rho)$ be a $5$--coloured knot in $S^{3}$, and let $B$
be a ball in $S^{3}$ such that $T:=B\cap K_{0}$ is the trivial
$(2,2)$--tangle $\mathbbm{1}$ (two unknotted untangled parallel
strands), where the colouring of the upper arc is $1$ and the
colouring of the lower arc is $0$.

$T$ is the product of two trivial $(2,2)$--tangles $T_{0}$ and
$T_{1}$. We first add and subtract $\frac{-3}{2}$ twists from
$T_{0}$ by a series of Reidemeister $RII$ moves, and then we extract
a copy of $l$ from these twists by
$$
\begin{minipage}{50pt}
\raisebox{7pt}{\includegraphics[width=50pt]{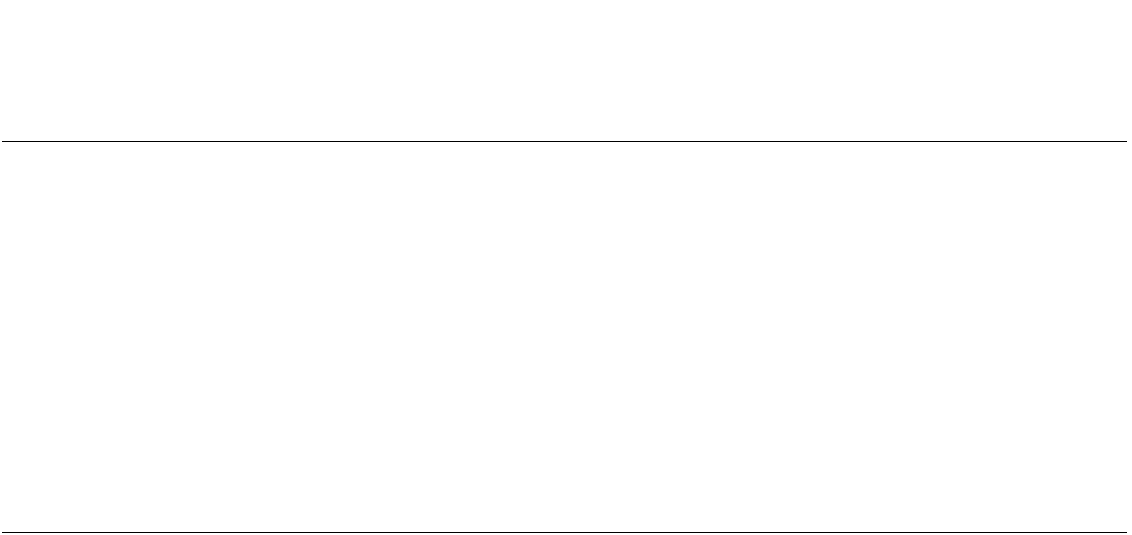}}
\end{minipage}
\quad=\quad
\begin{minipage}{80pt}
\includegraphics[width=80pt]{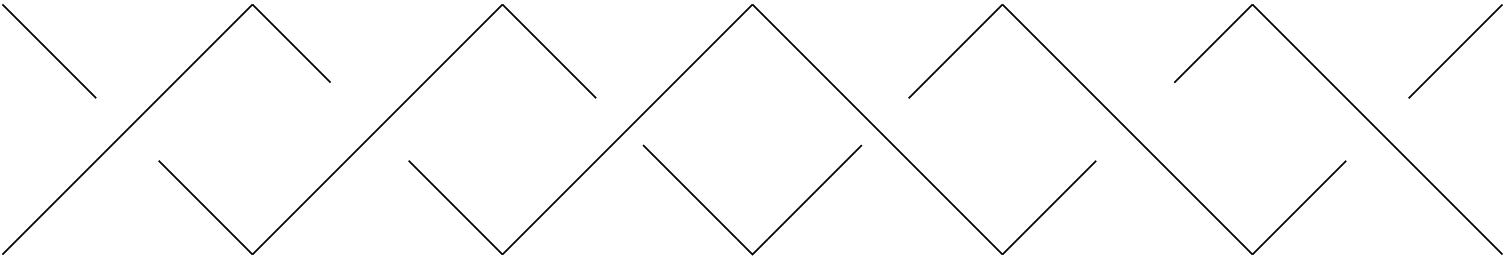}
\end{minipage}
\quad=\quad
\begin{minipage}{80pt}
\includegraphics[width=80pt]{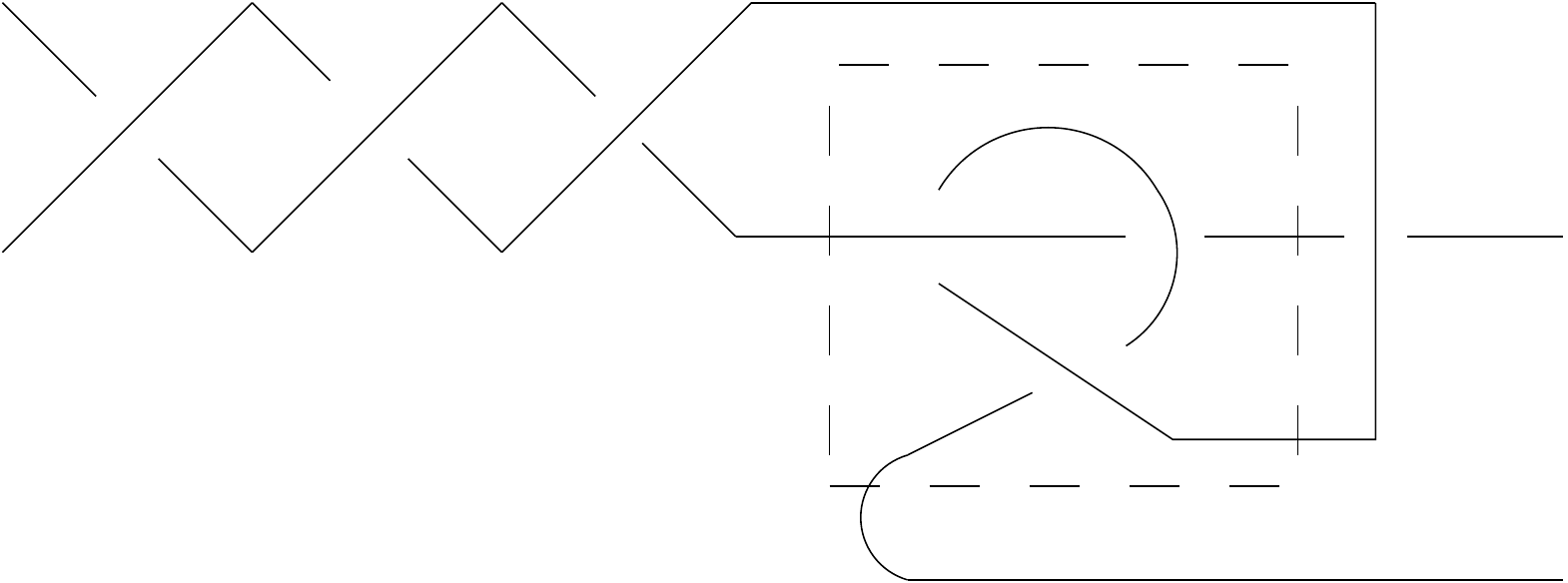}
\end{minipage}
$$
\noindent By \fullref{L:taus} this allows us to extract a copy of
$(5_{1},\rho_{1})$ from the modified $T_{0}$. We repeat this
procedure six times. Then by repeated use of Equations
\ref{E:51-32-1} and \ref{E:51-32-2} we have:
$$
(5_{1},\rho_{1})^{\hash 6}\sim
(\overline{5_{1}},\overline{\rho_{1}})^{\hash
6}\hash(5_{1},\rho_{1})^{\hash 2}\sim(5_{1},\rho_{1})^{\hash 11}
$$
By isotopy, we take five copies of $(5_{1},\rho_{1})$ into $T_{1}$
(making it isotopic to a copy of $(5_{1},\rho_{1})^{\hash 5}$),
while returning the remaining six copies of $(5_{1},\rho_{1})$ to
their initial positions in the modified $T_{0}$. Reversing the
process on what was $T_{0}$, we may make it once again into a
trivial $(2,2)$--tangle. Since we could have chosen $T$ to have any
colouring of its upper and lower arcs (by different choices of $B$),
we may determine the colouring of the $5_1$ knots which we obtain as
connect-summands.
\end{proof}

\begin{remi}
Two coloured knots are said to be $2p$--move equivalent if they are
related by a sequence of $2p$--moves:
$$
\begin{minipage}{30pt}
\includegraphics[width=30pt]{\figdir/ablines}
\end{minipage}
\quad\underset{\text{$2p$--move}}{\To}\quad
\begin{minipage}{60pt}
\labellist\small
\pinlabel {$p$} at 215 245
\endlabellist
\includegraphics[width=60pt]{\figdir/twist}
\end{minipage}
$$
\end{remi}

\begin{prop}\label{P:10move}
Let $K_{0}$ and $K_{1}$ be $10$--move equivalent $5$--coloured
knots. Then $K_{0}\seq K_{1}\hash (K')^{\hash 2}$ where
$K'=(5_{1},\rho_{1})$ or $(5_{1},\rho_{2})$ or
$(\overline{5_{1}},\overline{\rho_{1}})$ or
$(\overline{5_{1}},\overline{\rho_{2}})$.
\end{prop}

\begin{proof}
We may assume that $K_{1}$ is obtained from $K_{0}$ by a single
$10$--move, by adding ten half twists in two parallel strands of
$K_{0}$. By \fullref{L:51generator}, $K_{0}\seq K_{0}\hash
(5_{1},\rho_{1})^{\hash 5}$ (the proof is the same for different
colouring and handedness of the $5_{1}$--knot). We now look at a
trivial $(2,2)$--tangle $\mathbbm{1}$ in $K_{0}$ connect-summed with
$(5_{1},\rho_{1})^{\hash 2}$. By Equation \ref{E:10move} this is
equivalent modulo surgery in $\ker\rho$ to
$\tau^{5}\hash(\overline{5_{1}},\overline{\rho_{1}})$. So we have
converted $K_{0}$ to $K_{1}\hash (5_{1},\rho_{1})^{\hash 3}\hash
(\overline{5_{1}},\overline{\rho_{1}})$ by $\pm 1$--framed surgeries
in $\ker{\rho}$. By Equations \ref{E:51-32-1}--\ref{E:51-32-4} we
have
$$
(5_{1},\rho_{1})^{\hash 3}\hash
(\overline{5_{1}},\overline{\rho_{1}})\seq
(\overline{5_{1}},\overline{\rho_{1}})^{\hash 3}\seq
(5_{1},\rho_{1})^{\hash 2}
$$
\noindent which proves the proposition.
\end{proof}

\subsection{Reducing $S(m,n)$ to $5_{1}$}\label{SS:Smnreduce}

In the present section we prove equivalences between the $5_{1}$
knot and $5$--coloured knots of the form $S(m,n)$ with $-2\leq
m,n\leq 2$ (see \fullref{fig1}). By \fullref{P:prop9} there are
two such knots up to reflection--- $S(2,2)$, and $S(-1,1)$--- and
the $5$--colouring of each of them is unique (up to automorphism of
$D_{2p}$).

\begin{prop}\label{P:Knotequiv}
The $5$--coloured knots $S(2,2)$, $S(-1,1)$, and their mirror images
are all equivalent to connect-sums of $5_{1}$ knots modulo surgery
in $\ker \rho$.
\end{prop}

\begin{proof}
$S(-1,1)\seq 5_{1}$:
\begin{multline*}
\qua\begin{minipage}{57.9pt}
\labellist\small
\hair=1pt
\pinlabel {$0$} [b] at 260 420
\pinlabel {$4$} [t] at 260 60
\pinlabel {$1$} [br] at 230 280
\pinlabel {$2$} [br] at 40 305
\endlabellist
\includegraphics[width=0.8in]{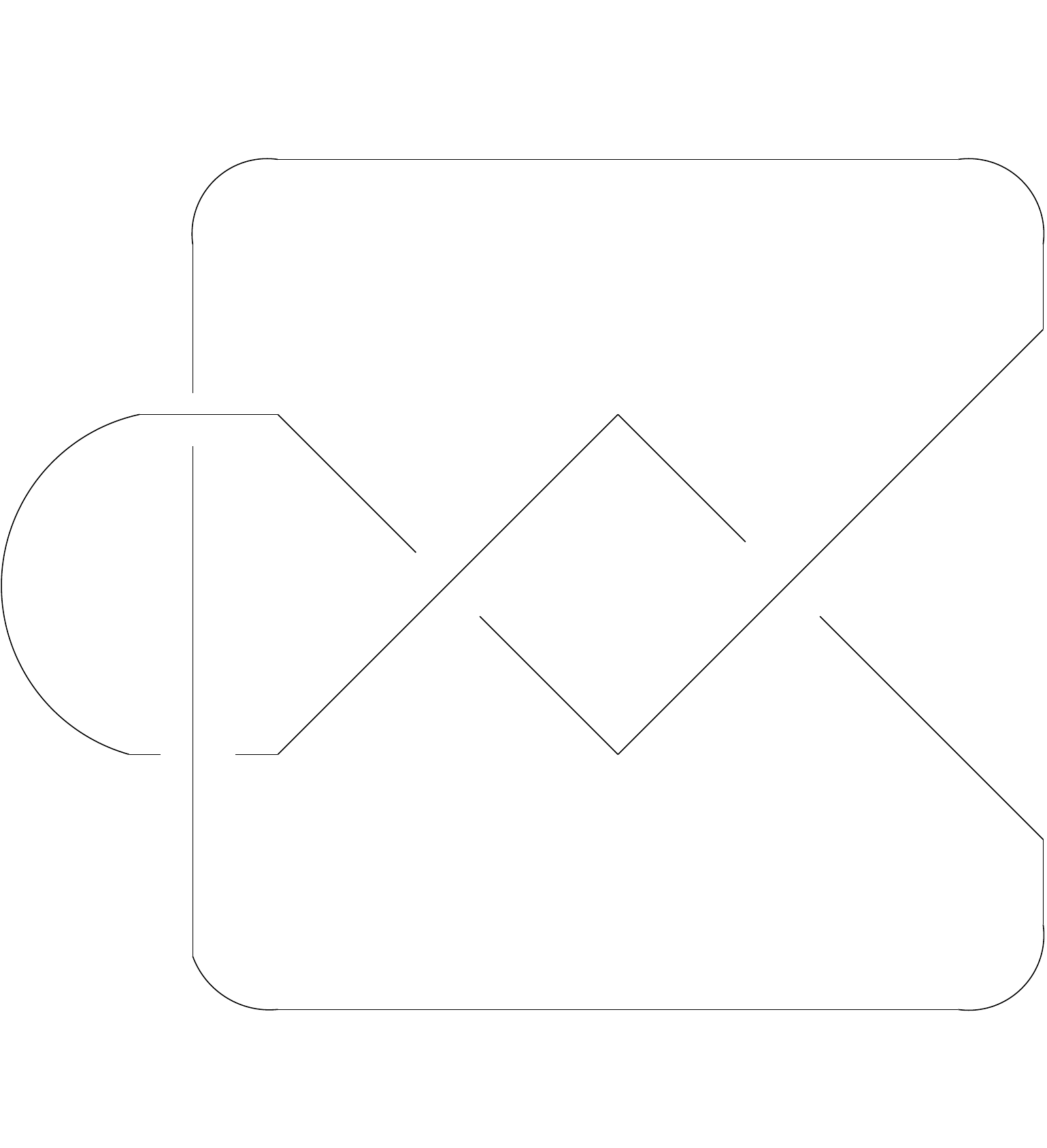}
\end{minipage}
\quad\overset{\text{\fullref{L:51generator}}}\seq\quad
\begin{minipage}{57.9pt}
\labellist\small
\hair=1pt
\pinlabel {$0$} [b] at 260 420
\pinlabel {$4$} [t] at 260 60
\pinlabel {$1$} [br] at 230 280
\pinlabel {$2$} [br] at 40 305
\endlabellist
\includegraphics[width=0.8in]{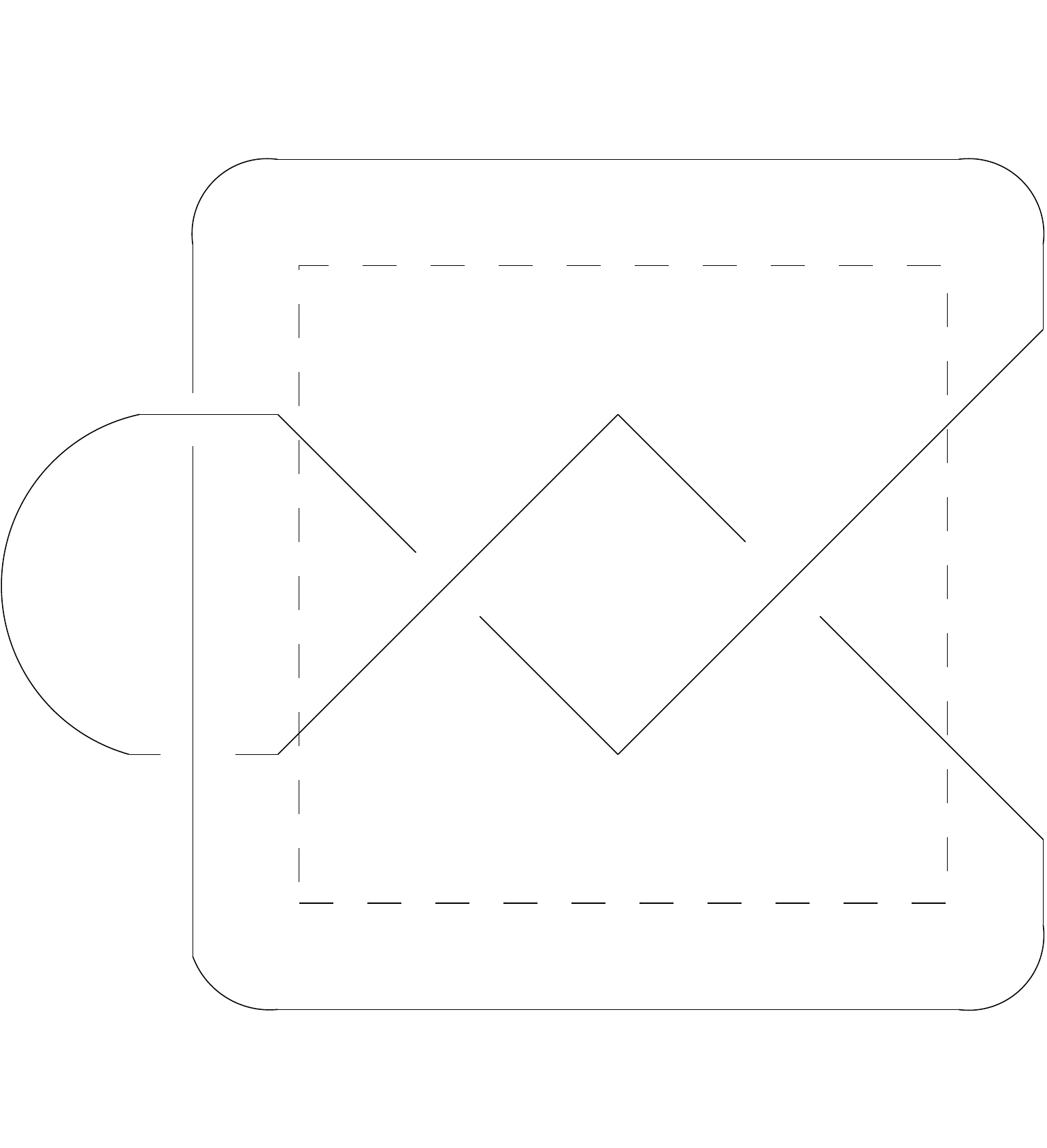}
\end{minipage}
\quad\hash\quad
\begin{minipage}{58pt}
\labellist\tiny
\pinlabel {$0$} [r] at 35 258
\pinlabel {$1$} [b] at 93 216
\pinlabel {$2$} [b] at 165 216
\pinlabel {$3$} [b] at 237 216
\pinlabel {$4$} [b] at 309 216
\endlabellist
\includegraphics[width=0.8in]{\figdir/r25torus}
\end{minipage}^{\quad\hash 5}\\
\seq\quad\begin{minipage}{61.5pt}
\labellist\small
\hair=2pt
\pinlabel {$0$} [b] at 260 420
\pinlabel {$4$} [t] at 260 60
\pinlabel {$2$} [br] at 40 305
\endlabellist
\includegraphics[width=0.85in]{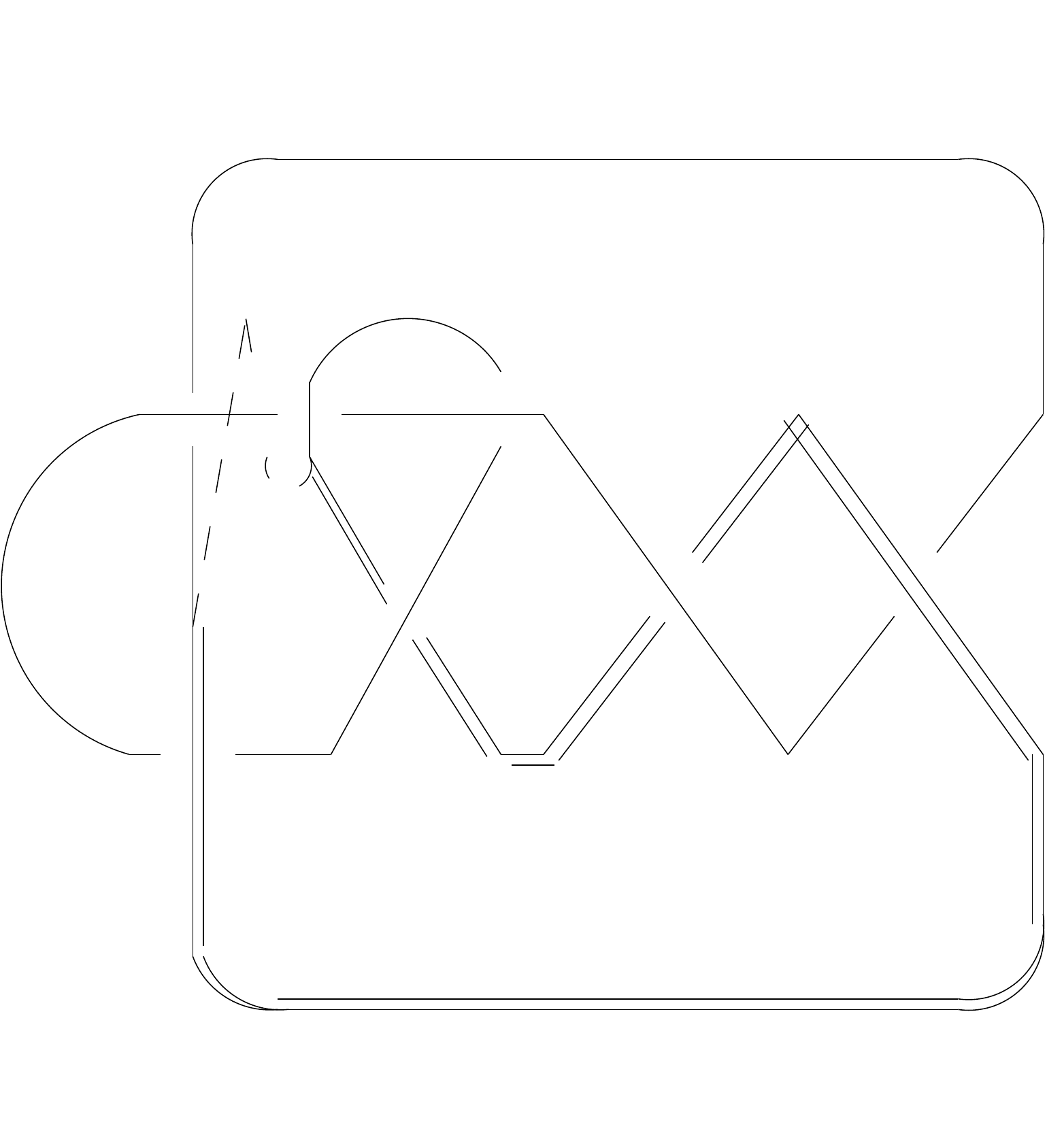}
\end{minipage}\hspace{5pt}\hash\hspace{7pt}
\begin{minipage}{58pt}
\labellist\tiny
\pinlabel {$0$} [r] at 35 258
\pinlabel {$1$} [b] at 93 216
\pinlabel {$2$} [b] at 165 216
\pinlabel {$3$} [b] at 237 216
\pinlabel {$4$} [b] at 309 216
\endlabellist
\includegraphics[width=0.8in]{\figdir/r25torus}
\end{minipage}^{\quad\hash 4}
\hspace{0.2pt}=\hspace{10pt}
\begin{minipage}{50.6pt}
\labellist\tiny
\pinlabel {$0$} [r] at 35 258
\pinlabel {$2$} [b] at 93 216
\pinlabel {$4$} [b] at 165 216
\pinlabel {$1$} [b] at 237 216
\pinlabel {$3$} [b] at 309 216
\endlabellist
\includegraphics[width=0.7in]{\figdir/r25torus}
\end{minipage}\hspace{5pt}\hash\hspace{5.5pt}
\begin{minipage}{50.6pt}
\labellist\tiny
\pinlabel {$0$} [r] at 35 258
\pinlabel {$1$} [b] at 93 216
\pinlabel {$2$} [b] at 165 216
\pinlabel {$3$} [b] at 237 216
\pinlabel {$4$} [b] at 309 216
\endlabellist
\includegraphics[width=0.7in]{\figdir/r25torus}
\end{minipage}^{\quad\hash
4} \\
=(5_{1},\rho_{2})\hash(5_{1},\rho_{1})^{\hash 4}
\end{multline*}
where the second equivalence follows from \fullref{L:taus} by considering the copy of $\tau^{-1}$ in $S(-1,1)$ in
the dotted box connect-summed with one of the copies of
$(5_{1},\rho_{1})$. We then have
$$\tau^{-1}\hash(5_{1},\rho_{1})=\tau^{-1}\cdot \tau\cdot l^{2}=r\cdot\tau$$
$S(2,2)\seq 5_{1}$:
\begin{multline*}
S(2,2)=\quad
\begin{minipage}{80pt}
\labellist\small
\hair=2pt
\pinlabel {$0$} [br] at 34 248
\pinlabel {$4$} [b] at 193 261
\pinlabel {$0$} [bl] at 343 236
\pinlabel {$1$} [tl] at 343 80
\pinlabel {$2$} [t] at 204 44
\pinlabel {$1$} [tr] at 36 51
\pinlabel {$3$} [bl] at 157 185
\endlabellist
\includegraphics[width=80pt]{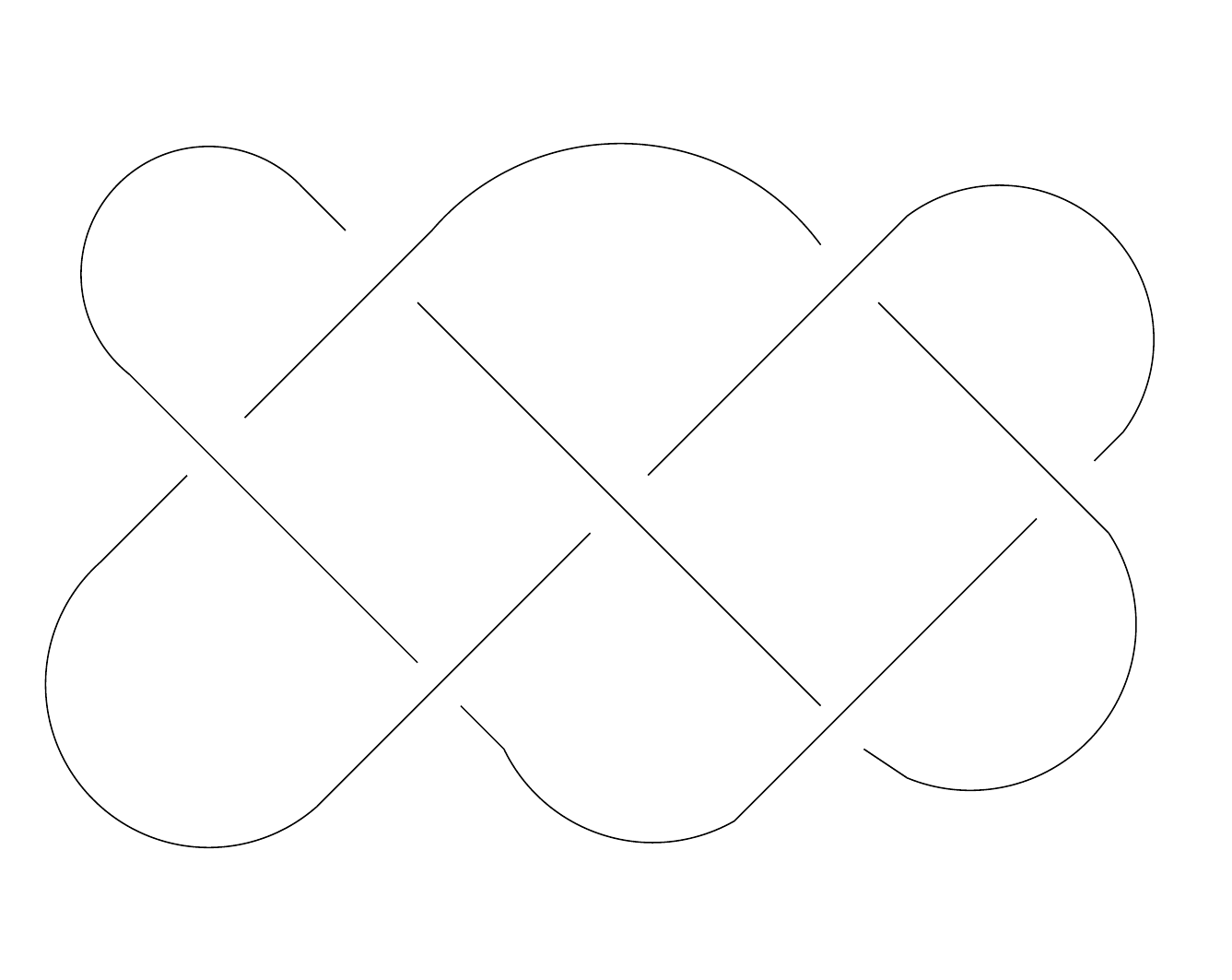}
\end{minipage}
\quad=\quad
\begin{minipage}{100pt}
\labellist\small
\pinlabel {$3$} [r] at 45 143
\pinlabel {$4$} [l] at 566 143
\pinlabel {$4$} [r] at 153 73
\pinlabel {$0$} [l] at 463 112
\pinlabel {$1$} [l] at 440 208
\endlabellist
\includegraphics[width=100pt]{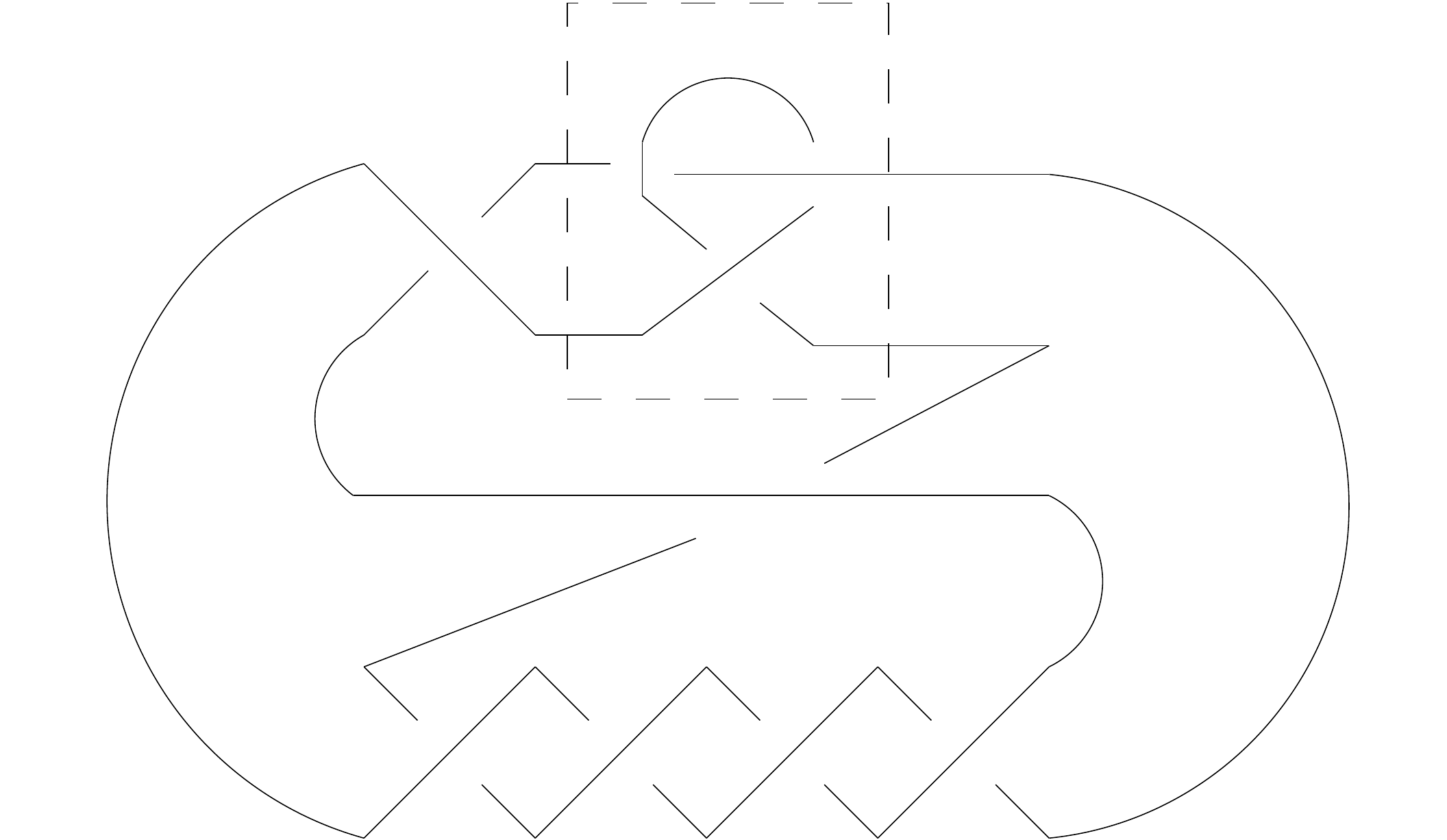}
\end{minipage} \\
\seq\quad
\begin{minipage}{100pt}
\labellist\small
\pinlabel {$3$} [r] at 45 143
\pinlabel {$4$} [l] at 566 143
\pinlabel {$4$} [r] at 153 73
\pinlabel {$0$} [l] at 463 112
\pinlabel {$1$} [l] at 440 208
\endlabellist
\includegraphics[width=100pt]{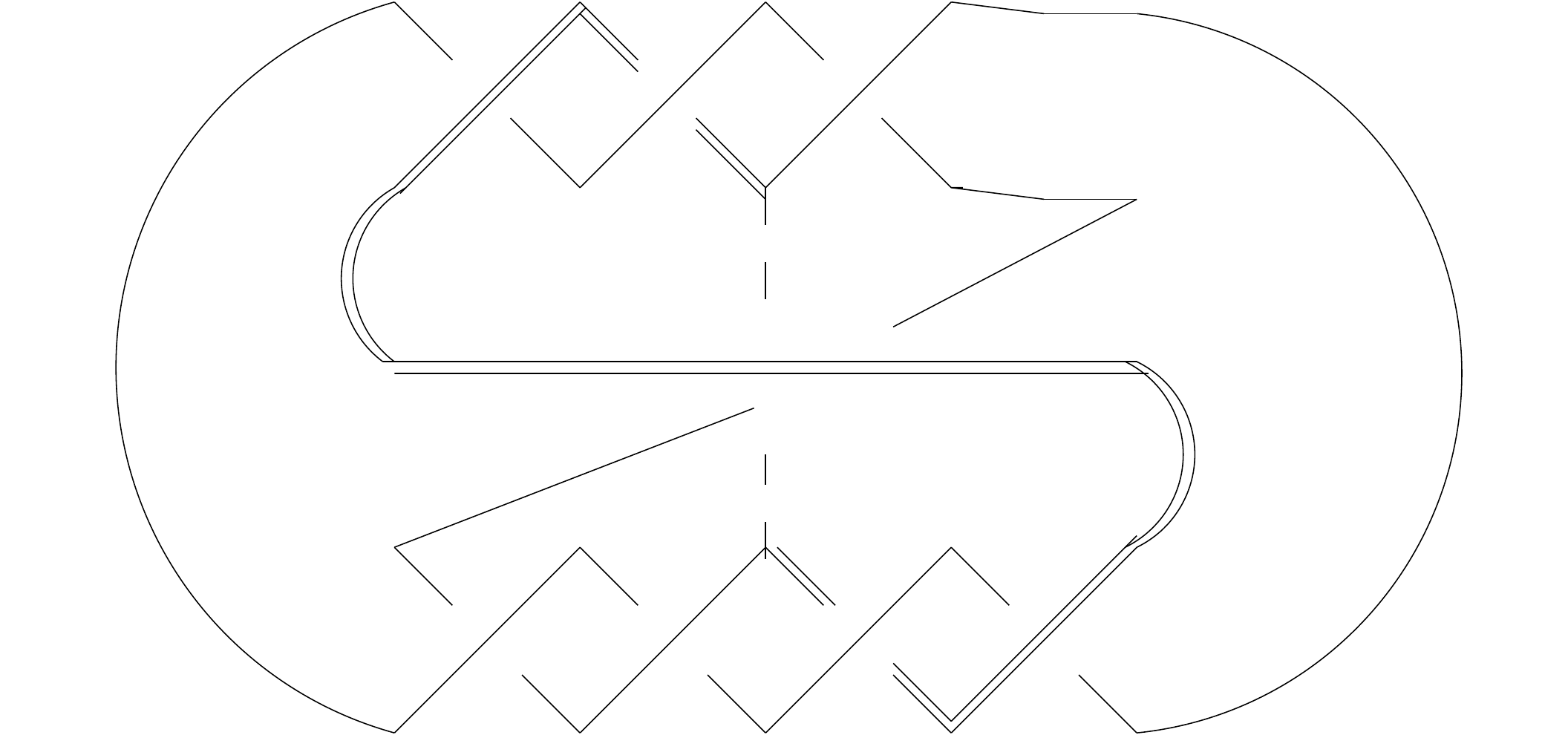}
\end{minipage}
\quad\hash\quad
\begin{minipage}{58pt}
\labellist\tiny
\pinlabel {$0$} [r] at 35 258
\pinlabel {$1$} [b] at 93 216
\pinlabel {$2$} [b] at 165 216
\pinlabel {$3$} [b] at 237 216
\pinlabel {$4$} [b] at 309 216
\endlabellist
\includegraphics[width=0.8in]{\figdir/r25torus}
\end{minipage} \\
=\quad
\begin{minipage}{0.8in}
\labellist\small
\pinlabel {$1$} [r] at 45 260
\pinlabel {$2$} [l] at 405 260
\pinlabel {$0$} [r] at 154 260
\pinlabel {$4$} [tr] at 155 52
\endlabellist
\includegraphics[width=0.8in]{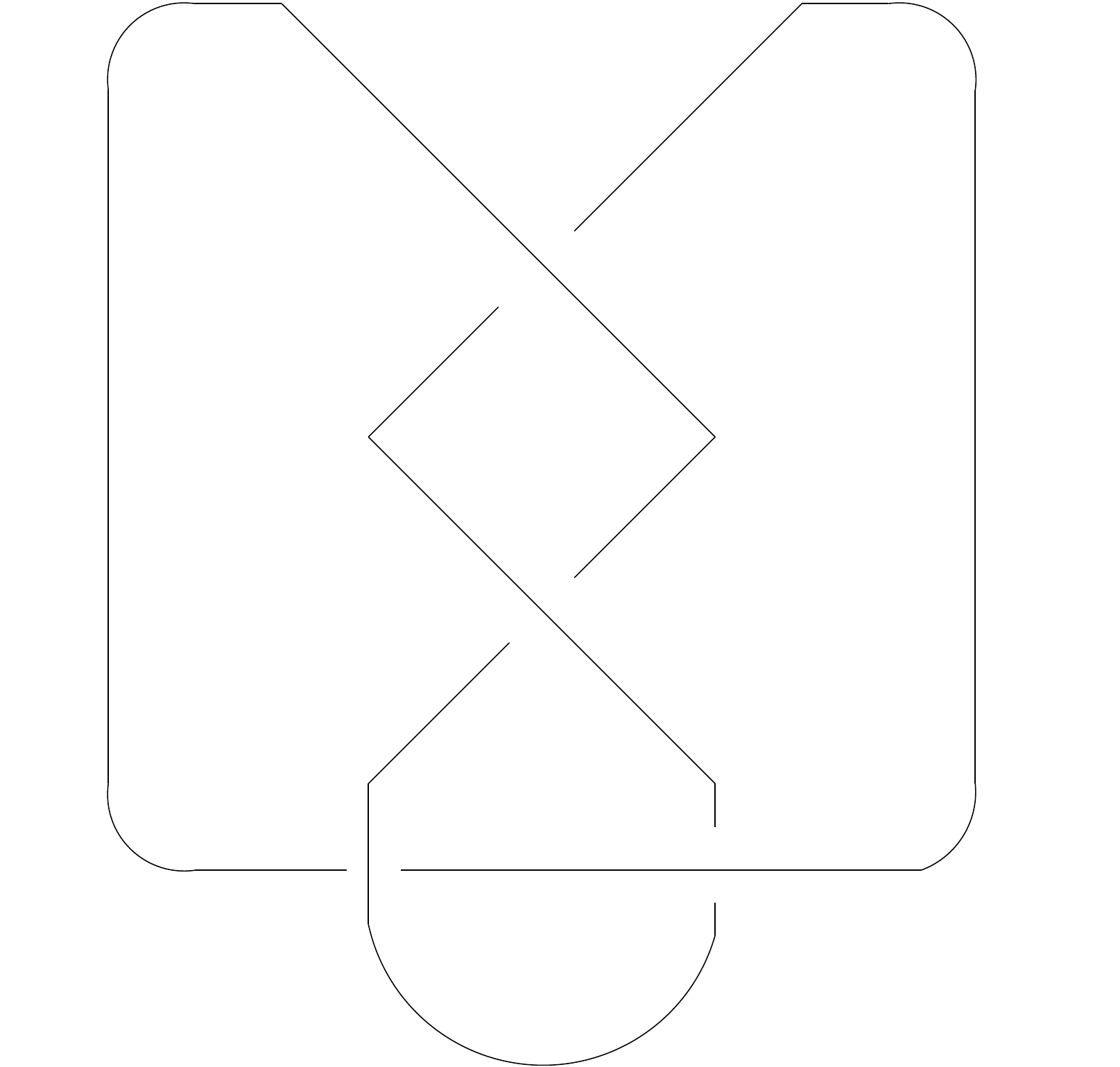}
\end{minipage}
\quad\hash\quad
\begin{minipage}{58pt}
\labellist\tiny
\pinlabel {$0$} [r] at 35 258
\pinlabel {$1$} [b] at 93 216
\pinlabel {$2$} [b] at 165 216
\pinlabel {$3$} [b] at 237 216
\pinlabel {$4$} [b] at 309 216
\endlabellist
\includegraphics[width=0.8in]{\figdir/r25torus}
\end{minipage}
\quad=S(-1,1)\hash(5_1,\rho_{1})
\end{multline*}
To get from the second diagram to  the third, use \fullref{L:taus}
to deduce
$$r=r\cdot\tau^{2}\cdot\tau^{-2}=l^{2}\cdot\tau\cdot
  \tau^{-2}=(5_{1},\rho_{a})\hash\tau^{-2}$$
where the copy of $r$ in the second diagram is drawn
inside a dotted box. In the final stage we have reduced to the
previous case.
\end{proof}

\subsection{Reducing the number of connect summands}\label{SS:51alg}

Let
$$(K,\rho)= (5_{1},\rho_{1})^{\hash g_{1}}\hash
(5_{1},\rho_{2})^{\hash g_{2}}\hash
(\overline{5_{1}},\overline{\rho_{1}})^{\hash g_{3}}\hash
(\overline{5_{1}},\overline{\rho_{2}})^{\hash g_{4}}$$
\noindent be a connect-sum of $5_{1}$--knots. The purpose of this
section is to prove that
$$(K,\rho)\seq (5_{1},\rho_{1})^{\hash n}\qquad\mbox{for some $n=1,2,3,4,5$.}$$
\noindent We achieve this by finding relations between connect-sums
of $5_{1}$ knots with different colouring and handedness.

\begin{lem}\label{L:csqab}\hfill
\begin{align}
(5_{1},\rho_{1})\hash(\overline{5_{1}},\overline{\rho_{1}})^{2}&\seq (5_{1},\rho_{2})\label{E:csqab-1}\\
(\overline{5_{1}},\overline{\rho_{1}})\hash(5_{1},\rho_{1})^{2}&\seq (\overline{5_{1}},\overline{\rho_{2}})\\
(5_{1},\rho_{2})\hash(\overline{5_{1}},\overline{\rho_{2}})^{2}&\seq (5_{1},\rho_{1})\\
(\overline{5_{1}},\overline{\rho_{2}})\hash(5_{1},\rho_{2})^{2}&\seq
(\overline{5_{1}},\overline{\rho_{1}})
\end{align}
\end{lem}

\begin{proof}
We prove Equation \ref{E:csqab-1}. The proofs of the other equations
in the lemma follow by reflection and automorphism of $D_{10}$.
Consider:
$$(5_{1},\rho_{1})\hash(\overline{5_{1}},\overline{\rho_{1}})^{2}=\hspace{7pt}
\begin{minipage}{58pt}
\labellist\tiny
\hair=1pt
\pinlabel {$0$} [r] at 30 258
\pinlabel {$1$} [bl] at 93 216
\pinlabel {$2$} [bl] at 165 216
\pinlabel {$3$} [bl] at 237 216
\pinlabel {$4$} [bl] at 309 216
\endlabellist
\includegraphics[width=0.8in]{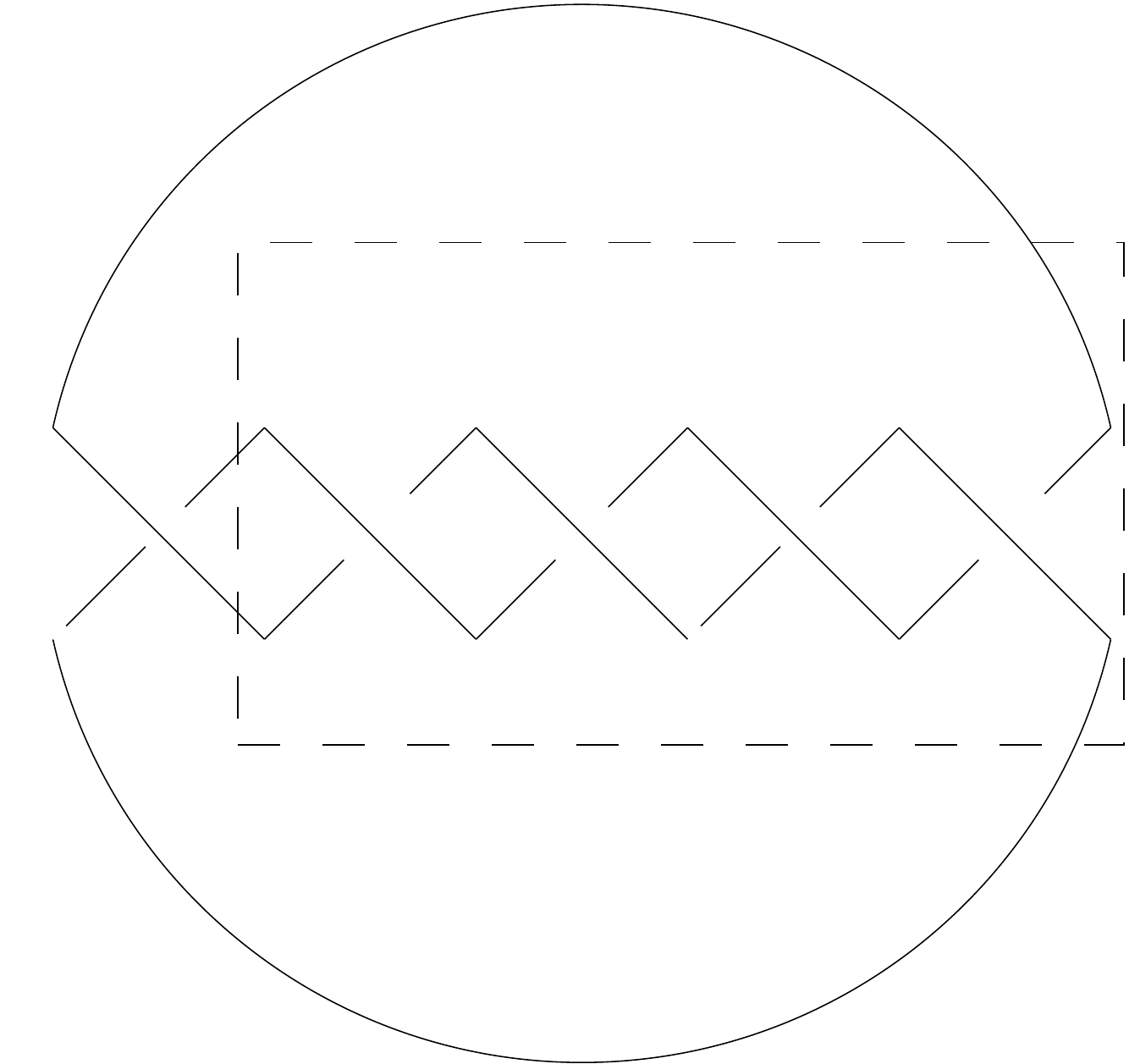}
\end{minipage} \hspace{15pt}\hash\hspace{7pt}
\begin{minipage}{58pt}
\labellist\tiny
\pinlabel {$0$} [r] at 35 258
\pinlabel {$1$} [b] at 93 216
\pinlabel {$2$} [b] at 165 216
\pinlabel {$3$} [b] at 237 216
\pinlabel {$4$} [b] at 309 216
\endlabellist
\includegraphics[width=0.8in]{\figdir/l25torus}
\end{minipage} \hspace{15pt}\hash\hspace{7pt}
\begin{minipage}{58pt}
\labellist\tiny
\pinlabel {$0$} [r] at 35 258
\pinlabel {$1$} [b] at 93 216
\pinlabel {$2$} [b] at 165 216
\pinlabel {$3$} [b] at 237 216
\pinlabel {$4$} [b] at 309 216
\endlabellist
\includegraphics[width=0.8in]{\figdir/r25torus}
\end{minipage}$$
Consider two of the twists (four of the half-twists) in one of the
copies of $(\overline{5_{1}},\overline{\rho_{1}})$ as the
$(2,2)$--tangle $\tau^{2}$, and the remaining
$(\overline{5_{1}},\overline{\rho_{1}})$ and $(5_{1},\rho_{1})$ as
being connect-summed to this $(2,2)$--tangle. We have:
$$
\tau^{2}\hash(\overline{5_{1}},\overline{\rho_{1}})\seq r^{2}\cdot
\tau^{-1}\cdot \tau^{2}=r^{2}\cdot \tau\seq l\cdot \tau^{-1}\cdot
\tau=l
$$
The equation above implies that:
$$
(5_{1},\rho_{1})\hash(\overline{5_{1}},\overline{\rho_{1}})^{2}\seq\hspace{7pt}
\begin{minipage}{58pt}
\labellist\small
\hair=2pt
\pinlabel {$0$} [r] at 40 267
\pinlabel {$1$} [r] at 50 80
\pinlabel {$4$} [b] at 105 218
\pinlabel {$2$} [b] at 200 258
\endlabellist
\includegraphics[width=0.8in]{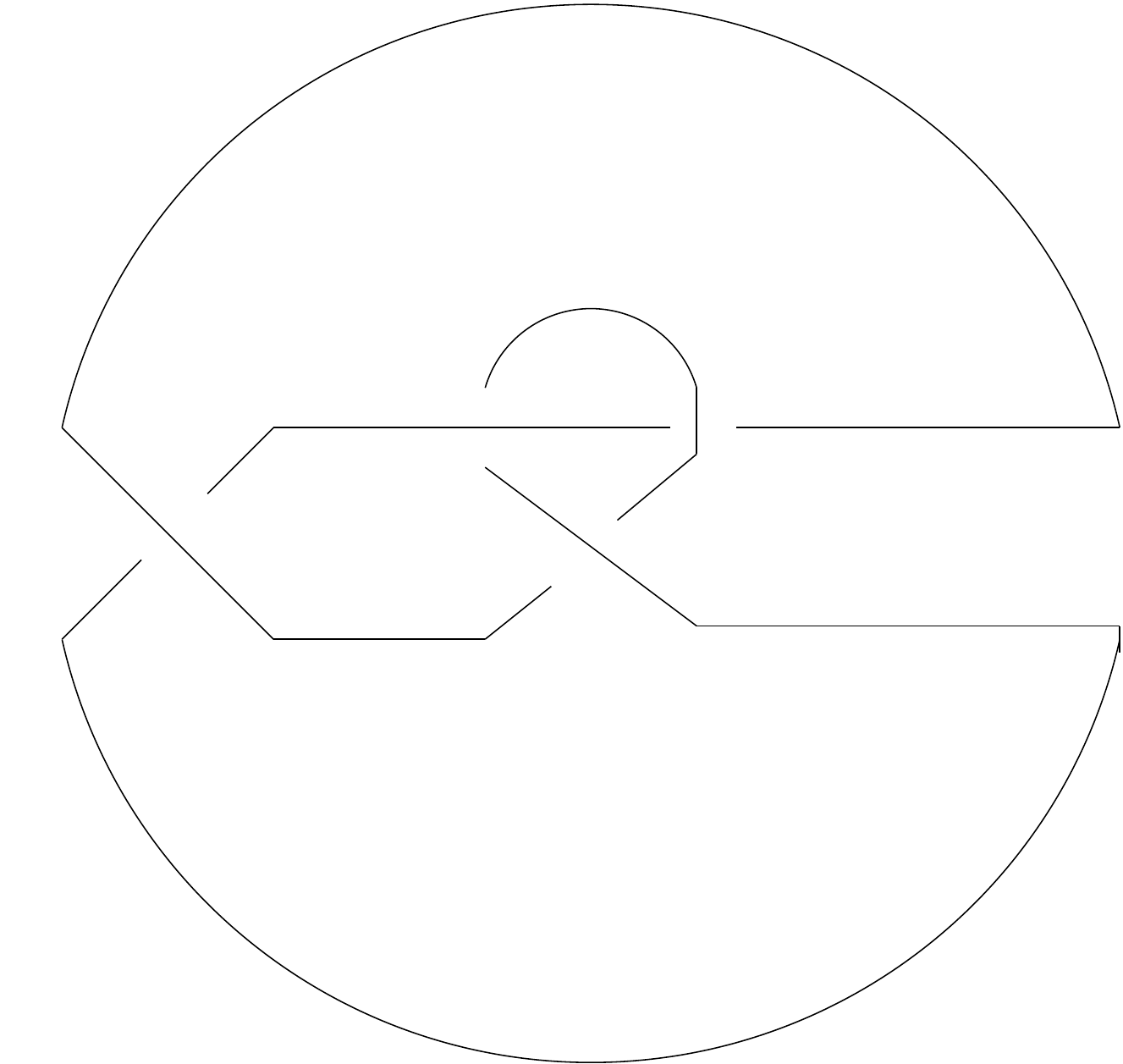}
\end{minipage}
\hspace{3pt}\hash\hspace{7pt}
\begin{minipage}{50.6pt}
\labellist\tiny
\pinlabel {$0$} [r] at 35 258
\pinlabel {$1$} [b] at 93 216
\pinlabel {$2$} [b] at 165 216
\pinlabel {$3$} [b] at 237 216
\pinlabel {$4$} [b] at 309 216
\endlabellist
\includegraphics[width=0.7in]{\figdir/r25torus}
\end{minipage}\hspace{7pt}=\hspace{6pt}
\begin{minipage}{58pt}
\labellist\small
\hair=1pt
\pinlabel {$0$} [b] at 260 420
\pinlabel {$4$} [t] at 260 60
\pinlabel {$1$} [br] at 230 280
\pinlabel {$2$} [br] at 40 305
\endlabellist
\includegraphics[width=0.8in]{\figdir/lying41}
\end{minipage}
\hspace{10pt}\hash\hspace{7pt} (5_{1},\rho_{1})
$$
\noindent where by the last three steps of the reduction of
$S(-1,1)$ to a connect-sum of $5_{1}$ knots in the proof of
\fullref{P:Knotequiv}, the knot at the end of this sequence
is equivalent to $(5_{2},\rho_{2})$ modulo surgery in $\ker\rho$.
\end{proof}

\begin{prop}\label{P:cs}
$$
(K,\rho)= (5_{1},\rho_{1})^{\hash g_{1}}\hash
(5_{1},\rho_{2})^{\hash g_{2}}\hash
(\overline{5_{1}},\overline{\rho_{1}})^{\hash g_{3}}\hash
(\overline{5_{1}},\overline{\rho_{2}})^{\hash g_{4}}\seq
(5_{1},\rho_{1})^{n}
$$
\noindent for $n=1,2,3,4,5$.
\end{prop}

\begin{proof}
Using \fullref{L:csqab} and Equations \ref{E:51-32-1}--
\ref{E:51-32-4} we have:
\begin{align*}
(5_{1},\rho_{2})&\seq (5_{1},\rho_{1})\hash
(\overline{5_{1}},\overline{\rho_{1}})^{\hash 2}\seq
(5_{1},\rho_{1})^{\hash 4}\\
\begin{split}
(\overline{5_{1}},\overline{\rho_{1}})&\seq (5_{1},\rho_{2})^{\hash
2}\hash (\overline{5_{1}},\overline{\rho_{2}})\seq
((5_{1},\rho_{1})^{\hash 2}\hash
(\overline{5_{1}},\overline{\rho_{1}}))\hash ((5_{1},\rho_{1})\hash
(\overline{5_{1}},\overline{\rho_{1}})^{\hash 2})^{\hash 2}\\&=
(5_{1},\rho_{1})^{\hash 4}\hash
(\overline{5_{1}},\overline{\rho_{1}})^{\hash 5}\seq
(5_{1},\rho_{1})^{\hash 9}\end{split}\\
(\overline{5_{1}},\overline{\rho_{2}})&\seq (5_{1},\rho_{1})^{\hash
2}\hash (\overline{5_{1}},\overline{\rho_{1}})\seq
(5_{1},\rho_{1})^{\hash 11}
\end{align*}
which gives
$$
(K,\rho)\seq (5_{1},\rho_{1})^{\hash n}
$$
\noindent for some positive integer $n$. We may take $n$ to be
between one and five since:
\begin{multline*}
(5_{1},\rho_{1})\seq (5_{1},\rho_{2})\hash
(\overline{5_{1}},\overline{\rho_{2}})^{\hash 2}\seq
((5_{1},\rho_{1})^{\hash 2}\hash
(\overline{5_{1}},\overline{\rho_{1}}))^{\hash 2}\hash
((5_{1},\rho_{1})\hash
(\overline{5_{1}},\overline{\rho_{1}})^{\hash 2})\\
= (5_{1},\rho_{1})^{\hash 5}\hash
(\overline{5_{1}},\overline{\rho_{1}})^{\hash 4}\seq
(5_{1},\rho_{1})^{\hash 2}\hash
(\overline{5_{1}},\overline{\rho_{1}})^{\hash 6}\seq
(5_{1},\rho_{1})^{\hash 6}
\end{multline*}
which completes the proof.
\end{proof}

\section{An invariant for $p$--coloured knots}\label{S:noteq}

In \fullref{S:3Col} and in \fullref{S:5Col}, it was shown that for
$p=3$ or $p=5$, any $p$--coloured knot $(K,\rho)$ can be reduced to a
connect-sum of $n$ left-hand $(p,2)$--torus knots with a given colouring
for some $1\leq n\leq p$ by a series of surgeries in $\ker\rho$. In the
present section we complete the proof of \fullref{T:3case} by showing
that this result cannot be improved.  Explicitly, for $p$ an odd prime
(in particular for $p\in\set{3,5}$) and for $r_{1}\neq r_{2}\bmod{p}$
the connect-sum of $r_{1}$ left-hand $(p,2)$--torus knots with a
given colouring is not equivalent modulo surgery in $\ker \rho$ to the
connect-sum of $r_{2}$ left-hand $(p,2)$--torus knots with the same
colouring. We prove this by finding a $\mathds{Z}/p\mathds{Z}$--valued
invariant for all $p$--coloured knots, invariant under surgery in
$\ker\rho$, which takes different values on each of these knots. Our
reference for algebraic topology is Spanier \cite{Spa66}.

We begin by recalling Cappell and Shaneson's notion of a \emph{mod
$p$ characteristic knot} \cite{CS75,CS84}. Recall the presentation
$D_{2p}:=\{t,s|t^{2}=s^{p}=1, tst=s^{p-1}\}$ for the dihedral group
of order $2p$. Let $(K,\rho)$ be a $p$--coloured knot with Seifert
surface $F$ and Seifert matrix $M$ with respect to a basis
$x_{1},x_{2},\ldots,x_{2n}$ of $H_{1}(F)$ with orientations induced
by the orientation of $F$. Let $(\xi_{1},\xi_{2},\ldots,\xi_{2n})$
be a basis for $H_{1}(S^{3}-F)$ oriented such that
$\mathrm{Link}(\xi_{i},x_{j})=\delta_{ij}$. Then there exists a link
$\mathrm{ch}(K,\rho)\subset F$ called a \emph{mod $p$ characteristic
link} such that
\begin{equation}\label{E:rhoxi}
\rho(\xi_{i})=s^{\mathrm{Link}(\mathrm{ch}(K,\rho),\xi_{i})}
\end{equation}
for all $1\leq i\leq 2n$. Recall also that $\rho(\tau)=t$
where $\tau$ denotes the generator of the infinite cyclic covering
group. In fact $\mathrm{ch}(K,\rho)$ may be chosen to be a knot, but
we do not require this fact here.

Let $\hC$ denote the $2$--fold covering of $S^{3}$ branched over $K$,
and let $pr$ denote the covering projection. The $p$--colouring
$\rho\co \pi_{1}(S^{3}-K)\to D_{2p}$ restricts in the
double-covering to a map $\rho\prime\co H_{1}(\hC;\mathds{Z})\to
\mathds{Z}/p\mathds{Z}$ which corresponds to a cohomology class
$a\in H^{1}(\hC;\mathds{Z}/p\mathds{Z})$  by the
universal-coefficient theorem for cohomology.

To simplify notation, we define $\alpha\ass\
pr^{-1}\mathrm{ch}(K,\rho)$. Since $\xi_{i}$ is in the complement of
the Seifert surface, its pre-image in $\hC$ is two disjoint circles,
while the mod $p$ characteristic link is contained in $F$ and its
lift therefore has only one component.  Let $[pr^{-1}\xi_{i}]\in
H_{1}(\hC;\mathds{Z})$ and $[\alpha]\in H_{1}(\hC;\mathds{Z})$ denote
the homology classes represented by $pr^{-1}\xi_{i}$ and $\alpha$
correspondingly. Equation \ref{E:rhoxi} tells us that
\begin{multline*}
a\frown [pr^{-1}\xi_{i}]= \langle a,[pr^{-1}\xi_{i}]\rangle=
s^{\mathrm{Link}(\alpha,pr^{-1}\xi_{i})}=\\ \langle
D\partial^{-1}[\alpha],\ [pr^{-1}\xi_{i}]\rangle \bmod p=
D\partial^{-1}[\alpha]\bmod p \frown [pr^{-1}\xi_{i}]
\end{multline*}
for all $1\leq i\leq 2n$, where $D$ denotes the
Poincar\'{e} duality isomorphism, and $\partial^{-1}[\alpha]$
denotes the pre-image of $[\alpha]$ under the Bockstein homomorphism
$$
\cdots \To H_{2}(\hC;\mathds{Z})\To H_{2}(\hC;\mathds{Z})\To
H_{2}(\hC;\mathds{Z}/p\mathds{Z})\overset{\partial}{\To}
H_{1}(\hC;\mathds{Z})\To\cdots
$$
\noindent where the above long exact sequence is associated to the
short exact sequence
$$
0\To\mathds{Z}\overset{-\times
p}{\To}\mathds{Z}\overset{\text{projection}}{\To}\mathds{Z}/p\mathds{Z}\To
0
$$
Thus
$$
a= D\partial^{-1}[\alpha]\bmod p
$$
The short exact sequence
$$
0\To\mathds{Z}\overset{-\times
p}{\To}\mathds{Z}\overset{\text{projection}}{\To}\mathds{Z}/p\mathds{Z}\To
0
$$
\noindent also gives rise to the long exact sequence on cohomology
$$
\cdots \To H^{1}(\hC;\mathds{Z})\To H^{1}(\hC;\mathds{Z})\To
H^{1}(\hC;\mathds{Z}/p\mathds{Z})\overset{\partial^{*}}{\To}
H^{2}(\hC;\mathds{Z})\To\cdots
$$
\noindent where $\partial^{*}$ is the Bockstein homomorphism on
cohomology. We now define the \emph{coloured untying invariant} of
$(K,\rho)$ as
$$
\mathrm{cu}(K,\rho)\ass\ \partial^{*}a\smile a\in
H^{3}(\hC;\mathds{Z}/p\mathds{Z})\cong \mathds{Z}/p\mathds{Z}
$$
This gives a non-trivial $\mathds{Z}/p\mathds{Z}$--valued invariant
of $p$--coloured knots.

\begin{prop}
The coloured untying invariant is invariant under $\pm1$--framed
surgery in $\ker \rho$.
\end{prop}

\begin{proof}
By Poincar\'{e} duality the dual of $\mathrm{cu}(K,\rho)$ is equal
to the algebraic intersection number of $D(a)$ with $\partial
D(a)$.

Let $L$ be a loop in $\ker\rho$. Since both components in the
pre-image of $L$ in $\hC$ vanish as elements of
$H_{1}(\hC;\mathds{Z}/p\mathds{Z})$, these components may be chosen
to be disjoint with $D(a)$. Thus performing surgery by $L$ does not
change $\mathrm{cu}(K,\rho)$.
\end{proof}

The next lemma gives a way to calculate the coloured untying
invariant.

\begin{lem}\label{L:vformula}
Let $v\ass
(v_{1},v_{2},\ldots,v_{2n})^{T}\in\mathds{Z}/p\mathds{Z}^{2n}$ be a
column vector such that $v_{i}=a \frown [p^{-1}\xi_{i}]$ for all
$1\leq i\leq 2n$. Then
$$
\mathrm{cu}(K,\rho)=\frac{2(v^{T}\cdot M\cdot v)}{p}\bmod p
$$
\end{lem}

\begin{proof}
Notice first that
$$
\frac{2(v^{T}\cdot M\cdot v)}{p}\bmod p=\frac{v^{T}\cdot
(M+M^{T})\cdot v}{p}\bmod p
$$
\noindent where $v^{T}\cdot (M+M^{T})\cdot v$ equals the linking
number of $\alpha$ with itself in $\hC$. By definition we have
\begin{multline*}
\tfrac{1}{p}\mathrm{Link}(\alpha,\alpha)\bmod
p=\tfrac{1}{p}\left\langle D\partial^{-1}[\alpha]\bmod p,\
[\alpha]\right\rangle=\\\tfrac{1}{p}D\partial^{-1}[\alpha]\bmod
p\frown [\alpha]= \tfrac{1}{p}(a\frown [\alpha]) = a\frown\partial
D(a)=\\ D\left(a\smile\partial^{*}a\right)=
D(\mathrm{cu}(K,\rho))=\mathrm{cu}(K,\rho)
\end{multline*}
which completes the proof.
\end{proof}

We may now prove the main result of this section, that there are at
least $p$ equivalence classes of $p$--coloured knots modulo surgery
in $\ker\rho$, and that these are represented by connect-sums of
left-hand $(p,2)$--torus knots with a given colouring.

\begin{prop}\label{P:r1neqr2}
For $r_{1}\neq r_{2}\bmod{p}$ the connect-sum of $r_{1}$ left-hand
$(p,2)$--torus knots with a given colouring is not equivalent modulo
surgery in $\ker \rho$ to the connect-sum of $r_{2}$ left-hand
$(p,2)$--torus knots with the same colouring.
\end{prop}

\begin{proof}
The proof is by comparing coloured untying invariants.

For $(K,\rho)$ a $p$--coloured left-hand $(p,2)$--torus knot, $\hC$
is a $(1,p)$ lens space. Its fundamental group and therefore also
its first homology is $\mathds{Z}/p\mathds{Z}$. Thus the
$p$--colouring $\rho$ is an automorphism of $\mathds{Z}/p\mathds{Z}$
and therefore corresponds to a non-trivial cohomology class $a\in
H^{1}(\hC; \mathds{Z}/p\mathds{Z})$. Since $\partial^{*}a$ is also
non-trivial $\mathrm{cu}(K,\rho)$ is a non-trivial element of
$H^{3}(\hC;\mathds{Z}/p\mathds{Z})\cong \mathds{Z}/p\mathds{Z}$.

The coloured untying invariant is additive, since taking the
connect-sum of two knots corresponds to taking the direct sum of
their Seifert matrices, and the characteristic link of the
connect-sum may be taken to be the disjoint union of the
characteristic link of the direct summands. An alternative proof of
additivity is that for $K=K_{1}\hash K_{2}$, the double covering
$\hC$ decomposes as $\hC^{1}\oplus \hC^{2}$, and therefore $D(a)$ may
be taken to be the disjoint union of two homology classes in
$H_{2}(\hC; \mathds{Z}/p\mathds{Z})$.
\end{proof}

\section{Applications}\label{S:appl}

\subsection{Surgery presentation for irregular dihedral
coverings of $S^{3}$ branched over knots}

The original motivation for this paper was a series of discussions
between Andrew Kricker and the author about finding surgery
presentations for irregular branched dihedral coverings of
knots.

\begin{figure}
\centering
\labellist\small
\pinlabel {$L$} at 218 208
\pinlabel {$\vdots$} at 19 175
\endlabellist
\label{Fi:lclpres}
\centerline{\includegraphics[scale=0.31]{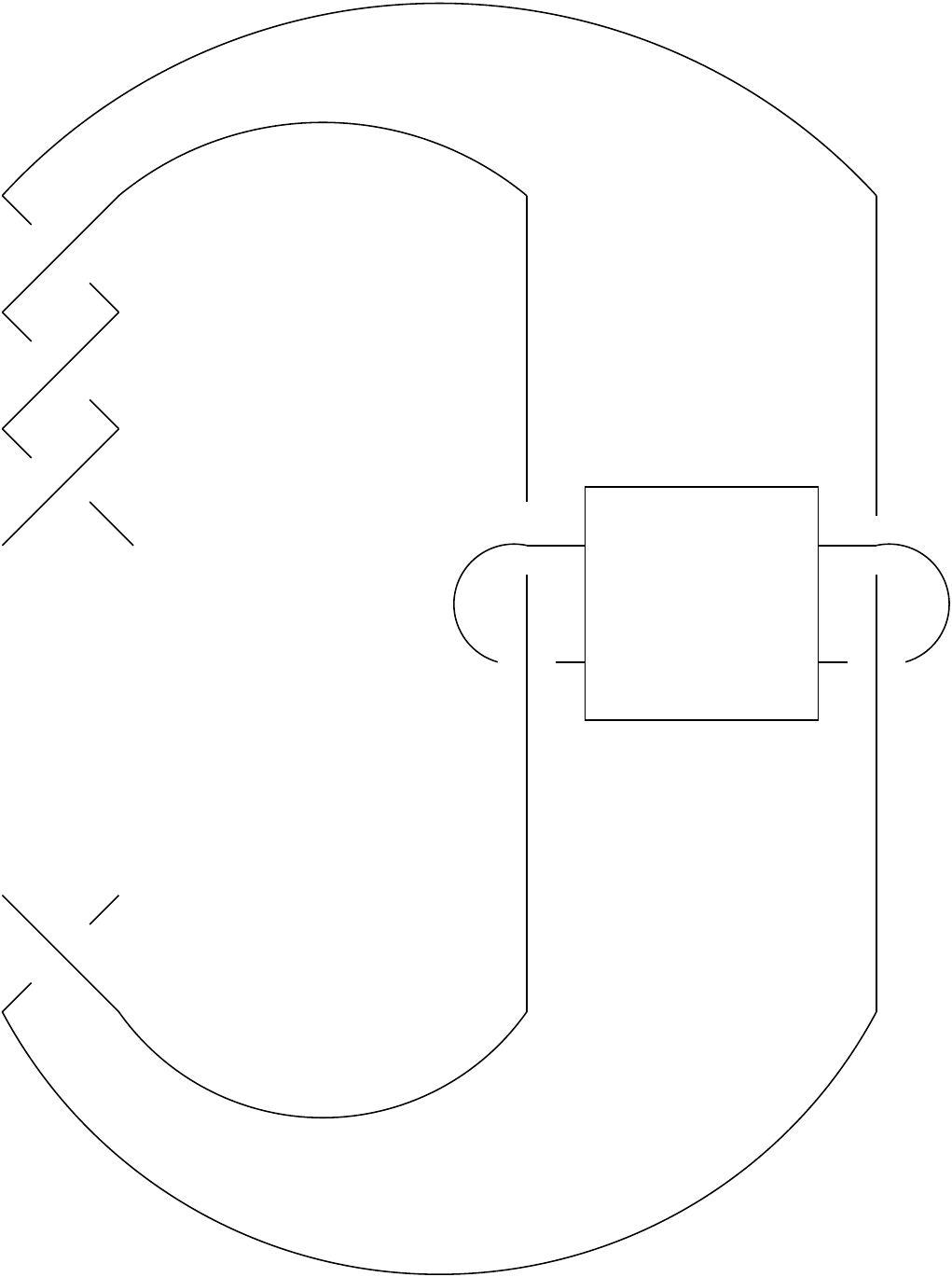}}
\caption{An untying link}
\end{figure}

\begin{figure}
\label{Fi:lklft}
\centering
\labellist\small
\pinlabel {$L$} at 60 36
\pinlabel {\rotatebox{180}{$L$}} at 166 36
\pinlabel {$L$} at 276 36
\pinlabel {$\cdots$} at 383 36
\pinlabel {\rotatebox{180}{$L$}} at 488 36
\pinlabel {$L$} at 600 36
\endlabellist
\centerline{\includegraphics[scale=0.27]{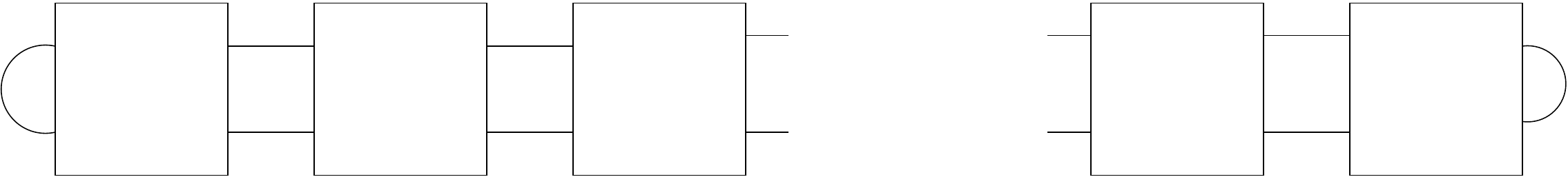}}
\caption{A surgery presentation for the irregular dihedral cover}
\end{figure}

Knowing how to reduce any $p$--coloured knot ($p=3,5$) to a left-hand
$(p,2)$--torus knot (or a connect-sum of such) by surgery in the kernel
of the $p$--colouring gives us a new algorithm for translating between
covering presentations and surgery presentations of closed orientable
$3$--manifolds. First, by $\pm1$--framed surgeries in the kernel of its
colouring, a $p$--coloured knot $K$ may be presented as a framed link
in the complement of a connect-sum of $(p,2)$--torus knots. Lifting the
knot to the irregular branched dihedral covering space lifts also the
link that is in its complement. Explicitly, in \cite{Mos07a} we show
that the surgery presentation of the irregular dihedral cover of a
knot $K$ which is obtained from surgery by $L$ on the $(p,2)$--torus
knot as in \fullref{Fi:lclpres} (or a connect-sum of such) is given
by \fullref{Fi:lklft}.

This method generalizes Yamada's algorithm \cite{Yam02}, and does
not require use of the $3$--move in the $3$--colour case.

\subsection{Surgery presentation of $D_{2p}$--periodic
maps on compact $3$--manifolds}

Montesinos \cite{Mon77} showed that a closed orientable
$3$--manifold is a double branched covering of $S^{3}$ if and only
if this manifold is obtained by rational surgery on a strongly
invertible link $L$ in $S^{3}$. This allows us to `visualize' the
covering involution of such a manifold $M$, since it is conjugate to
the involution of $M$ induced by the involution of $S^{3}$
preserving $L$.

This result has recently been generalized by Przytycki and Sokolov
\cite{PrzS01} and later by Sakuma \cite{Sak01} to all cyclic
branched covering spaces. For a closed orientable $3$--manifold $M$
which admits an orientation-preserving periodic diffeomorphism $f$,
Sakuma showed that $M$ is obtained by integral surgery on a link $L$
in $S^{3}$ which is invariant under a standard $\frac{2\pi}{n}$
rotation $\varphi_{n}$ around a trivial knot, and $f$ is conjugate
to the periodic diffeomorphism of $M$ induced by $\varphi_{n}$.

The key fact used in these papers which prevents the proofs there
from directly carrying over to the dihedral case is that any knot
$K$ can be transformed into an unknot by $\pm 1$--surgeries on a
trivial knot whose linking number with $K$ is $0$, but such
surgeries may not be in the kernel of a $p$--colouring of the knot.
As Makoto Sakuma pointed out to the author, by substituting our main
theorem in this note for this fact, we may generalize the above
result to the case in which the cyclic group is replaced by the
dihedral group $D_{2p}$ with $p\in\set{3,5}$ and in which $M$ is a
regular dihedral $p$--fold covering space. It seems an interesting
problem to generalize this result to a wider class of
$3$--manifolds.

\bibliographystyle{gtart}
\bibliography{link}

\end{document}